\documentclass[10pt]{article}
\usepackage{mathrsfs}
\usepackage{amsmath,amssymb}
\usepackage{amsthm}
\usepackage[colorlinks,linkcolor=red,anchorcolor=blue,citecolor=green]{hyperref}
\usepackage{graphics,graphicx}
\usepackage{subfigure}
\usepackage{color}
\usepackage{enumerate}
 \usepackage[titletoc]{appendix}
 \usepackage{bbm}
\allowdisplaybreaks

\usepackage{tikz}
\usepackage{indentfirst}

\newtheorem{thm}{Theorem}[section]
\newtheorem{lem}[thm]{Lemma}

\newtheorem{remark}[thm]{Remark}

\newtheorem{prob}[thm]{Problem}
\newtheorem{conj}[thm]{Conjecture}

\def\qed{\hfill \rule{4pt}{7pt}}

\def\pf{\noindent {\it{Proof.}\hskip 2pt}}

\def\R{{\mathbb{R}}}
\def\E{{\mathbb{E}}}
\def\N{{\mathbb{N}}}
\def\P{{\mathbb{P}}}

\def\Z{{\mathbb{Z}}}

\def\RW{{\mathrm{RW}}}

\begin{document}
\begin{center}
{{\large\bf Phase transition of disordered random networks on quasi-transitive graphs}\footnote{The project is supported partially by CNNSF (No.~11671216) and by
Hu Xiang Gao Ceng Ci Ren Cai Ju Jiao Gong Cheng-Chuang Xin Ren Cai (No. 2019RS1057).}}
\end{center}

\begin{center}
Liu Yuelin$^{a}$\ \& \ Xiang Kainan$^b$
\vskip 1mm
\footnotesize{$^a$Department of Mathematics, Tianjin University of Finance and Economics\\
            Tianjin City 300222, P. R. China}\\
\footnotesize{$^b$ Hunan Key Laboratory for Computation and Simulation in Science and Engineering \&\\
Key Laboratory of Intelligent Computing and Information Processing of Ministry of Education \&\\
School of Mathematics and Computational Science, Xiangtan University\\
Xiangtan City 411105, Hunan Province, P. R. China}\\
\footnotesize{Emails:  \texttt{liuyuelinmath@qq.com} (Liu)\ \ \ \ \  \texttt{kainan.xiang@xtu.edu.cn} (Xiang)}
\end{center}

\begin{center}
\emph{This paper is dedicated to the memory of Vladas Sidoravicius.}
\end{center}

\begin{abstract}
Given a quasi-transitive infinite graph $G$ with volume growth rate ${\rm gr}(G),$ a transient biased electric network $(G,\,\mathbf{c}_1)$ with bias $\lambda_1\in (0,\,{\rm gr}(G))$ and a recurrent biased one $(G,\,\mathbf{c}_2)$ with bias $\lambda_2\in ({\rm gr}(G),\infty).$ Write $G(p)$ for the Bernoulli-$p$ bond percolation on $G$, and define percolation process $(G(p))_{p\in [0,\, 1]}$ by the grand coupling. Let $(G,\, \mathbf{c}_1,\, \mathbf{c}_2,\, p)$ be the following biased disordered random network: Open edges $e$ in $G(p)$ take the conductance $\mathbf{c}_1(e)$, and closed edges $g$ in $G(p)$ take the conductance $\mathbf{c}_2(g)$. We mainly study recurrence/transience phase transition for $(G,\, \mathbf{c}_1,\, \mathbf{c}_2,\, p)$ when $p$ varies from $0$ to $1$, and our main results are as follows:
 \begin{enumerate}[{\bf (i)}]
 \item On connected quasi-transitive infinite graph $G$ with percolation threshold $p_c\in (0,\, 1),$ the biased disordered random network $(G,\, \mathbf{c}_1,\, \mathbf{c}_2,\, p)$ has a non-trivial recurrence/transience phase transition such that the threshold $p_{c}^{*}\in (0,\, 1)$ is deterministic, and almost surely $(G,\, \mathbf{c}_1,\, \mathbf{c}_2,\, p)$ is recurrent for any $p<p_c^*$ and transient for any $p>p_c^*.$
     On any Cayley graph $G$ of any group which is virtually $\Z$, there is no non-trivial recurrence/transience phase transition for $(G,\, \mathbf{c}_1,\, \mathbf{c}_2,\, p)$, i.e. $p_{c}^{*}=p_c= 1.$
     Note $p_c<1$ for an infinite finitely generated group if and only if it is not virtually $\Z.$ Thus there is a non-trivial recurrence/transience phase transition for $(G,\, \mathbf{c}_1,\, \mathbf{c}_2,\, p)$ with $G$ being a Cayley graph if and only if the corresponding group is not virtually $\Z$.

 \item On $\Z^d$ for any $d\geq 1,$ $p_c^{*}= p_c$ (note $p_c=1$ if and only if $d=1$). And on $d$-regular trees $\mathbb{T}^d$ with $d\geq 3$, $p_c^{*}=(\lambda_1\vee 1) p_c$, and thus $p_c^{*}>p_c$ for any $\lambda_1\in (1,\,{\rm gr}(\mathbb{T}^d)).$ Critical $\left(\Z^d,\,\mathbf{C}_{\lambda_1},\,\mathbf{C}_{\lambda_2},\,p_c\right)$ with $0<\lambda_1<1<\lambda_2$ and $d=2$ or $0<\lambda_1\leq 1<\lambda_2$ and $d\geq 11$ is recurrent almost surely, and so is critical $\left(\mathbb{T}^d,\,\mathbf{C}_{\lambda_1},\,\mathbf{C}_{\lambda_2},\,\frac{1}{d-1}\right)$ with $0<\lambda_1\leq 1<d-1<\lambda_2$ and $d\geq 3.$
     Generally, we propose a conjecture characterizing the $p_c^*.$
  \end{enumerate}
As a contrast, we also consider phase transition of having unique currents or not for $(\Z^d,\, \mathbf{c}_1,\, \mathbf{c}_2,\, p)$ with $d\geq 2$ when $p$ varies from $0$ to $1$ (the case $d=1$ is trivial due to $p_c^*=1$), and prove that
almost surely $(\Z^2,\, \mathbf{c}_1,\, \mathbf{c}_2,\, p)$ with $\lambda_1<1\leq\lambda_2$ has unique currents for any $p\in [0,1]$ (and thus has no current uniqueness/non-uniqueness phase transition), and conjecture that the same conclusion holds for $d\geq 3.$

\vskip 3mm

\vskip 3mm

\end{abstract}

\section{Introduction}
\setcounter{equation}{0}
\noindent
Let $G=(V,\, E)$ be a locally finite infinite connected graph with vertex set $V$ and edge set $E$, and fixed root $o\in V.$ When two vertices $x$ and $y$ of $G$ are adjacent, write $x\sim y$ and denote by $\{x,y\}$ (resp. $xy$) the corresponding undirected edge (resp. directed edge from $x$ to $y$). Write $\Z$ (resp. $\N$) for the set of all integers (resp. natural numbers). Let each $B_n(o)$ be the closed ball in $G$ centered at $o$ with radius $n,$ and $\vert A\vert$ the cardinality of a set $A.$ Define the lower (volume) growth rate and the (volume) growth rate of graph $G$ respectively by
\begin{eqnarray}\label{eq-growth-rate}
\underline{{\rm gr}}(G)=\liminf\limits_{n\rightarrow\infty}\sqrt[n]{\vert B_n(o)\vert}\ \mbox{and}\
{\rm gr}(G)=\lim\limits_{n\rightarrow\infty}\sqrt[n]{\vert B_n(o)\vert}\ (\mbox{if exists}).
\end{eqnarray}
Call an edge weighted graph $(G,\,\mathbf{c})$ is a network (or an electrical network), and $\mathbf{c}:\ E\rightarrow\R_{+}=[0,\infty)$ the conductance function and its reciprocal $\mathbf{r}=1/\mathbf{c}$ the resistance function.  Recall the random walk associated to a network $(G,\,\mathbf{c})$ is a random walk $(X_n)_{n\geq 0}$ on graph $G=(V,\,E)$ with transition probability $\mathbf{p}(\cdot, \cdot)$ such that
\[
\mathbf{p}(x,y):= \frac{c(\{x,y\})}{\sum\limits_{x\in e} \mathbf{c}(e)},\ x,y\in V,\ x\sim y.
\]
Say network $(G,\,\mathbf{c})$ is transient (resp. recurrent) if so is its associated random walk $(X_n)_{n\geq 0}$. For any $\lambda\in (0,\,\infty),$ let
$$\mathbf{C}_\lambda(e)=\lambda^{-\vert e\vert},\ e=\{x,y\}\in E,$$
where with ${\rm dist}(\cdot,\cdot)$ being the graph distance on $G$,
$$\vert e\vert={\rm dist}(x,o)\wedge {\rm dist}(y,o)={\rm dist}(e,o).$$
Say $(G,\,\mathbf{C}_\lambda)$ is a biased network with bias $\lambda$ and its associated random walk ${\rm RW}_\lambda$ on $G$ a biased random walk with bias $\lambda$.

Suppose $(G,\, \mathbf{c}_1)$ and $(G,\, \mathbf{c}_2)$ are two electrical networks. Introduce Bernoulli bond percolation process $\omega=(\omega_p)_{p\in [0,\,1]}:=(G(p))_{p\in [0,\,1]}$ on $G$ by the grand coupling: Let $(U_e)_{e\in E}$ be an i.i.d.\,family of the uniform distribution on $[0,\,1].$ An edge $e$ is open in Bernoulli-$p$ bond percolation $G(p)$ if $U_e\leq p$ and closed in $G(p)$ otherwise. Namely, for any $p\in [0,\,1],$ $\omega_p(e)=I_{\{U_e\leq p\}},\ e\in E.$ Write $\mathbb{P}_p$ for the law of $\omega_p$. Define the following disordered random network process $\left((G,\,\mathbf{c}_1,\,\mathbf{c}_2,\,p)\right)_{p\in [0,1]}$ on $G$: For any $p\in [0,\,1],$ $(G,\,\mathbf{c}_1,\,\mathbf{c}_2,\,p)$
is a random network such that each open edge $e$ in $G(p)$ takes the conductance $\mathbf{c}_1(e)$, while each closed edge $g$ in $G(p)$
takes the conductance $\mathbf{c}_2(g)$. Specially for any $0<\lambda_1<\lambda_2<\infty,$ each $(G,\,\mathbf{C}_{\lambda_1},\,\mathbf{C}_{\lambda_2},\,p)$ is called a biased disordered random network (with biases $\lambda_1$ and $\lambda_2$). When $(G,\, \mathbf{c}_1)$ is transient and $(G,\, \mathbf{c}_2)$ is recurrent, each $(G,\, \mathbf{c}_1,\, \mathbf{c}_2,\, p)$ with $p\in (0,1)$ is called a competing disordered random network in the sense that $(G,\, \mathbf{c}_1)$ wins $(G,\, \mathbf{c}_2)$ if $(G,\, \mathbf{c}_1,\, \mathbf{c}_2,\, p)$ is transient, and otherwise $(G,\, \mathbf{c}_2)$ wins $(G,\, \mathbf{c}_1).$

This paper mainly studies recurrence/transience phase transition for competing disordered random networks $(G,\, \mathbf{c}_1,\, \mathbf{c}_2,\, p)$ when $p$ varies from $0$ to $1$. To state our main results, recall the following preliminaries:
\begin{enumerate}[{\bf (i)}]
\item

Let ${\rm Aut}(G)$ be the group consisting of all automorphisms of graph $G$. $G$ is quasi-transitive (resp. transitive) if
there are only finitely many orbits (resp. is only one orbit) under group action of ${\rm Aut}(G)$. When $G$ is quasi-transitive, $\underline{{\rm gr}}(G)={\rm gr}(G),$
and the critical parameter $\lambda_c(G)$ such that ${\rm RW}_\lambda$ is transient for $\lambda<\lambda_c(G)$ and recurrent for $\lambda>\lambda_c(G)$ is just ${\rm gr}(G)$ (which can be proved similarly to \cite[Theorem 1.1]{RL1995}).

\item A group $\Gamma$ is an extension of a group $H$ by $Q$ if there is a short exact sequence
\[
1 \longrightarrow Q \stackrel{f}\longrightarrow \Gamma\stackrel{g}\longrightarrow H\longrightarrow 1
\]
such that $f,g$ are group homomorphisms and ${\rm Im}(f)={\rm Ker}(g)$, equivalently $Q$ is a normal subgroup of $\Gamma$ and $H$ is isomorphic to quotient group $\Gamma/Q$. If $Q$ is a finite group, $\Gamma$ is called a finite extension of $H$ or is virtually $H.$ In other words, $\Gamma$ is a finite extension of $H$ if $H$ is a subgroup of $\Gamma$ with a finite index $[\Gamma:\,H].$
%
\end{enumerate}
Then our main results, Theorems \ref{generalgraph01}, \ref{recurthm}, \ref{generald}-\ref{thm-current-unique}, are summarized as follows:
 \begin{enumerate}[{\bf (i)}]
 \item On connected quasi-transitive infinite graph $G$ with percolation threshold $p_c=p_c(G)\in (0,\, 1),$ for any $0<\lambda_1<\lambda_c(G)<\lambda_2,$ $(G,\, \mathbf{C}_{\lambda_1},\, \mathbf{C}_{\lambda_2},\, p)$ has a non-trivial recurrence/transience phase transition such that the threshold $p_{c}^{*}\in (0,\, 1)$ is deterministic, and almost surely $(G,\, \mathbf{C}_{\lambda_1},\, \mathbf{C}_{\lambda_2},\, p)$ is recurrent for any $p<p_c^*$ and transient for any $p>p_c^*.$ On any Cayley graph $G$ of any group which is virtually $\Z$, there is no non-trivial recurrence/transience phase transition for $(G,\, \mathbf{C}_{\lambda_1},\, \mathbf{C}_{\lambda_2},\, p)$, i.e. $p_{c}^{*}=p_c= 1.$
     Note $p_c<1$ for an infinite finitely generated group if and only if it is not virtually $\Z.$ Thus there is a non-trivial recurrence/transience phase transition for $(G,\, \mathbf{C}_{\lambda_1},\, \mathbf{C}_{\lambda_2},\, p)$ with $G$ being a Cayley graph if and only if the corresponding group is not virtually $\Z$.

 \item On $\Z^d$ for any $d\geq 1,$ $p_c^{*}= p_c$ (note $p_c=1$ if and only if $d=1$). And on $d$-regular trees $\mathbb{T}^d$ with $d\geq 3$, $p_c^{*}=(\lambda_1\vee 1) p_c$, and thus $p_c^{*}>p_c$ for any $\lambda_1\in (1,\,\lambda_c(\mathbb{T}^d)).$ Generally, we propose Conjecture \ref{conj-p_c^*} to characterize the $p_c^*$.  Critical $\left(\Z^d,\,\mathbf{C}_{\lambda_1},\,\mathbf{C}_{\lambda_2},\,p_c\right)$ with $0<\lambda_1<1<\lambda_2$ and $d=2$ or $0<\lambda_1\leq 1<\lambda_2$ and $d\geq 11$ is recurrent almost surely, and so is critical $\left(\mathbb{T}^d,\,\mathbf{C}_{\lambda_1},\,\mathbf{C}_{\lambda_2},\,\frac{1}{d-1}\right)$ with $0<\lambda_1\leq 1<d-1<\lambda_2$ and $d\geq 3$; and moreover we have Conjectures
     \ref{conj-critical-recurrent/transient} and \ref{conj-critical-recurrent-tree}.

 \item As a contrast, for $(\Z^d,\,\mathbf{C}_{\lambda_1},\, \mathbf{C}_{\lambda_2},\, p)$ with $d\geq 2$ (the case $d=1$ is trivial due to $p_c^*=1$), we also consider phase transition of having unique currents or not when $p$ varies from $0$ to $1$, and prove that almost surely $(\Z^2,\,\mathbf{C}_{\lambda_1},\, \mathbf{C}_{\lambda_2},\, p)$ with $\lambda_1<1\leq\lambda_2$ has unique currents for any $p\in [0,1]$ (no current uniqueness/non-uniqueness phase transition!), and think that the same conclusion holds for $d\geq 3$ (Conjecture \ref{conj-current-unique}).
\end{enumerate}

Now we are in the position to describe backgrounds, motivations and interests on disordered random networks and our main results.

Disordered random network is one of the most important models in discrete probability theory and is also widely applied in physics and biology. The natural physical background of disordered random network is to study the effective conductance in the doped semiconductors where each edge has different resistance which decays along temperature (\cite{MA1960, AHL1971, POF2013}). Percolation theory plays an important role in analysis of above models. In biology, disordered random network can be seen in several statistical biology models such as DNA-unzipping experiments or DNA-polymerase phenomenon (\cite{BBCEM2006, BBCEM2007, HFR2009, KSJW2002}). Recall from \cite[pp.\,6-7, pp.\,380-382]{GG1999}, disordered random network $(G,\,\mathbf{c}_1,\,\mathbf{c}_2,\,p)$ on finite (and infinite) graphs $G$ is a mathematical modelling of a disordered mixture of two conductor materials $A$ and $B;$ and effective resistance $\mathscr{R}_i$ of
disordered random network $\left(\{0,\,1,\,\ldots,\,i\}^d,\, 1,\, 0,\, p\right)$ between the bottom and top sides of $\{0,\,1,\,\ldots,\,i\}^d$ satisfies that for a constant $p_c(d)\in (0,\, 1),$
$$
\mathscr{R}_i= \infty\ \mbox{a.s.\,for all large $i$}\ \mbox{if}\ p< p_c(d)\ \mbox{and}\ \mbox{a.s.}\,
\lim\limits_{i\rightarrow \infty}\frac{\mathscr{R}_{i}}{i^{2-d}}\in (0,\,\infty)\ \mbox{exists if}\ p>p_c(d)\ (\cite{RK1983,VJ1994}).
$$

Theoretically, Chernov \cite{AAC1967} introduced the idea of random walk generated in a random environment in 1967 as a mathematical model to study the transport in a random media in biology. Random walk in random environment (RWRE) has become one of the most popular probability models in recent decades (\cite{AAC1967, FS1975, YS1982, RY1995, SS2004, AS2004, OZ2006}). Typical RWRE on Euclidean lattice $\Z^d$ can be defined as follows: Suppose $\omega=\{p_x\}_{x\in\Z^d}$ is an i.i.d.\,family of random probability measures on $\mathscr{S}_d=\{\pm e_i,\, 1\leq i\leq d\},$
where each $e_i$ is the $i$th standard unit vector in $\Z^d.$ Given random environment $\omega,$ define a nearest-neighbour random walk $(X_n)_{n=0}^{\infty}$ on $\Z^d$ by
$$\P_\omega[X_{n+1}=x+y\,\vert\, X_n=x]=p_x(y),\ y\in\mathscr{S}_d,\ n\in\Z_{+}:=\{0,\,1,\,2,\,\ldots\}.$$
Such a model was firstly defined on $\Z$ by Solomon 1975 \cite{FS1975} (in this case it is also known as Sinai's simple random walk in random environment \cite{YS1982}). The above definition can be extended to more general (random) graphs with random environment being ergodic (e.g. translation invariant independent random environment on transitive graphs). See \cite{AS2004, OZ2006} and \cite[pp.\,56-57]{LP2017}.
There are two layers of randomness for RWRE which makes the model very interesting: the first is the random environment; the second is the random walk in a given random environment. As written in \cite[p.\,56]{LP2017}, ``The topic of RWRE with any i.i.d.\,transition probabilities, is quite natural and extensive but only partially understood, except on trees." Secondly, there is a class of interesting random walks in inhomogeneous random environment such that the law of the random environment is stationary and ergodic with respect to space-time shifts (see \cite{MVK1985}, \cite{BR2017} and \cite{OZ2006}). Thirdly, a special class of RWRE models, which are reversible Markov chain in the environment, is given by the class of nearest-neighbour random conductance models (RCMs). In RCM, random environment is generally translation-invariant and usually random conductance function on edges is an i.i.d.\,family. See \cite{MTB2004, SS2004, BB2007, BP2007, MP2007, BD2010, ABDH2013} and particularly survey \cite{MB2011}.

For fixed $p\in (0,\,1),$ disordered random network $(G,\,\mathbf{c}_1,\,\mathbf{c}_2,\,p)$ is a random electrical network (hence a random conductance model and a random resistance one) and can be viewed naturally as a RWRE; and different with usual RWRE, RCM models, we remove the assumption of stationarity, ergodicity and translation-invariance of random environments (i.e. random conductances) which calls for new techniques and more precise estimation on percolation structure to study a disordered random network. In a certain sense, $(G,\,\mathbf{c}_1,\,\mathbf{c}_2,\,p)$ is a new type of RWRE models. When $p$ evolves from $0$ to $1$, disordered random networks $(G,\,\mathbf{c}_1,\,\mathbf{c}_2,\,p)$ is an interesting interpolation of two deterministic networks $(G,\,\mathbf{c}_1)$ and $(G,\,\mathbf{c}_2);$ and to understand typical probability behaviours varying in $p$ of the interpolation between two networks (or their associated random walks), is an original motivation to study disordered random network process $((G,\,\mathbf{c}_1,\,\mathbf{c}_2,\,p))_{p\in [0,\,1]}.$ As said before, $((G,\,\mathbf{c}_1,\,\mathbf{c}_2,\,p))_{p\in [0,\,1]}$ is also
a competing stochastic process (in fact, a new competing stochastic process) when $(G,\, \mathbf{c}_1)$ is transient and $(G,\, \mathbf{c}_2)$ is recurrent in the sense that $(G,\, \mathbf{c}_1)$ wins $(G,\, \mathbf{c}_2)$ if $(G,\, \mathbf{c}_1,\, \mathbf{c}_2,\, p)$ is transient, and otherwise $(G,\, \mathbf{c}_2)$ wins $(G,\, \mathbf{c}_1).$ The interpolation and competition lead to a natural featured topic of disordered random networks: recurrence/transience phase transitions for $((G,\,\mathbf{c}_1,\,\mathbf{c}_2,\,p))_{p\in [0,\,1]}$ with $p$ varying from $0$ to $1.$
For other phase transitions related to disordered random networks, see Problem \ref{prob-PhaseTransition}. We hope that the percolation theory can lead to a sequence of profound and interesting results for disordered random networks, and conversely disordered random networks can provide new interesting topics (insights) to the percolation theory.

Note interpolation and competition are interesting topics for stochastic processes. Recall that \cite{WLE2018} introduced a $p$-rotor walk on $\Z$ which is an interpolation between simple random walk and deterministic rotor walk, and proved an invariance principle such that the limiting process is a doubly perturbed Brownian motion multiplying constant $\sqrt{\frac{1-p}{p}}$. Here the interpolation is to choose random transition probability through site percolation (refer to Subsection \ref{sec-disorderedRW} for such a similar interpolation). Additionally, there are some models studying the competing behaviour such as competing frogs model \cite{MTF2019}, and competing first passage percolation (Richardson model) \cite{DR1973, HP1998, HP2000} and so on.

Our aforementioned main results show that for disordered random networks, recurrence vs transience phase transition and current uniqueness vs non-uniqueness one may present different phase transition vs no phase transition phenomena, and there are interesting universal properties; and disordered random networks can provide new interesting topics to the percolation theory (for this viewpoint see also Section \ref{sec-concluding}). To study systematically disordered random networks is our future goal.

Finally we need to explain the reason for choosing biased conductances $\mathbf{C}_{\lambda}$ to study disordered random networks.
Recall an original motivation for introducing ${\rm RW}_\lambda$ on graphs $G$ by Berretti and Sokal \cite{BA-SA-1985} in 1985 is to design a new Monte Carlo algorithm for self-avoiding walks, see \cite{LG-SA1988,SA-JM1989, RD1994} for refinements of this idea. And ${\rm RW}_\lambda$ has received much attention recently, see \cite{BG-FA2016}, \cite{LP2017} and references therein.
When $G$ is a locally finite quasi-transitive infinite graph, ${\rm RW}_\lambda$s capture geometric information on $G$: notably critical parameter $\lambda_c(G)$, such that ${\rm RW}_\lambda$ is transient for $\lambda<\lambda_c(G)$ and recurrent for $\lambda>\lambda_c(G)$, is just the volume growth rate ${\rm gr}(G)$ for $G.$ While growth of groups is an important area for group theory (\cite{BE2014, HH2015}). Secondly when $G$ is a random graph (e.g. Galton-Watson tree), ${\rm RW}_\lambda$ has close relation with trapping phenomenon of RWREs (\cite{BG-FA2016}). Thirdly, networks $(G,\,\mathbf{C}_\lambda)$ ($\lambda\not=1$) are not transitive, and may provide a very useful setting to check some properties for probability models in a non-ergodic situation. All these facts will make geometry of percolation play an important role in studying recurrence/transience phase transition of disordered random network process $((G,\,\mathbf{C}_{\lambda_1},\,\mathbf{C}_{\lambda_2},\,p))_{p\in [0,1]}$
and can lead to some interesting results of the mentioned phase transition.\\

\noindent{\bf Notations.} For any graph $G$, use $V=V(G)$ and $E=E(G)$ to denote its vertex and edge sets respectively, and let $\overrightarrow{E}$ be the set of all directed edges of graph $G.$ For any $e\in \overrightarrow{E},$ write $e_{-}$ and $e_{+}$ for its tail and head respectively. Note Bernoulli bond percolation process $\omega=(\omega_p)_{p\in [0,\,1]}:=(G(p))_{p\in [0,\,1]}$ on $G$ is defined by the grand coupling,
and $p_c=p_c(G)$ is the corresponding percolation threshold for infinite $G$. And when $G$ is infinite and quasi-transitive, $\lambda_c(G)={\rm gr}(G).$

For two nonnegative functions $f$ and $g$ defined on a set, denote
$f\asymp g$ if for two positive constants $c_1$ and $c_2$, $c_1g\leq f\leq c_2g$; and denote $f(x)\asymp g(x)$ as $x\rightarrow x_0$ if
for two positive constants $c_1$ and $c_2,$ $c_1g(x)\leq f(x)\leq c_2g(x)$ for $x$ sufficiently close to $x_0.$ Recall $\Z$ (resp. $\N$) is the set of all integers (resp. natural numbers),
and $\Z_+=\{0,\,1,\,2,\ldots\}.$

\section{Main results}
\setcounter{equation}{0}
\noindent

\begin{thm}\label{generalgraph01}
On connected quasi-transitive locally finite infinite graph $G$ with percolation threshold $p_c\in (0,\, 1),$ for any $0<\lambda_1<\lambda_c(G)<\lambda_2,$ $((G,\, \mathbf{C}_{\lambda_1},\, \mathbf{C}_{\lambda_2},\, p))_{p\in [0,\,1]}$ has a non-trivial recurrence/transience phase transition such that the threshold $p_{c}^{*}\in (0,\, 1)$ is deterministic, and almost surely $(G,\, \mathbf{C}_{\lambda_1},\, \mathbf{C}_{\lambda_2},\, p)$ is recurrent for any $p<p_c^*$ and transient for any $p>p_c^*.$ On any Cayley graph $G$ of any group which is virtually $\Z$, there is no non-trivial recurrence/transience phase transition for $((G,\, \mathbf{C}_{\lambda_1},\, \mathbf{C}_{\lambda_2},\, p))_{p\in [0,\,1]}$, i.e. $p_{c}^{*}=p_c= 1.$
%
%
%
%
%
%
%
\end{thm}
\vskip 2mm

\begin{remark}\label{remark-PT}
{\bf (i)} Note $p_c<1$ holds for an infinite finitely generated group if and only if it is not virtually $\Z$ (\cite[Theorem~1.3]{DGRSY2018}). Thus there is a non-trivial recurrence/transience phase transition for $((G,\, \mathbf{C}_{\lambda_1},\, \mathbf{C}_{\lambda_2},\, p))_{p\in [0,\,1]}$ with $G$ being a Cayley graph if and only if the corresponding group is not virtually $\Z$.
Additionally, when $\lambda_1,\lambda_2< \lambda_c(G)$ (resp. $\lambda_1,\lambda_2>\lambda_c(G)$), from the Rayleigh's monotonicity principle, almost surely every $(G,\, \mathbf{C}_{\lambda_1},\, \mathbf{C}_{\lambda_2},\, p)$ is transient (resp. recurrent).


{\bf (ii)} There are two thresholds, $p_c^*$ and $\widehat{p}_c^*,$ for recurrence/transient phase transition of
disordered random networks $((G,\, \mathbf{C}_{\lambda_1},\, \mathbf{C}_{\lambda_2},\, p))_{p\in [0,\,1]}$ on quasi-transitive infinite graph $G:$
\begin{eqnarray*}
&&p_c^*=\sup\left\{p\in [0,\,1]:\ \mbox{almost surely},\ (G,\, \mathbf{C}_{\lambda_1},\, \mathbf{C}_{\lambda_2},\, q)\ \mbox{is recurrent
            for all}\ q\in [0,\,p)\right\},\\
&&\widehat{p}_c^*=\sup\left\{p\in [0,\,1]:\ (G,\, \mathbf{C}_{\lambda_1},\, \mathbf{C}_{\lambda_2},\, p)\ \mbox{is almost surely recurrent}\right\}.
\end{eqnarray*}
Due to ergodicity of Bernoulli bond percolation on $G$ (\cite[Proposition~7.3]{LP2017}) and the Rayleigh's monotonicity principle,
it is easy to check that $p_c^*=\widehat{p}_c^*.$ See {\bf (i)} in proving Theorem \ref{generalgraph01}.

{\bf (iii)} To prove recurrence/transience of random networks, a known approach is to take the average network (see \cite[Exercises~2.96-2.97]{LP2017}): Suppose $R$ (resp. $C$) is a random resistance (resp. conductance) function on graph $G$ such that for any $e\in E,$
$$r(e)=\mathbb{E}[R(e)]\in [0,\infty)\ (\mbox{resp.}~c(e)=\mathbb{E}[C(e)]\in [0,\infty)).$$
If $(G,\,r)$ is transient (resp. $(G,\,c)$ is recurrent), then $(G,\,R)$ is a.s. transient (resp. $(G,\,C)$ is a.s. recurrent).
The method lose effect for competing disordered random networks $(G,\, \mathbf{c}_1,\, \mathbf{c}_2,\, p)$ with $p\in (0,\,1):$ When taking average for random conductance function $C(\cdot)$, we have $(G,\,c)$ is a transient network; while when taking average for random resistance function $R(\cdot)$, we see $(G,\, r)$ is a recurrent network.

To get the phase transition of recurrence/transience and value of critical parameter $p_c^{*}$, we need to use one or more ingredients such as the Nash-Williams criterion, the Rayleigh's monotonicity principle, the energy transience/recurrence criterion, geodesic spanning tree, rough embeddings, and more delicate properties with respect to structure of percolations in the cases of Theorems \ref{generalgraph01}, \ref{recurthm}, \ref{generald} and \ref{regulartree1}.

%

\end{remark}
\vskip 2mm

\begin{thm}\label{recurthm}
{\bf (i)} Given any two networks $(\Z,\,\mathbf{c}_1)$ and $(\Z,\,\mathbf{c}_2)$ such that $(\Z,\,\mathbf{c}_2)$ is recurrent. Then almost surely, all disordered random networks $(\Z,\,\mathbf{c}_1,\,\mathbf{c}_2,\, p)$ with $p\in [0,\,1)$ are recurrent.

{\bf (ii)} Let $G$ be a Cayley graph of $\Z$, and $(G,\,\mathbf{c}_1)$ and $(G,\,\mathbf{c}_2)$ two networks with
$c=\sup\limits_{e\in E}\{\mathbf{c}_2(e)\}<\infty.$ Then almost surely, all $(G,\,\mathbf{c}_1,\,\mathbf{c}_2,\, p)$ with $p\in [0,\,1)$ are recurrent.

{\bf (iii)} Consider graph $\Z\times G$ with $G$ being a finite connected graph, two connected networks $(\Z\times G,\, \mathbf{c}_1)$ and $(\Z\times G,\, \mathbf{c}_2)$ such that $\mathbf{c}_1$ and $\mathbf{c}_2$ are positive functions, and network on $\Z$ with the conductance function
$$\mathbf{c}_2^\prime(\{x,x+1\})=\max\limits_{y\in G}\left\{\mathbf{c}_2\left(\{(x,y), (x+1,y)\}\right)\right\},\ x\in\Z$$
is recurrent. Then almost surely, for any $p\in [0,\, 1),$ $(\Z\times G,\,\mathbf{c}_1,\,\mathbf{c}_2,\,p)$ is recurrent.

{\bf (iv)} Let $\Gamma$ be a finite extension of group $\Z$ and $G$ a Cayley graph of $\Gamma.$ Assume $(G,\,\mathbf{c}_1)$ and $(G,\,\mathbf{c}_2)$ are two connected networks such that $\mathbf{c}_1$ and $\mathbf{c}_2$ are positive functions, and $c=\sup\limits_{e\in E}\{\mathbf{c}_2(e)\}<\infty.$ Then almost surely, all $(G,\,\mathbf{c}_1,\,\mathbf{c}_2,\,p)$ with $p\in [0,1)$ are recurrent.
\end{thm}
\vskip 2mm

\begin{prob}\label{prob-Z-Cayley-Graph}
Do there exist a Cayley graph $\widehat{\Z}$ of $\Z$, a transient network $\left(\widehat{\Z},\,\mathbf{c}_1\right)$ and a recurrent one $\left(\widehat{\Z},\,\mathbf{c}_2\right)$ such that $\left(\left(\widehat{\Z},\,\mathbf{c}_{1},\,\mathbf{c}_{2},\, p\right)\right)_{p\in [0,\,1]}$ has a nontrivial recurrence/transience phase transition?
\end{prob}
\vskip 2mm

When $p_c^*\in (0,1),$ does $p_c^*=p_c$ hold? We will see in Theorem \ref{generald}\ (i) and Theorem \ref{regulartree1}\,(i) that both $p_c^*=p_c$ (on $\Z^d,\ d\geq 2$) and $p_c^*>p_c$ (on $\mathbb{T}^d,\ d\geq 3$) may be true; and generally we propose Conjecture \ref{conj-p_c^*} to characterize $p_c^*.$ Additionally, whether critical $(G,\,\mathbf{C}_{\lambda_1},\,\mathbf{C}_{\lambda_2},\,p_c^*)$ on $\Z^d$ and $\mathbb{T}^d$ is recurrent almost surely or not, we have Theorem \ref{generald}\ (ii) and Theorem \ref{regulartree1}\,(ii), and Conjectures \ref{conj-critical-recurrent/transient} and \ref{conj-critical-recurrent-tree}.

\begin{thm}\label{generald}
\begin{enumerate}[{\bf (i)}]
\item For $\left(\left(\Z^d,\,\mathbf{C}_{\lambda_1},\,\mathbf{C}_{\lambda_2},\,p\right)\right)_{p\in [0,\,1]}$ with $0<\lambda_1\leq 1<\lambda_2$ and $d\geq 3$ or $0<\lambda_1<1\leq \lambda_2$ and $d=2$, $p_c^*$ is just $p_c.$

\item
    Critical $\left(\Z^d,\,\mathbf{C}_{\lambda_1},\,\mathbf{C}_{\lambda_2},\,p_c\right)$ with $0<\lambda_1<1<\lambda_2$ and $d=2$ or $0<\lambda_1\leq 1<\lambda_2$ and $d\geq 11$ is recurrent almost surely.

\end{enumerate}
\end{thm}
\vskip 2mm


For $\mathbb{T}^d$ with $d\geq 3,$ $p_c=\frac{1}{d-1},$ $\lambda_c(\mathbb{T}^d)={\rm gr}(\mathbb{T}^d)= d-1$.

\begin{thm}\label{regulartree1}
\begin{enumerate}[{\bf (i)}]
\item
For $\left(\left(\mathbb{T}^d,\,\mathbf{C}_{\lambda_1},\,\mathbf{C}_{\lambda_2},\, p\right)\right)_{p\in [0,\,1]}$ with $0<\lambda_1<d-1\leq \lambda_2$ and $d\geq 3,$
$$p_c^*=\frac{(\lambda_1\vee 1)}{d-1}=(\lambda_1\vee 1)p_c.$$
\item
 Critical $\left(\mathbb{T}^d,\,\mathbf{C}_{\lambda_1},\,\mathbf{C}_{\lambda_2},\,\frac{1}{d-1}\right)$ with $0<\lambda_1\leq 1<d-1<\lambda_2$ and $d\geq 3$ is recurrent almost surely.

 \end{enumerate}
\end{thm}
\vskip 2mm

For any current $i$ on network $(G,\,\mathbf{c})$, define
 \[
 {\rm d}^*i(x)= \sum\limits_{y\sim x} i(xy), \ x\in V.
 \]
Say currents are unique on $(G,\,\mathbf{c})$ if for any currents $i, i'$ satisfying $d^*i= d^*i'$, we have $i= i'$.
As a contrast to recurrence/transience phase transition, we also consider phase transition of having unique currents or not for $(\Z^d,\, \mathbf{C}_{\lambda_1},\, \mathbf{C}_{\lambda_2},\, p)$ with $d\geq 2$ when $p$ varies from $0$ to $1$ (the case $d=1$ is trivial due to $p_c^*=1$), and prove that
almost surely $(\Z^2,\, \mathbf{C}_{\lambda_1},\, \mathbf{C}_{\lambda_2},\, p)$ with $\lambda_1<1\leq\lambda_2$ has no current uniqueness/non-uniqueness phase transition, and think that the same conclusion holds for $d\geq 3$ (Conjecture \ref{conj-current-unique}).
Note disordered random walk in Subsection \ref{sec-disorderedRW} can have different features from those of random networks $(G,\,\mathbf{C}_{\lambda_1},\,\mathbf{C}_{\lambda_2},\,p)$ (e.g. Theorem \ref{thm-sz1}).

\begin{thm}\label{thm-current-unique}
Almost surely, for any $p\in [0,\,1],$ $(\Z^2,\, \mathbf{C}_{\lambda_1},\, \mathbf{C}_{\lambda_2},\, p)$ with $0<\lambda_1<\lambda_2<\infty$ has unique currents.
\end{thm}

\section{Proofs of main results}
\setcounter{equation}{0}
\noindent
In this section, we firstly introduce some necessary preliminaries in Subsection \ref{pre}, then prove Theorem \ref{generalgraph01}, Theorems \ref{recurthm}, \ref{generald} and \ref{thm-current-unique},
and Theorem \ref{regulartree1} in respectively Subsection \ref{quasi}, Subsection \ref{rzd}, and Subsection \ref{rrt}.

\subsection{Preliminaries}\label{pre}
\noindent
Given two subsets $A$ and $Z$ of $V,$ call $v:\ V\rightarrow\mathbb{R}$ is a voltage function if it is harmonic at any $x\notin A\cup Z$.
Call a function $\theta:\ \overrightarrow{E}\rightarrow\mathbb{R}$ is a flow between $A$ and $Z$ if
$$\theta(xy)=-\theta(yx),\ \forall\,x,y\in V,\ x\sim y,\ \mbox{and}\ \sum_{w\sim z}\theta(zw)= 0,\ \forall z\notin A\cup Z.$$
For an antisymmetric function $\theta$ on $\overrightarrow{E}$, define its energy to be $\mathscr{E}(\theta)=\frac{1}{2}\sum_{xy\in \overrightarrow{E}} \theta^2(xy)\mathbf{r}(\{x,y\})$.
Call a flow $i(\cdot)$ is a current on network $(G,\,\mathbf{c})$ between $A$ and $Z$ if there is a voltage function $v$ satisfying the Ohm's law:
$$\mbox{for any}\ x\sim y,\ v(x)- v(y) = i(xy)/\mathbf{c}(\{x,y\}):=i(xy)\mathbf{r}(\{x,y\}).$$
Since $\sum_{x\sim a} i(ax)$ is the total amount of current flowing into the circuit at vertex $a$, one can regard the entire circuit between $a$ and $Z$ as a single conductor with effective conductance
\[C_{\rm eff}:= \sum\limits_{x\sim a} \mathbf{c}(\{a, x\})\P\left[a\rightarrow Z\right] =: \mathscr{C}(a\leftrightarrow Z)= \mathscr{C}_{\mathbf{c}}(a\leftrightarrow Z)= \mathscr{C}_{\mathbf{c}}(a\leftrightarrow Z;\,G),\]
where $\P(a\rightarrow Z)$ is the probability that $(X_n)_{n\geq 0}$, the random walk associated to network $(G,\,\mathbf{c})$ starting at $a$, hits $Z$ before visiting $a$ again. Define the effective resistance between $a$ and $Z$ as
$$\mathscr{R}(a\leftrightarrow Z)=\mathscr{R}(a\leftrightarrow Z;\,G)=\frac{1}{\mathscr{C}(a\leftrightarrow Z)}.$$
When $A$ is not a singleton, define $\mathscr{C}(A\leftrightarrow Z)$ to be $\mathscr{C}(a\leftrightarrow Z)$ by identifying $A$ to a single vertex $a,$ and $\mathscr{R}(A\leftrightarrow Z)=\frac{1}{\mathscr{C}(A\leftrightarrow Z)}.$

To define $\mathscr{C}(a\leftrightarrow \infty)$, take a sequence $(G_n)_n$ of finite subgraphs of $G$ exhausting $G,$ i.e., $G_n\subseteq G_{n+1}$ and $G=\bigcup\limits_{n}G_n.$ Let $Z_n$ be the vertex set of $G\setminus G_n$ and $G_n^W$ the graph obtained from $G$ by identifying $Z_n$ to a single vertex $z_n$ and removing loops (but keeping multiple edges). Call
$$\mathscr{C}(a\leftrightarrow \infty)=\mathscr{C}_{\mathbf{c}}(a\leftrightarrow \infty):=\lim\limits_{n\rightarrow\infty}\mathscr{C}_{\mathbf{c}}\left(a\leftrightarrow z_n;\,G_n^W\right)$$
the effective conductance from $a$ to $\infty$ in $G,$ and its reciprocal $\mathscr{R}(a\leftrightarrow\infty)=\frac{1}{\mathscr{C}(a\leftrightarrow \infty)}$ the effective resistance.
Recall that on connected network $(G,\,\mathbf{c}),$
$$(X_n)_{n\geq 0}\ \mbox{is transient (resp.\,recurrent)}\ \Longleftrightarrow \mathscr{C}(x\leftrightarrow \infty)>0\ (\mbox{resp}.\,=0)\ \mbox{for any vertex}\ x.$$

\begin{lem}[Rayleigh's monotonicity principle]\label{lem-ray}
Let $G$ be a connected graph with two conductances $\mathbf{c}$ and $\mathbf{c}'$ such that $\mathbf{c}(e)\leq \mathbf{c}'(e),\,e\in E$.
\begin{enumerate}[{\bf (i)}]
\item For finite $G$ and any its two disjoint vertex subsets $A$ and $Z$,
\[\mathscr{C}_{\mathbf{c}}(A\leftrightarrow Z) \leq \mathscr{C}_{\mathbf{c}'}(A\leftrightarrow Z).\]

\item For infinite $G$ and any its vertex $a,$
\[\mathscr{C}_{\mathbf{c}}(a\leftrightarrow \infty) \leq \mathscr{C}_{\mathbf{c}'}(a\leftrightarrow \infty).\]
In particular, $(G,\,\mathbf{c})$ is transient implies so is $(G,\,\mathbf{c}^\prime)$ (equivalently, the recurrence of $(G,\,\mathbf{c}')$ implies that of $(G,\,\mathbf{c})$).
\end{enumerate}
\end{lem}
\vskip 2mm

\begin{lem}[The Nash-Williams inequality and recurrence criterion]\label{lem-Nash}
For any distinct vertices $a$ and $z$ separated by pairwise disjoint cutsets $\Pi_1,\, \cdots,\, \Pi_{n}$ in a finite network,
\[
\mathscr{R}(a\leftrightarrow z)\geq \sum\limits_{k=1}^{n} \left(\sum\limits_{e\in \Pi_k}\mathbf{c}(e)\right)^{-1}.
\]
For any sequence $\left\{\Pi_n\right\}_n$ of pairwise disjoint finite cutsets in an infinite locally finite network $G$ such that
each $\Pi_n$ separates $a$ from $\infty,$
\[
\mathscr{R}(a\leftrightarrow \infty)\geq \sum\limits_{n}^{\infty} \left(\sum\limits_{e\in \Pi_n}\mathbf{c}(e)\right)^{-1};
\]
and particularly $G$ is recurrent when the right-hand side is $\infty.$
\end{lem}
\vskip 2mm

\begin{lem}[Energy transience criterion \mbox{\cite[Theorem 2.11]{LP2017}}]\label{lem-energy}
Connected infinite network $(G,\, \mathbf{c})$ is transient if and only if there exists an unit flow from some (every) vertex to $\infty$ with finite energy.
\end{lem}
\vskip 2mm

Given two networks $G=((V,\,E),\, \mathbf{c})$ and $G'=((V',\,E'),\,\mathbf{c}')$. Say $G$ can be roughly embedded into $G'$ if there exists a map $\phi:\,V\longmapsto V'$ such that there are constants $\alpha,\beta<\infty$ and
a map $\Phi$ mapping oriented edges $xy$ in $G$ to a non-empty simple oriented path $\Phi(xy)$ in $G'$ from $\phi(x)$ to $\phi(y)$ such that
\begin{eqnarray*}
&&\mbox{$\sum\limits_{e'\in \Phi(xy)} \mathbf{r}'(e')\leq \alpha \mathbf{r}(\{x,y\})$ and $\Phi(yx)$ is the reverse of $\Phi(xy)$; and for any $e'\in E'$,}\\
&&\mbox{there are no more than $\beta$ edges in $G$ whose image under $\Phi$ contains $e'$.}
\end{eqnarray*}
Call $G$ and $G'$ are roughly equivalent if $G$ and $G'$ can be roughly embedded into each other.

\begin{lem}[Rough embeddings and transience \mbox{\cite[Theorem 2.17]{LP2017}}] For two roughly equivalent connected networks $G$ and $G'$, $G$ is transient iff so is $G'$. In fact, if there is a rough embedding from $G$ to $G'$, then $G$ is transient implies so is $G'$ (equivalently, the recurrence of $G'$ implies the recurrence of $G$).
\end{lem}

\subsection{Proof of Theorem \ref{generalgraph01}}\label{quasi}
\noindent
To begin, define a geodesic spanning tree $\mathcal{T}$ on a quasi-transitive graph $G$ as follows:

{\bf (i)} Define an order for oriented edges adjacent to each vertex in $G$. Under group action of automorphism group ${\rm Aut}(G)$ of $G$, there are only finitely many orbits $\{\mathcal{O}_i \}_{i=1}^k$. Let $\mathcal{O}_1=\{y_1,\, y_2,\,\ldots\}$. Choose an order `$<$' for all oriented edges $y_1\cdot$ starting at $y_1$, and a sequence $\{\phi_i\}_{i=2}^{\infty}\subseteq {\rm Aut}(G)$ such that $\phi_i(y_1)= y_i\in \mathcal{O}_1$ for any $i\geq 2.$ Then there is a natural way to define an order of oriented edges starting at $\phi_i(y_1)$, namely $\phi_i(y_1)\phi_i(u)< \phi_i(y_1)\phi_i(v)$ if and only if $y_1u<y_1v$. Here $y_1\sim u$ and $y_1\sim v$. Then define similarly `$<$' for
oriented edges starting at vertices in other orbits. Finally, at each vertex of $G$, all oriented edges starting at this vertex have a well-defined order `$<$'.

{\bf (ii)} Notice $o$ is the root of $G$. Then each $x\in V$ can be uniquely denoted by lexicographically minimal finite words of vertices in $G$ as $x= \gamma_0\gamma_1\gamma_2 \cdots\gamma_{|x|}$. Here $\gamma_0= o$ and $\gamma_i\in V,\, 1\leq i\leq |x|$, and $\vert x\vert$ is the graph distance between $x$ and $o.$ {\it In fact}, for any finite words $x= \beta_0\beta_1\beta_2 \cdots \beta_{|x|}$ such that
$$\beta_0=o,\ \{\beta_0,\, \beta_1,\, \beta_2,\, \ldots,\, \beta_{|x|}\}\neq \{\gamma_0,\, \gamma_1,\, \gamma_2,\, \ldots,\, \gamma_{|x|}\},$$
there must be some $0< s\leq |x|$ such that $\gamma_i= \beta_i$ for $0\leq i \leq s-1$ and $\gamma_s\neq \beta_s$. Lexicographical
minimality implies that $\gamma_{s-1}\gamma_s<\beta_{s-1}\beta_{s}.$ Denote lexicographically minimal finite words representation of $x$ by $w_x= w_x(0)w_x(1)\cdots w_x(\vert x\vert).$

{\bf (iii)} A geodesic spanning tree $\mathcal{T}$ of $G$ is a subgraph of $G$ with no loop and contains all vertices in $G$ such that there is an edge between any two vertices $x$ and $y$ of $\mathcal{T}$ iff $|x|= |y|+1$ and $w_y(j)=w_x(j)$ for $0\leq j\leq |y|$ or $|y|= |x|+1$ and $w_x(j)=w_y(j)$ for $0\leq j\leq |x|.$ Now the construction of $\mathcal{T}$ is done.\\

Given a locally finite infinite tree $T$ with root $o.$ Recall branching number of $T$ is defined as
$${\rm br}(T)=\sup\limits\left\{\lambda\geq 0:\ \exists\ \mbox{a nonzero flow}\ \theta\ \mbox{on}\ T\ \mbox{such that}\ \vert\theta\vert (e)\leq \lambda^{-\vert e\vert},\ \forall\, \mbox{directed edge}\ e\right\};$$
and by the max-flow min-cut theorem,
$${\rm br}(T)=\sup\left\{\lambda\geq 0:\ \inf\limits_{\Pi}\sum\limits_{e\in\Pi}\lambda^{-\vert e\vert}>0\right\},$$
where the $\inf$ is over all cutsets $\Pi$ separating $o$ from $\infty.$ By \cite[Theorem 3.5]{LP2017}, ${\rm RW}_{\lambda}$ on $T$ is transient if $\lambda<{\rm br}(T)$ and recurrent if $\lambda>{\rm br}(T);$ and by \cite[Theorem 5.15]{LP2017},
\begin{eqnarray}\label{eq-p_c-br-tree}
p_c(T)=\frac{1}{{\rm br}(T)}.
\end{eqnarray}

\begin{lem}\label{brgr}
For a locally finite quasi-transitive infinite connected graph $G$ with a geodesic spanning tree $\mathcal{T}$,
\begin{equation}\label{quasi-relation}
\lambda_c(G)={\rm br}(\mathcal{T})= {\rm gr}(\mathcal{T})={\rm gr}(G).
\end{equation}
\end{lem}
\pf
By quasi-transitivity, the geodesic spanning tree $\mathcal{T}$ of $G$ is a sub-periodic tree. Then \cite[Theorem 3.8]{LP2017} implies that growth rate ${\rm gr}(\mathcal{T})$ exists and ${\rm gr}(\mathcal{T})={\rm br}(\mathcal{T})$. Note that the graph distances between $o$ and any vertex are the same in tree $\mathcal{T}$ and in original graph $G$. Thus ${\rm gr}(\mathcal{T})={\rm gr}(G)$. Recall $\lambda_c(G)={\rm gr}(G)$. We obtain the lemma immediately.
\qed
\vskip 2mm
\begin{lem}[\mbox{\cite[Proposition 6.1]{RL1990}}]\label{gpercolation}
Assume $G$ is a locally finite connected infinite graph and $G({\omega}_p)$ the open subgraph of $G$ in Bernoulli-$p$ bond percolation $\omega_p$ with $p\in [0,\,1].$ Given $\omega_p,$ let
$$p_c\left(G({\omega}_p)\right)=\sup\left\{q:\ \mathbb{P}\left[\mbox{Bernoulli-$q$ bond percolation on}\ G(\omega_p)\ \mbox{has an infinite cluster}\right]=0\right\}.$$
Then
\[
p_c\left(G({\omega}_p)\right)=\left(p_c(G)/p\right)\wedge 1 \ \mbox{a.s.}
\]
\end{lem}
\vskip 2mm

For a tree $\Gamma$ with root $o$, let $\Gamma^{\sigma}= \{\pi\in \Gamma: \sigma\leq \pi\}$ denote the subtree of $\Gamma$ with $\sigma$ and all its descendents. Recall from \cite[Corollary 6.3]{RL1990} that if $K_{\sigma}(\omega_p)$ denotes the cluster of $\sigma$ in $\omega_p$, then when $p> ({\rm br}(\Gamma))^{-1},$
\begin{equation}\label{ppercolation}
\sup\limits_{\sigma\in \Gamma}{\rm br}(K_{\sigma}(\omega_p))= \sup\limits_{\sigma\in \Gamma} {\rm br}(K^{\sigma}(\omega_p))= p \cdot {\rm br}(\Gamma) \ \ a.s.,
\end{equation}
where $K^{\sigma}(\omega_p)=\Gamma^{\sigma}\cap K_{\sigma}(\omega_p)$ and the branching number of a finite tree is regarded as zero. Therefore,
$${\rm ess} \sup {\rm br}(K_{o}(\omega_p))= p\cdot {\rm br}(\Gamma).$$
\vskip 2mm

\noindent{\bf Proof of Theorem \ref{generalgraph01}.} {\bf (i)} Fix $p\in [0,\,1].$ Define an event on percolation configuration space:
\[
A_p=\{(G,\, \mathbf{C}_{\lambda_1},\, \mathbf{C}_{\lambda_2},\, p)\ \mbox{is recurrent}\}.
\]
Then $A_p$ is an invariant event under group action of ${\rm Aut}(G)$, and thus $\mathbb{P}(A_p)\in\{0,1\}$ by the ergodicity of Bernoulli percolation $\omega_p$ (\cite[Proposition~7.3]{LP2017}),
namely either $(G,\, \mathbf{C}_{\lambda_1},\, \mathbf{C}_{\lambda_2},\, p)$ is a.s.\,recurrent or a.s.\,transient.
{\it In fact}, let $\mathbf{C}_{\omega_p}$ be the conductance function of $(G,\, \mathbf{C}_{\lambda_1},\, \mathbf{C}_{\lambda_2},\, p),$ then for any $\gamma\in {\rm Aut}(G),$ when $e$ is open,
\begin{eqnarray*}
\mathbf{C}_{\gamma\omega_p}(\gamma e)=\lambda_1^{-\vert\gamma e\vert}\in\lambda_1^{-\vert e\vert}\left[(\lambda_1\wedge 1)^{\vert\gamma o\vert}\wedge(\lambda_1\vee 1)^{-\vert\gamma o\vert} ,\, (\lambda_1\vee 1)^{\vert\gamma o\vert}\vee(\lambda_1\wedge 1)^{-\vert\gamma o\vert}\right],
\end{eqnarray*}
and when $e$ is closed,
\begin{eqnarray*}
\mathbf{C}_{\gamma\omega_p}(\gamma e)=\lambda_2^{-\vert\gamma e\vert}\in\lambda_2^{-\vert e\vert}\left[\lambda_2^{-\vert\gamma o\vert},\, \lambda_2^{\vert\gamma o\vert}\right].
\end{eqnarray*}
Write $\mathbf{C}^\prime_{\gamma\omega_p}(e)=\mathbf{C}_{\gamma\omega_p}(\gamma e),\ e\in E.$ Then
$$\mathbf{C}^\prime_{\gamma\omega_p}(e)\asymp \mathbf{C}_{\omega_p}(e),\ e\in E,$$
namely networks $(G,\,\mathbf{C}^{\prime}_{\gamma\omega_p})$ and $(G,\,\mathbf{C}_{\omega_p})$ are equivalent.
Notice $(G,\,\mathbf{C}^{\prime}_{\gamma\omega_p})$ is the image network of $(G,\,\mathbf{C}_{\gamma\omega_p})$ under the automorphism $\gamma.$
Therefore, $(G,\,\mathbf{C}_{\gamma\omega_p})$ is recurrent iff so is $(G,\,\mathbf{C}_{\omega_p}),$ and further $A_p$ is an ${\rm Aut}(G)$-invariant event.

Note that for any $0\leq p\leq q\leq 1,$
$$\omega_p(e)\leq\omega_q(e),\ \mathbf{C}_{\omega_p}(e)\leq\mathbf{C}_{\omega_q}(e),\ \forall\,e\in E;$$
and by the Rayleigh's monotonicity principle (Lemma \ref{lem-ray}), $(G,\,\mathbf{C}_{\omega_p})$ is recurrent if so is  $(G,\,\mathbf{C}_{\omega_q}).$
On one hand, remembering $(G,\,\mathbf{C}_{\omega_0})=(G,\,\mathbf{C}_{\lambda_2})$ is recurrent, we have that
\begin{eqnarray*}
p_c^*=\sup\{p\geq 0:\ \mathbb{P}[(G,\,\mathbf{C}_{\omega_p})\ \mbox{is recurrent}]=1\}
\end{eqnarray*}
satisfies that for any $p<p_c^*$, $(G,\,\mathbf{C}_{\omega_p})$ is a.s. recurrent; and for any $p>p_c^*,$  $(G,\,\mathbf{C}_{\omega_p})$ is a.s. transient.
On the other hand, again by the Rayleigh's monotonicity principle (Lemma \ref{lem-ray}), for any percolation environment $\omega,$
there exists $p_c^{*}(\omega)\in[0,\, 1]$ such that $(G,\, \mathbf{C}_{\omega_p})$ is recurrent for any $p< p_c^{*}(\omega)$ and is transient for any $p> p_c^{*}(\omega)$. So
\begin{eqnarray*}
&&\mbox{almost surely, $p_c^*=p_c(\omega),$ and $(G,\, \mathbf{C}_{\omega_p})$ is recurrent for any $p< p_c^{*}$}\\
&&\mbox{and is transient for any $p> p_c^{*}$.}
\end{eqnarray*}

{\it Indeed}, for any rational number $p\in [0,\,p_c^*)$ (if exists), almost surely $(G,\,\mathbf{C}_{\omega_p})$ is recurrent; and hence almost surely $p_c^*(\omega)\geq p.$ Let $p\uparrow p_c^*,$ we get that $p_c^*(\omega)
\geq p_c^*$ almost surely. Additionally, for any rational number $q\in (p_c^*,\,1]$ (if exists), almost surely $(G,\,\mathbf{C}_{\omega_q})$ is transient; and thus almost surely $p_c^*(\omega)\leq q.$ Let $q\downarrow p_c^*,$
we have that $p_c^*(\omega)\leq p_c^*$ almost surely. Therefore, almost surely, $p_c^*(\omega)=p_c^*.$
\vskip 2mm

{\bf (ii)} Assume $p_c=p_c(G)\in (0,1).$ To prove $p_c^*\in (0,1).$

Let $\mathcal{T}$ be a geodesic spanning tree of quasi-transitive graph $G$. Then by (\ref{quasi-relation}),
\begin{equation*}
\lambda_c(G)={\rm br}(\mathcal{T})={\rm gr}(\mathcal{T})={\rm gr}(G).
\end{equation*}

Assume firstly $0<\lambda_1<\lambda_c(G)=1.$ Then for any $p>p_c$, almost surely, there is an infinite open cluster $K$ in $\omega_p.$ By the energy transience criterion (Lemma \ref{lem-energy}),
$(K,\,\mathbf{C}_{\omega_p})=(K,\,\mathbf{C}_{\lambda_1})$ is transient. Thus from the Rayleigh's monotonicity principle (Lemma \ref{lem-ray}),
$(G,\,\mathbf{C}_{\omega_p})$ is transient; and further
$$p_c^*\leq p_c<1.$$

Suppose $0<\lambda_1<\lambda_c(G)$ and $\lambda_c(G)>1.$ Take $p>\frac{\lambda_1\vee 1}{\lambda_c(G)}$ and $\varepsilon\in (0,1)$ such that
$$(p-\varepsilon)\lambda_c(G)>\lambda_1\vee 1.$$
Then by Lemma \ref{gpercolation} and \eqref{ppercolation}, almost surely, there exists an infinite open cluster $K_{\sigma}(\omega_{p,\mathcal{T}})$ of some $\sigma\in\mathcal{T}$ in percolation $\omega_{p,\mathcal{T}}$, which is the restriction of $\omega_p$ to $\mathcal{T},$ such that
$${\rm br}\left(K_{\sigma}(\omega_{p,\mathcal{T}})\right)>(p-\varepsilon){\rm br}(\mathcal{T})=(p-\varepsilon)\lambda_c(G)>\lambda_1\vee 1.$$
Fix such an $\omega_{p}$. Note $K_{\sigma}(\omega_{p,\mathcal{T}})$ is a tree. Hence $\left(K_{\sigma}(\omega_{p,\mathcal{T}}),\,\mathbf{C}_{\lambda_1}\right)=\left(K_{\sigma}(\omega_{p,\mathcal{T}}),\,\mathbf{C}_{\omega_p}\right)$ is transient.
Then by the Rayleigh's monotonicity principle (Lemma \ref{lem-ray}), $(G,\,\mathbf{C}_{\omega_p})$ is transient; and further
$$p_c^*\leq\frac{\lambda_1\vee 1}{\lambda_c(G)}<1.$$

Now we are in the position to prove $p_c^*>0.$ Fix any $0<\varepsilon <\lambda_2- \lambda_c(G)$ and let
$$L_{k}:= \{x\in G:\ |x|= k\},\ k\in\mathbb{N}.$$
Use $A\stackrel{\omega_p}{\longleftrightarrow}B$ to denote vertex sets $A$ and $B$ are connected to each other in $\omega_p.$ Given any vertex $x\in G.$ Let $a_n(x)$ be the number of self-avoiding walks (paths) on $G$ with length $n$ starting at $x.$ Then due to $G$ is quasi-transitive,
$$\mu=\lim\limits_{n\rightarrow\infty}\sqrt[n]{a_n(x)}\in [1,\,\infty)\ \mbox{exists and independ of}\ x;$$
and call $\mu$ the connective constant of $G.$ Choose constant $C\in (0,\,\infty)$ such that
$$\vert B_n(o)\vert\leq C(\lambda_c(G)+\varepsilon)^n,\  a_n(x)\leq C(\mu+\varepsilon)^n,\ n\in\mathbb{N},\ x\in G.$$
Then for any $\alpha\in (1,\,\infty)$ with $\lambda_2^{1/\alpha}>\lambda_c(G)+\varepsilon,$ when $p\in \left[0,\,(\mu+\varepsilon)^{-1}(\lambda_c(G)+\varepsilon)^{-\frac{\alpha}{\alpha-1}}\right),$
\begin{equation*}
\begin{split}
\sum\limits_{n= n_0}^{\infty} \P\left(L_{\alpha^n}\stackrel{\omega_p}{\longleftrightarrow}L_{\alpha^{n+1}}\right)&\leq \sum\limits_{n=n_0}^{\infty} \sum\limits_{x\in L_{\alpha^{n+1}}} \P\left(x
\stackrel{\omega_p}{\longleftrightarrow}L_{\alpha^n}\right)\\
&\leq  \sum\limits_{n=n_0}^{\infty}C(\lambda_c(G)+\varepsilon)^{\alpha^{n+1}}\sum\limits_{j\geq \alpha^{n+1}-\alpha^n-1}p^jC(\mu+\varepsilon)^j \\
&=\frac{C^2}{1-p(\mu+\varepsilon)}\sum\limits_{n= n_0}^{\infty}\left[(p(\mu+\varepsilon))^{1-\alpha^{-1}-\alpha^{-(n+1)}}(\lambda_c(G)+\varepsilon)\right]^{\alpha^{n+1}}\\
&<\infty,
\end{split}
\end{equation*}
where $\mathbb{N}\ni n_0>\left(-\frac{\log(\alpha-1)}{\log\alpha}\right)\vee 0,$ each $L_{\alpha^n}$ is viewed as $L_{\lfloor\alpha^n\rfloor}$ with $\lfloor\alpha^n\rfloor$ being the integer part of
$\alpha^n.$ By the Borel-Cantelli lemma, almost surely, we can find a sequence $\{\Pi_n\}_{n=n_1}^{\infty}$ of minimum closed cutsets such that each $\Pi_n$ is between $L_{\alpha^n}$ and $L_{\alpha^{n+1}}$,
where $n_1$ is a large random natural number. Fix such a percolation configuration. By the Nash-Williams recurrence criterion (Lemma \ref{lem-Nash}),
\begin{equation*}
\begin{split}
\mathscr{R}(o\leftrightarrow \infty)&\geq \sum\limits_{n=n_1}^{\infty} \left(\sum\limits_{e\in \Pi_n} \mathbf{C}_{\lambda_2}(e)\right)^{-1}\geq \sum\limits_{n=n_1}^{\infty}\left(\lambda_2^{-\lfloor\alpha^n\rfloor} d C(\lambda_c(G)+\varepsilon)^{\alpha^{n+1}}\right)^{-1}\\
&\geq  \sum\limits_{n=n_1}^{\infty}\left(\lambda_2^{-\alpha^n+1} d C(\lambda_c(G)+\varepsilon)^{\alpha^{n+1}}\right)^{-1}\\
&= \frac{1}{\lambda_2 dC}\sum\limits_{n=n_1}^{\infty}\left(\frac{\lambda_2^{1/\alpha}}{\lambda_c(G)+\varepsilon}\right)^{\alpha^{n+1}}= \infty,
\end{split}
\end{equation*}
and $(G,\,\mathbf{C}_{\omega_p})$ is recurrent. Here we have used that $\vert\Pi_n\vert$ is no more than the number of edges in $B_{\lfloor\alpha^{n+1}\rfloor}(o)\setminus B_{\lfloor\alpha^n\rfloor}(o)$,
and clearly the latter is at most $d\left\vert B_{\lfloor\alpha^{n+1}\rfloor}(o)\right\vert$ with $d$ being the maximum of vertex degrees of $G.$
Therefore, $p_c^*>0.$
\vskip 2mm

{\bf (iii)} On any Cayley graph $G$ of any group which is virtually $\Z$, by Theorem \ref{recurthm} (iv), there is no non-trivial recurrence/transience phase transition for $((G,\, \mathbf{C}_{\lambda_1},\, \mathbf{C}_{\lambda_2},\, p))_{p\in [0,\,1]}$, i.e. $p_{c}^{*}=p_c= 1.$
\qed
\vskip 2mm

\subsection{Proofs of Theorems \ref{recurthm}, \ref{generald} and \ref{thm-current-unique}}\label{rzd}
\noindent\textbf{Proof of Theorem \ref{recurthm}\,(i)-(ii).}
{\bf (i)} Write $e_{i}^{+}$ and $e_{i}^{-}$ for the directed edges from $i$ to $i+1$ and $i-1$
respectively for any $i\in\Z.$ Note that any unit flow $\theta$ from $i_0$ to infinity on $\Z$ must have the following form:
\begin{equation}\label{eq-flow-structure}
\begin{split}
&\mbox{For some constant $a\in\mathbb{R}$,}\
\theta(e_{i_0+i}^{+})=a,\ \theta(e_{i_0-i}^{-})=1-a,\ i\in \mathbb{Z}_+;\ \mbox{and $\theta$ is}\\
&\mbox{the only unit flow on $i_0+\Z_{+}$ (resp. $i_0-\Z_{+}$) if $a=1$ (resp. $a=0$).}
\end{split}
\end{equation}
For any unit flow $\theta$ on $\Z$ and any $p\in [0,\,1],$ the energy $\mathscr{E}_p(\theta)$ of $\theta$ on $(\Z,\,\mathbf{c}_1,\,\mathbf{c}_2,\,p)$ satisfies that
\begin{eqnarray*}
\mathscr{E}_p(\theta)= \sum\limits_{e\in E(\Z)} \frac{\theta^2(e)}{\mathbf{c}_1(e)}I_{\{U_e\leq p\}}+\sum\limits_{e\in E(\Z)} \frac{\theta^2(e)}{\mathbf{c}_2(e)}I_{\{U_e>p\}}
  \geq \sum\limits_{e\in E(\Z)} \frac{\theta^2(e)}{\mathbf{c}_2(e)}I_{\{U_e>p\}}.
\end{eqnarray*}
And $\sum\limits_{e\in E(\Z)} \frac{\theta^2(e)}{\mathbf{c}_2(e)}I_{\{U_e>p\}}$ is a decreasing function in $p.$ So to prove that
\begin{eqnarray}\label{eq-flow-energy-Z}
\mbox{almost surely},\ \mbox{for any unit flow}\ \theta\ \mbox{on}\ \Z\ \mbox{and any}\ p\in [0,1),\ \mathscr{E}_p(\theta)=\infty,
\end{eqnarray}
it suffices to prove that for any fixed $p\in [0,\,1),$ almost surely,
\begin{eqnarray}\label{eq-energy-Z}
\sum\limits_{i\in i_0+\Z_+} \frac{1}{\mathbf{c}_2(\{i,i+1\})}I_{\{U_{\{i,i+1\}}>p\}}=\sum\limits_{i\in i_0-\Z_+} \frac{1}{\mathbf{c}_2(\{i,i-1\})}I_{\{U_{\{i,i-1\}}>p\}}=\infty,\ i_0\in\Z.
\end{eqnarray}

Once \eqref{eq-energy-Z} is true, then \eqref{eq-flow-energy-Z} holds; and let $\mathbf{c}_p$ be the conductance function of $(\Z,\, \mathbf{c}_1,\, \mathbf{c}_2,\, p)$ and
$\Z(\mathbf{c}_p)$ the graph on $\Z$ induced by $\{e\in E(\Z):\ \mathbf{c}_p(e)>0\}$ for any $p\in [0,1);$ and fix a percolation environment satisfying \eqref{eq-flow-energy-Z}.
Note a connected component of network $(\Z,\,\mathbf{c}_p)$ means a connected component of $\Z(\mathbf{c}_p).$ Trivially the associated random walk on every finite connected component of network $(\Z,\,\mathbf{c}_p)$ is recurrent. While for the associated random walk on every infinite connected component of network $(\Z,\,\mathbf{c}_p),$ note \eqref{eq-flow-structure} and \eqref{eq-flow-energy-Z}, by Lemma \ref{lem-energy},
it is recurrent. Namely, assuming \eqref{eq-energy-Z}, almost surely, for any $p\in [0,\,1),$ $(\Z,\, \mathbf{c}_1,\, \mathbf{c}_2,\, p)$ is recurrent. We prove \eqref{eq-energy-Z} in two steps.

{\bf (i.1)} We prove firstly
\begin{eqnarray}\label{eq-energy-Z-2}
\sum\limits_{i\in i_0+\Z_+} \frac{1}{\mathbf{c}_2(\{i,i+1\})}=\sum\limits_{i\in i_0-\Z_+} \frac{1}{\mathbf{c}_2(\{i,i-1\})}=\infty,\ i_0\in\Z.
\end{eqnarray}

Let $\Z({\mathbf{c}_2})$ be the graph on $\Z$ induced by $\{e\in E(\Z):\ \mathbf{c}_2(e)>0\}.$ When $(\Z,\,\mathbf{c}_2)$ is connected (i.e., $\Z(\mathbf{c}_2)$ is the graph $(\Z,\,E(\Z))$),
by \eqref{eq-flow-structure}, recurrence of connected network $(\Z,\,\mathbf{c}_2)$ and Lemma \ref{lem-energy}, \eqref{eq-energy-Z-2} is true.
Clearly \eqref{eq-energy-Z-2} holds if all connected components of $\Z(\mathbf{c}_2)$ are finite. If $\Z(\mathbf{c}_2)$ is not connected and has only one infinite connected component, say $i_1+\Z_+,$
then trivially
\begin{eqnarray*}
\sum\limits_{i\in i_0-\Z_+} \frac{1}{\mathbf{c}_2(\{i,i-1\})}=\infty,\ i_0\in\Z;
\end{eqnarray*}
and by recurrence of connected network $(i_1+\Z_+,\,\mathbf{c}_2)$, Lemma \ref{lem-energy} and \eqref{eq-flow-structure},
\begin{eqnarray}\label{eq-energy-Z-3}
\sum\limits_{i\in i_1+\Z_+} \frac{1}{\mathbf{c}_2(\{i,i+1\})}=\infty,\ \mbox{and further}\ \sum\limits_{i\in i_0+\Z_+} \frac{1}{\mathbf{c}_2(\{i,i+1\})}=\infty,\ i_0\in \Z.
\end{eqnarray}
If $\Z(\mathbf{c}_2)$ is not connected and has two infinite connected components, $i_1+\Z_+$ and $i_2-\Z_{+}$ with $i_2<i_1,$ similarly to \eqref{eq-energy-Z-3}, one can prove
\eqref{eq-energy-Z-2}.

{\bf (i.2)} By \eqref{eq-energy-Z-2}, for any $i_0\in \Z,$
\begin{eqnarray*}
\sum\limits_{i\in i_0+\Z_+}\left(1-e^{-1/\mathbf{c}_2(\{i,i+1\})}\right)=\infty
=\sum\limits_{i\in i_0-\Z_+}\left(1-e^{-1/\mathbf{c}_2(\{i,i-1\})}\right),
\end{eqnarray*}
and further for any $p\in [0,\,1),$
\begin{eqnarray*}
\prod\limits_{i\in i_0+\Z_+}\left\{(1-p)\left(e^{-1/\mathbf{c}_2(\{i,i+1\})}-1\right)+1\right\}=0
  =\prod\limits_{i\in i_0-\Z_+}\left\{(1-p)\left(e^{-1/\mathbf{c}_2(\{i,i-1\})}-1\right)+1\right\}.
\end{eqnarray*}
Namely for any $p\in [0,\,1),$
\begin{eqnarray*}
\mathbb{E}\left[\prod\limits_{i\in i_0+\Z_+}e^{-\frac{1}{\mathbf{c}_2(\{i,i+1\})}I_{\left\{U_{\{i,i+1\}}>p\right\}}}\right]=0=
\mathbb{E}\left[\prod\limits_{i\in i_0-\Z_+}e^{-\frac{1}{\mathbf{c}_2(\{i,i-1\})}I_{\left\{U_{\{i,i-1\}}>p\right\}}}\right],
\end{eqnarray*}
where we have used that
\begin{eqnarray*}
&&\mathbb{E}\left[\prod\limits_{i\in i_0+\Z_+}e^{-\frac{1}{\mathbf{c}_2(\{i,i+1\})}I_{\left\{U_{\{i,i+1\}}>p\right\}}}\right]=
   \prod\limits_{i\in i_0+\Z_+}\left\{(1-p)\left(e^{-1/\mathbf{c}_2(\{i,i+1\})}-1\right)+1\right\},\\
&&\mathbb{E}\left[\prod\limits_{i\in i_0-\Z_+}e^{-\frac{1}{\mathbf{c}_2(\{i,i-1\})}I_{\left\{U_{\{i,i-1\}}>p\right\}}}\right]
=\prod\limits_{i\in i_0-\Z_+}\left\{(1-p)\left(e^{-1/\mathbf{c}_2(\{i,i-1\})}-1\right)+1\right\}.
\end{eqnarray*}
Therefore, \eqref{eq-energy-Z} is true.
\vskip 2mm

{\bf (ii)} Since $\Z$ is Abelian, so any its Cayley graph corresponds to a certain symmetric generating set. Assume the generating set of $G$ is
$$S=\{\pm a_i:\ 1\leq i\leq\ell,\ 0\leq a_1<a_2<\ldots<a_\ell\}.$$
\begin{figure}[!htp]
\centering
\includegraphics[width=9cm, height= 1.5cm]{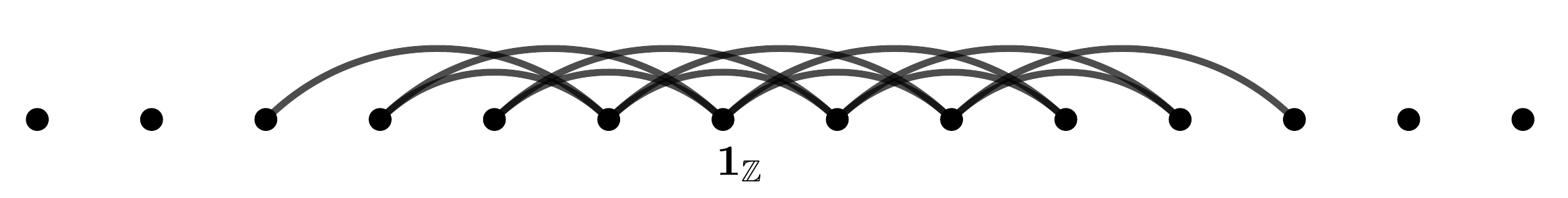}
\caption{A Cayley graph on $\Z$ with generating set $\{\pm 2, \pm 3\}$.}
\label{gcayleyz}
\end{figure}

Given any $p\in (0,\,1).$ For any $x\in\Z_+,$ let
$$A_x(p)=\left\{\forall\,y\in (x,\,x+a_{\ell}],\ \{y,y+a_i\}\ \mbox{is closed in}\ \omega_p\ \mbox{for any}\ 1\leq i\leq\ell\right\};$$
and for any negative $x\in\Z,$ let
$$A_x(p)=\left\{\forall\,y\in [x-a_{\ell},\,x),\ \{y,y-a_i\}\ \mbox{is closed in}\ \omega_p\ \mbox{for any}\ 1\leq i\leq\ell\right\}.$$
Clearly, $\{A_{ka_{\ell}}(p)\}_{k\in\Z}$ is an i.i.d.\,sequence of events and $\mathbb{P}[A_0(p)]>0.$ By the law of large numbers, almost surely,
$$\frac{1}{n+1}\sum\limits_{k=0}^n I_{A_{ka_{\ell}}(p)}\rightarrow \mathbb{P}[A_0(p)]>0,$$
and thus infinitely many $A_{ka_{\ell}}(p)$s, say $A_{k_ja_{\ell}}(p),\ j=1,2,\ldots,$ occur. Similarly, almost surely, infinitely many $A_{-ka_{\ell}}(p)$s ($k>0$), say $A_{-s_ja_{\ell}}(p),\ j=1,2,\ldots,$ occur. Fix a percolation environment $\omega=(\omega_q)_{q\in [0,\,1]}$ such that these events hold. Then
\begin{eqnarray*}
&&\Pi_n=\left\{\{y,y+a_i\}:\ 1\leq i\leq \ell,\ y\in (k_na_{\ell},\,(k_n+1)a_{\ell}]\right\}\\
&&\ \ \ \ \ \ \ \ \ \ \cup \left\{\{y,y-a_i\}:\ 1\leq i\leq \ell,\ y\in [-(s_n+1)a_{\ell},\,-s_na_{\ell})\right\}
\end{eqnarray*}
is a set of closed edges in $\omega_p$ (and all $\omega_q$ with $q\leq p$) and a cutset of $G$ separating $0$ from infinity. Then effective resistance $\mathscr{R}_q(0\leftrightarrow\infty)$ between $0$ and infinity of $(G,\,\mathbf{c}_1,\,\mathbf{c}_2,\,q)$ with $q\leq p$ satisfies that
\begin{eqnarray*}
\mathscr{R}_q(0\leftrightarrow\infty)\geq \sum\limits_{n=1}^{\infty}\left(\sum\limits_{e\in\Pi_n}\mathbf{c}_2(e)\right)^{-1}
   \geq\sum\limits_{n=1}^{\infty}\left(\sum\limits_{e\in\Pi_n}c\right)^{-1}
   =\frac{1}{c}\sum\limits_{n=1}^{\infty}\frac{1}{\vert\Pi_n\vert}=\frac{1}{c}\sum\limits_{n=1}^{\infty}\frac{1}{\vert\Pi_1\vert}=\infty.
\end{eqnarray*}
By the Nash-Williams criterion, for the percolation environment $\omega,$ all networks $(G,\,\mathbf{c}_1,\,\mathbf{c}_2,\,q)$ with $q\leq p$ are recurrent. Due to $p$ is arbitrary, we are done.
\qed
\vskip 2mm


\begin{figure}
\centering
\subfigure{

\includegraphics[width=7cm, height= 3cm]{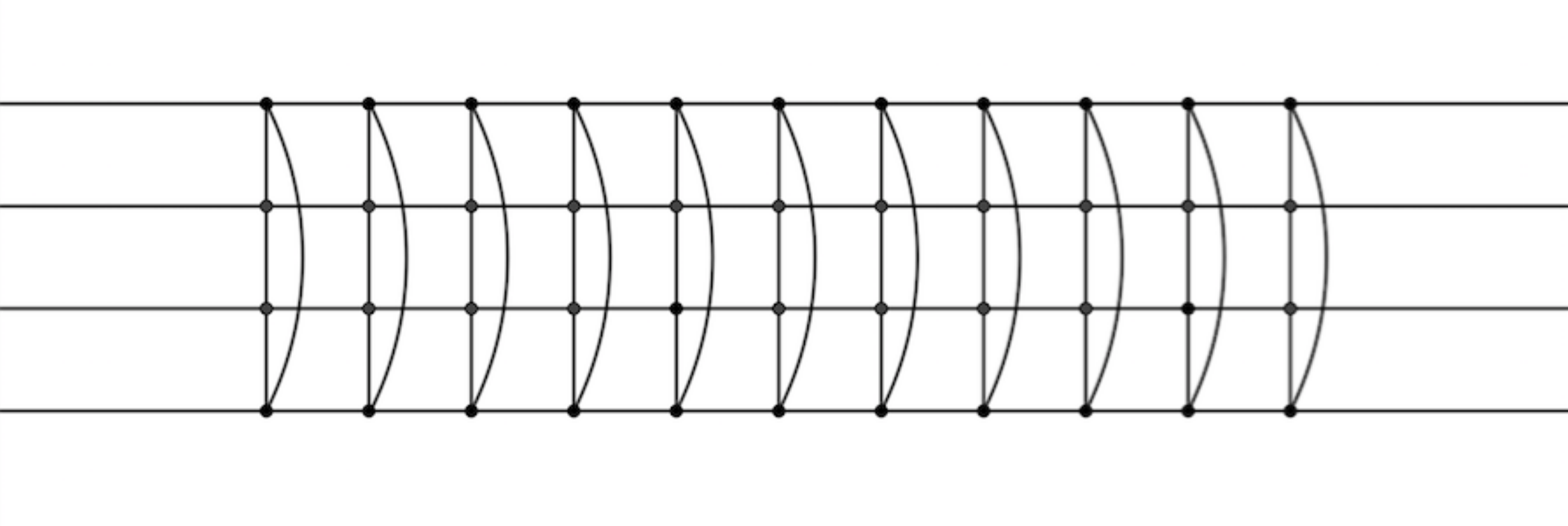}}
\hspace{0cm}
\subfigure{

\includegraphics[width=3cm]{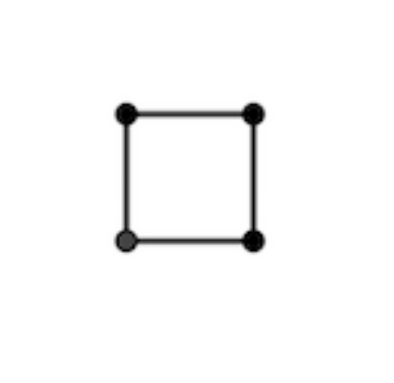}}
\caption{An example of a ladder graph $\Z \times G$; $G$ is isomorphic to a square.}
\label{graphladder}
\end{figure}



\noindent\textbf{Proof of Theorem \ref{recurthm}\,(iii).}
{\bf Step 1.} Consider an infinite connected multi-graph $H=(V(H),\,E(H))$ with loop edges. Let $(H,\,\mathbf{c})$ and $(H,\,\mathbf{c}^*)$ be two networks such that $\mathbf{c}^*(e)=\mathbf{c}(e)$ for any non-loop edge $e$ and otherwise $\mathbf{c}^*(e)=0,$ and
$\sum\limits_{u\in e}\mathbf{c}(e)I_{\{e\ \mbox{is not a loop}\}}>0$ for any vertex $u\in H.$ Then $(H,\,\mathbf{c})$ is recurrent iff so is $(H,\,\mathbf{c}^*).$

{\it Indeed}, let $(X_n)_{n=0}^{\infty}$ be the random walk associated to $(H,\,\mathbf{c})$ starting at $o\in V(H).$ Define
$$\tau_0=0,\ \tau_{k+1}=\inf\{n>\tau_k:\ X_n\not=X_{\tau_{k}}\},\ k\in\Z_+.$$
Note that for any $k\in\Z_{+},$ given $X_{\tau_k}=u,$ $\tau_{k+1}-\tau_k$ obeys a geometric distribution with parameter
$$q_u=\frac{\sum\limits_{u\in e}\mathbf{c}(e)I_{\{e\ \mbox{is not a loop}\}}}{\sum\limits_{u\in e}\mathbf{c}(e)}>0.$$
This implies that each $\tau_k$ is finite a.s.. It is easy to see that $(X_{\tau_k})_{k=0}^{\infty}$ is just a random walk associated to $(H,\,\mathbf{c}^*)$ starting at $o$. Thus $(X_n)_{n=0}^{\infty}$ is recurrent iff so is $(X_{\tau_k})_{k=0}^{\infty},$ the claim holds.
\vskip 2mm

{\bf Step 2.} Write $\mathbf{c}_p$ for the conductance function of $(\Z\times G,\,\mathbf{c}_1,\,\mathbf{c}_2,\,p).$
For any $x\in \Z$, identify vertices $(x,y),\, y\in G$ as one point $x$. Then we obtain naturally a multi-graph $\widehat{\Z}$ on $\Z$ with loop edges such that there are $|G|$ parallel edges between any two adjacent points in $\Z,$ and a disordered random network $\left(\widehat{\Z},\,\mathbf{\hat{c}}_1,\,\mathbf{\hat{c}}_2,\,p\right)=:\left(\widehat{\Z},\,\mathbf{\hat{c}}_p\right)$ for any $p\in [0,\,1].$ Here $\mathbf{\hat{c}}_1,\,\mathbf{\hat{c}}_2$ and $\mathbf{\hat{c}}_p$ inherit naturally from $\mathbf{c}_1,\,\mathbf{c}_2$ and $\mathbf{c}_p$ respectively. Let each $\mathbf{\hat{c}}_p^*$ be the restriction of $\mathbf{\hat{c}}_p$ to all non-loop edges. Then
by Step 1,
\begin{eqnarray}\label{eq-recurrent-equiv-1}
\mbox{each}\ \left(\widehat{\Z},\,\mathbf{\hat{c}}_p\right)\ \mbox{is recurrent iff so is each}\ \left(\widehat{\Z},\,\mathbf{\hat{c}}_p^*\right).
\end{eqnarray}

Denote by $\widetilde{\Z}$ the multi-graph obtained from $\widehat{\Z}$ by deleting all loop edges of $\widehat{\Z}.$ Clearly, each $\left(\widehat{\Z},\,\mathbf{\hat{c}}_p^*\right)$ is just the network $\left(\widetilde{\Z},\,\mathbf{\hat{c}}_p^*\right).$ Since $\mathbf{c}_1$ and $\mathbf{c}_2$ are positive functions, we have that $\mathbf{c}_p$ and $\mathbf{\hat{c}}_p$ are also positive functions, both
$(\Z\times G,\,\mathbf{c}_1,\,\mathbf{c}_2,\,p)=(\Z\times G,\,\mathbf{c}_p)$ and $\left(\widehat{\Z},\,\mathbf{\hat{c}}_1,\,\mathbf{\hat{c}}_2,\,p\right)=\left(\widehat{\Z},\,\mathbf{\hat{c}}_p\right)$ are connected
networks. Notice
$$(\Z\times G,\,\mathbf{c}_1,\,\mathbf{c}_2,\,p)=(\Z\times G,\,\mathbf{c}_p)\longmapsto \left(\widehat{\Z},\,\mathbf{\hat{c}}_1,\,\mathbf{\hat{c}}_2,\,p\right)=\left(\widehat{\Z},\,\mathbf{\hat{c}}_p\right)$$
is a rough embedding (see \cite{LP2017} p.\,44), by \eqref{eq-recurrent-equiv-1} and the version of \cite[Theorem 2.17]{LP2017} on networks with multiple edges and loop ones (which can be proved identically to \cite[Theorem 2.17]{LP2017}), if
\begin{eqnarray}\label{eq-recurrent-equiv-2}
\mbox{almost surely},\ \left(\widehat{\Z},\,\mathbf{\hat{c}}_p^*\right)=\left(\widetilde{\Z},\,\mathbf{\hat{c}}_p^*\right)\ \mbox{is recurrent for all}\ p\in [0,1),
\end{eqnarray}
then
\begin{eqnarray*}
\mbox{almost surely},\ \mbox{for any}\ p\in [0,\,1),\ \left(\Z\times G,\,\mathbf{c}_p\right)\ \mbox{is recurrent}.
\end{eqnarray*}
Therefore, it suffices to prove \eqref{eq-recurrent-equiv-2}.
\vskip 2mm

{\bf Step 3.} This step devotes to prove \eqref{eq-recurrent-equiv-2}.

For any $x\in\Z$ and $y\in G,$ write $\widetilde{e}_{x,y}$ for the parallel edge of $\widetilde{\Z}$ between $x$ and $x+1$, which comes from the edge $e_{x,y}=\{(x,y),(x+1,y)\}$ of $\Z\times G.$ Clearly, for any $p\in [0,\,1),$
$\left(\widetilde{\Z},\,\mathbf{\hat{c}}_p^*\right)$ is recurrent iff so is $(\Z,\,\mathbf{\tilde{c}}_p)$, where
$U_{x,y}:=U_{e_{x,y}},$ and $B_{x}(p)=\sum\limits_{y\in G}I_{\{U_{x,y}\leq p\}},$
\begin{eqnarray*}
\mathbf{\tilde{c}}_p(\{x,x+1\})&:=&\sum\limits_{y\in G}\mathbf{\hat{c}}_p^*(\widetilde{e}_{x,y})=
  \sum\limits_{y\in G}\left\{\mathbf{c}_1(e_{x,y})I_{\{U_{x,y}\leq p\}}+\mathbf{c}_2(e_{x,y})I_{\{U_{x,y}>p\}}\right\}\\
  &\leq & \sum\limits_{y\in G}\mathbf{c}_1(e_{x,y})I_{\{U_{x,y}\leq p\}}+\mathbf{c}_2^\prime(\{x,x+1\})(\vert G\vert-B_x(p)),\
   x\in\Z.
\end{eqnarray*}
Recall the proof of Theorem \ref{recurthm}\,(i). For any unit flow $\theta$ from $0$ to infinity on $\Z$ and any $p\in [0,\,1],$ the energy $\widetilde{\mathscr{E}}_p(\theta)$ of $\theta$ on $(\Z,\,\mathbf{\tilde{c}}_p)$ satisfies that
\begin{eqnarray*}
\widetilde{\mathscr{E}}_p(\theta)= \sum\limits_{i\in \Z} \frac{\theta^2(\{i,i+1\})}{\mathbf{\tilde{c}}_p(\{i,i+1\})}
  \geq \sum\limits_{i\in \Z} \frac{\theta^2(\{i,i+1\})}{\vert G\vert \mathbf{c}_2^\prime(\{i,i+1\})}I_{\{B_i(p)=0\}}.
\end{eqnarray*}
Note network $(\Z,\,\mathbf{\tilde{c}}_p)$ is connected for any $p\in [0,1).$ Similarly to prove Theorem \ref{recurthm}\,(i), to prove that
$$\mbox{almost surely},\ \mbox{for all}\ p\in [0,1),\ (\Z,\,\mathbf{\tilde{c}}_p)\ \mbox{is recurrent},$$
it only needs to prove that for any fixed $p\in [0,\,1),$ almost surely,
\begin{eqnarray}\label{eq-energy-Z-1}
\sum\limits_{i\in \Z_+} \frac{1}{\mathbf{c}_2^\prime(\{i,i+1\})}I_{\{B_{i}(p)=0\}}=\sum\limits_{i\in \Z_+} \frac{1}{\mathbf{c}_2^\prime(\{-i,-i-1\})}I_{\{B_{-i-1}(p)=0\}}=\infty.
\end{eqnarray}
Note that connected network $(\Z,\,\mathbf{c}_2^\prime)$ is recurrent by the assumption, and
$$\P[B_i(p)=0]=(1-p)^{\vert G\vert}\in (0,1],\ i\in\Z.$$
Similarly to \eqref{eq-energy-Z}, one can verify \eqref{eq-energy-Z-1}, which implies \eqref{eq-recurrent-equiv-2}.
\qed
\vskip 2mm

\noindent{\bf Proof of Theorem \ref{recurthm}\,(iv).}
{\bf Step 1.} By the assumption, there is a finite normal subgroup $Q$ of $\Gamma$ such that
$$\Gamma/Q\cong\Z\ \mbox{and}\ \Gamma=Q\ltimes \Gamma/Q\cong Q\ltimes \Z. $$
Assume $T=\{t_1,\,\ldots,\,t_k\}$ is the generating set corresponding to $G.$ Then multiple set
$$\langle T\rangle=\{t_1Q,\,\ldots,\,t_kQ\}$$
is a generating set of $\Gamma/Q\cong\Z.$ For convenience, write
$$\langle T\rangle=\{a_1,\,\ldots,\,a_k\}\subset\Z\ \mbox{with}\ a_1\leq a_2\leq\ldots\leq a_k.$$
Due to $\Z$ is Abelian, the Cayley graphs on $\Z$ for $\langle T\rangle$ and $\langle T\rangle_{\pm}=\{\pm a_1,\,\ldots,\,\pm a_k\}$ are identical. Therefore, the image graph $\langle G\rangle$ of $G$ under the map
$$\pi:\ \gamma\in\Gamma\longmapsto \gamma Q\in \Gamma/Q\cong\Z$$
is a multi-graph on $\Z$ possibly with loop edges such that for any $x,y\in\Z,$ $x\sim y$ iff $x-y\in\langle T\rangle_{\pm}.$ With abusing notations, let $\left(\langle G\rangle,\,\mathbf{c}_1\right)$ and $\left(\langle G\rangle,\,\mathbf{c}_2\right)$ be the image networks of $(G,\,\mathbf{c}_1)$ and $(G,\,\mathbf{c}_2)$ respectively under $\pi,$ and
$$(\langle G\rangle,\,\mathbf{c}_1,\,\mathbf{c}_2,\, p)\ \mbox{the image network of}\ (G,\,\mathbf{c}_1,\,\mathbf{c}_2,\,p)\ \mbox{for any}\ p\in [0,1]\  \mbox{under}\ \pi.
$$

Write $\langle G\rangle_*$ for the multi-graph obtained from $\langle G\rangle$ by deleting all loop edges, $\mathbf{c}_1^*$ and $\mathbf{c}_2^*$ respectively for the restrictions of $\mathbf{c}_1$ and $\mathbf{c}_2$ to all non-loop edges of $\langle G\rangle.$
Then each $(\langle G\rangle_*,\,\mathbf{c}_1^*,\,\mathbf{c}_2^*,\, p)$ is just the restriction of $(\langle G\rangle,\,\mathbf{c}_1,\,\mathbf{c}_2,\, p)$ to $\langle G\rangle_*.$ And by Step 1 in proving Theorem \ref{recurthm}\,(iii), for any $p\in [0,\,1),$
\begin{eqnarray}\label{eq-recurrent-equiv-3}
(\langle G\rangle_*,\,\mathbf{c}_1^*,\,\mathbf{c}_2^*,\, p)=:(\langle G\rangle_*,\,\mathbf{c}_p^*)\ \mbox{is recurrent iff so is}\
  (\langle G\rangle,\,\mathbf{c}_1,\,\mathbf{c}_2,\, p)=:(\langle G\rangle,\,\mathbf{c}_p).
\end{eqnarray}

Notice that
$$(G,\,\mathbf{c}_1,\,\mathbf{c}_2,\,p)=:(G,\,\mathbf{c}_p)\longmapsto (\langle G\rangle,\,\mathbf{c}_1,\,\mathbf{c}_2,\, p)
  =\left(\langle G\rangle,\,\mathbf{c}_p\right)\ \mbox{is a rough embedding}.
$$
Similarly to Step 2 in proving Theorem \ref{recurthm}\,(iii), we have that if
\begin{eqnarray}\label{eq-recurrent-equiv-4}
\mbox{almost surely},\
\left(\langle G\rangle_*,\,\mathbf{c}_p^*\right)\ \mbox{is recurrent for all}\ p\in [0,1),
\end{eqnarray}
then
\begin{eqnarray*}
\mbox{almost surely},\ \mbox{for any}\ p\in [0,\,1),\ \left(G,\,\mathbf{c}_p\right)\ \mbox{is recurrent}.
\end{eqnarray*}
Therefore, it suffices to prove \eqref{eq-recurrent-equiv-4}.
\vskip 2mm

{\bf Step 2.} Clearly the generating set of $\langle G\rangle_*$ is the multiple set
$S=\{\pm a_i:\ a_i\not=0,\ 1\leq i\leq k\}.$ Removing multiplicities of elements in $S$, we get a set
$\widetilde{S}=\{\pm b_j:\ 0<b_1<\ldots<b_{\ell}\}.$

For any $x\in \widetilde{S},$ let
$$\mathbf{\tilde{c}}_p(\{0,x\})=\sum\limits_{e\in\pi^{-1}(\{0,x\})}\mathbf{c}_p(e),\ p\in [0,1),$$
and identify all parallel edges in $\pi^{-1}(\{0,x\})$ to the edge $\{0,x\}.$ Then we get the Cayley graph $\langle G\rangle_{\sim}$
on $\Z$ corresponding to $\widetilde{S},$ and a network process $\left((\langle G\rangle_{\sim},\,\mathbf{\tilde{c}}_p)\right)_{p\in [0,\,1)}.$
Clearly
\begin{eqnarray*}
\mbox{each}\ (\langle G\rangle_{\sim},\,\mathbf{\tilde{c}}_p)\ \mbox{is recurrent iff so is}\ \left(\langle
  G\rangle_*,\,\mathbf{c}_p^*\right).
\end{eqnarray*}
So \eqref{eq-recurrent-equiv-4} boils down to that
\begin{eqnarray}\label{eq-recurrent-equiv-5}
\mbox{almost surely},\ (\langle G\rangle_{\sim},\,\mathbf{\tilde{c}}_p)
\ \mbox{is recurrent for all}\ p\in [0,1).
\end{eqnarray}
\vskip 2mm

{\bf Step 3.} This step is to prove \eqref{eq-recurrent-equiv-5}.

Given any $p\in (0,\,1).$ For any $x\in\Z_+,$ let
$$A_x(p)=\left\{\forall\,y\in (x,\,x+b_{\ell}],\ \pi^{-1}(\{y,y+b_i\})\ \mbox{is a closed edge set in}\ \omega_p\ \mbox{for any}\ 1\leq i\leq\ell\right\};$$
and for any negative $x\in\Z,$ let
$$A_x(p)=\left\{\forall\,y\in [x-b_{\ell},\,x),\ \pi^{-1}(\{y,y-b_i\})\ \mbox{is a closed edge set in}\ \omega_p\ \mbox{for any}\ 1\leq i\leq\ell\right\}.$$
Clearly, $\{A_{kb_{\ell}}(p)\}_{k\in\Z}$ is an i.i.d.\,sequence of events and $\mathbb{P}[A_0(p)]>0.$ By the law of large numbers, almost surely,
$$\frac{1}{n+1}\sum\limits_{k=0}^n I_{A_{kb_{\ell}}(p)}\rightarrow \mathbb{P}[A_0(p)]>0,$$
and thus infinitely many $A_{kb_{\ell}}(p)$s, say $A_{k_jb_{\ell}}(p),\ j=1,\,2,\,\ldots,$ occur. Similarly, almost surely, infinitely many $A_{-kb_{\ell}}(p)$s ($k>0$), say $A_{-s_jb_{\ell}}(p),\ j=1,\,2,\,\ldots,$ occur. Fix a percolation environment $\omega=(\omega_q)_{q\in [0,\,1]}$ such that these events hold. Then
\begin{eqnarray*}
&&\Pi_n=\left\{\{y,y+b_i\}:\ 1\leq i\leq \ell,\ y\in (k_nb_{\ell},\,(k_n+1)b_{\ell}]\right\}\\
&&\ \ \ \ \ \ \ \ \ \ \cup \left\{\{y,y-b_i\}:\ 1\leq i\leq \ell,\ y\in [-(s_n+1)b_{\ell},\,-s_nb_{\ell})\right\}
\end{eqnarray*}
is a cutset of $\langle G\rangle_{\sim}$ separating $0$ from infinity such that $\pi^{-1}(e)$ is a closed edge set in $\omega_p$ (and all $\omega_q$ with $q\leq p$) for any $e\in\Pi_n.$ Notice that for any edge $e$ in $\Pi_n$ and $q\leq p,$
\begin{eqnarray*}
\mathbf{\tilde{c}}_q(e)&=&\sum\limits_{e^\prime\in\pi^{-1}(e)}\mathbf{c}_q(e^\prime)=\sum\limits_{e^\prime\in\pi^{-1}(e)}\left\{
   \mathbf{c}_1(e^\prime)I_{\{U_{e^\prime}\leq q\}}+ \mathbf{c}_2(e^\prime)I_{\{U_{e^\prime}>q\}}\right\}\\
   &=&\sum\limits_{e^\prime\in\pi^{-1}(e)}\mathbf{c}_2(e^\prime)\leq c\left\vert\pi^{-1}(e)\right\vert
     \leq c\max\limits_{e\in\Pi_n}\left\{\left\vert\pi^{-1}(e)\right\vert\right\}
     =c\max\limits_{e\in\Pi_1}\left\{\left\vert\pi^{-1}(e)\right\vert\right\}:=C.
\end{eqnarray*}
Thus effective resistance $\mathscr{R}_q(0\leftrightarrow\infty)$ between $0$ and infinity of
$(\langle G\rangle_{\sim},\,\mathbf{\tilde{c}}_q)$ with $q\leq p$ satisfies that
\begin{eqnarray*}
\mathscr{R}_q(0\leftrightarrow\infty)\geq \sum\limits_{n=1}^{\infty}\left(\sum\limits_{e\in\Pi_n}\mathbf{\tilde{c}}_q(e)\right)^{-1}
   \geq\sum\limits_{n=1}^{\infty}\left(\sum\limits_{e\in\Pi_n}C\right)^{-1}
   =\frac{1}{C}\sum\limits_{n=1}^{\infty}\frac{1}{\vert\Pi_n\vert}=\frac{1}{C}\sum\limits_{n=1}^{\infty}\frac{1}{\vert\Pi_1\vert}=\infty.
\end{eqnarray*}
By the Nash-Williams criterion, for the percolation environment $\omega,$ all networks $(\langle G\rangle_{\sim},\,\mathbf{\tilde{c}}_q)$
with $q\leq p$ are recurrent. Due to $p$ is arbitrary, we are done.
\qed
\vskip 2mm


The following Lemmas \ref{clustersize} and \ref{supdense} are necessary preliminaries for proving Theorem \ref{generald}. To state them, for any $n\in \N$, let
\[
R_n:= \{x\in \Z^d:\ \vert x\vert \leq n\},\
\partial R_n:= \{x\in \Z^d:\ \vert x\vert= n\}.
\]
Here $\vert \cdot \vert$ denotes $\ell_{1}$ norm (graph distance) on $\Z^d$. Recall $A\stackrel{\omega_p}{\longleftrightarrow}B$ denotes vertex sets $A$ and $B$ are connected to each other in $\omega_p.$
When $p$ is fixed, write $A{\longleftrightarrow}B$ for $A\stackrel{\omega_p}{\longleftrightarrow}B.$

\begin{lem}\label{clustersize}
{\bf (i)} {\rm \cite[Theorem 5.4]{GG1999}.} For Bernoulli-$p$ bond percolation $\omega_p$ on $\Z^d$ with $p< p_{c}$, there is a constant $\psi(p,d)>0$ such that for any $n,$
\[
\P_p\left(\mathbf{0}\longleftrightarrow \partial R_n\right)\leq e^{-n\psi(p,d)}.
\]

{\bf (ii)} {\rm \cite[Theorem 11.89]{GG1999}, \cite[Corollary 4.7]{HD2018}.} For critical Bernoulli bond percolation $\omega_{1/2}$ on $\Z^2$, there is a constant $\alpha>0$ such that
\[
\frac{1}{2n}\leq \P_{1/2}\left(\mathbf{0}\longleftrightarrow \partial R_n\right)\leq \frac{1}{n^{\alpha}},\ \forall\,n\geq 1.
\]

{\bf (iii)} {\rm \cite[Theorem 11.5]{HH2016}.} On $\Z^d$ with $d\geq 11,$ for critical Bernoulli bond percolation $\omega_{p_c},$
\[
\P_{p_c}\left(\mathbf{0}\longleftrightarrow \partial R_n\right)\asymp n^{-2}.
\]
\end{lem}
\vskip 2mm


The next lemma is a slight improvement of estimation in \cite[Appendix~\uppercase\expandafter{\romannumeral1}, Theorem~5]{FDTT2007}. Let $C_m$ be a box with size $m\times \kappa\log m$ whose boundaries are all edges in $\Z^2$,
where $\kappa\log m$ is understood as its integer part $\lfloor \kappa\log m\rfloor.$ Now consider Bernoulli bond percolation on $\Z^2$. A horizontal open crossing of $C_m$ is a simple path within $C_{m}$ connecting its left boundary with its right boundary.

\begin{lem}\label{supdense}
Consider Bernoulli bond percolation on $\Z^2.$ For any $p>1/2$ and $\kappa>\frac{2}{\psi\left(3/4-p/2\right)}$, there exists a $\delta(\kappa,p)> 0$ such that almost surely, the number of edge-disjoint horizontal open crossings of $C_m$ is no less than $\delta(\kappa,p)\log m$ when $m$ is large enough. Here $\psi(\cdot)$ is given in Lemma \ref{clustersize}\,(i).
\end{lem}
\pf Let ${\rm LR}(B_{m,n})$ be the event that there exists an open crossing from left boundary to right boundary of a box $B_{m,n}$ with size $m\times n$ in Bernoulli-$p$ bond percolation $\omega_p$ on $\Z^2.$ Define similarly the vertical crossing event ${\rm UD}(B_{m,n})$. Note that $C_m= B_{m,\kappa\log m}$. Write $B_n= B_{n,\, n}$. Recall a natural coupling of two bond percolations on $\Z^2$ and its dual lattice $(\Z^2)^{*}$: For any edge $e$ of $\Z^2$, write $e_*$ for its dual edge. Define the dual percolation $\omega_p^*$ on $(\Z^2)^{*}$ by
$$\omega_p^*(e_*)=1-\omega_p(e),\ \forall\,e\in E(\Z^2).$$
Clearly $\omega_p^*$ is a Bernoulli-$(1-p)$ bond percolation on $(\Z^2)^{*}$. Then the event that there is an open horizontal crossing in box $B_{m,n}$ in $\Z^2$ is equivalent to the event that there is no open vertical crossing in $B^{*}_{m,n}$ in $(\Z^2)^*$.

{\it In fact}, define the vertices in $B_{m,n}$ which is within or can be connected to left boundary by an open path to be `black vertices' and the other vertices to be `white vertices'. When ${\rm LR}(B_{m,n})$ does not occur, in $B^*_{m,n}$, we can find a unique interface which is a path separating `black vertices' from `white vertices' such that each type of vertices distribute in the same side of the interface, and the interface is an open vertical crossing in $B^*_{m,n}$, namely ${\rm UD}(B^{*}_{m,n})$ occurs. And clearly when ${\rm LR}(B_{m, n})$ occurs, ${\rm UD}(B^{*}_{m,n})$ in $(\Z^2)^*$ does not occur.

\begin{figure}[!htp]
\centering
\includegraphics[width=6.5cm, height= 6.5cm]{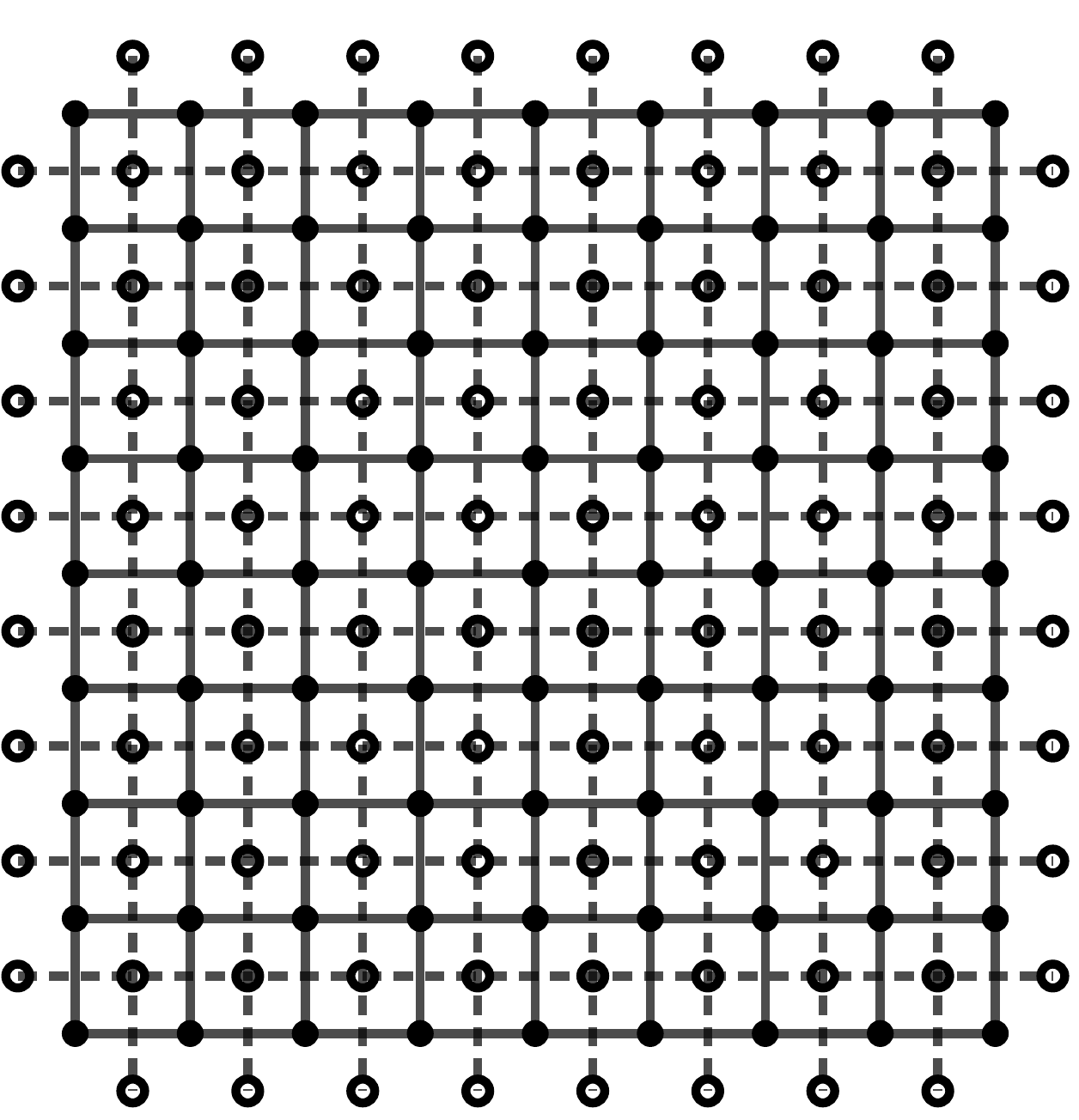}
\caption{$B_{n}$ is in the real line and filled dot; and $B_{n}^{*}$ is in the dashed line and hollow dot.}
\label{graphduallattice}
\end{figure}

Suppose $1/2< p_1< p$. By Lemma \ref{clustersize},
\begin{eqnarray*}
\P_{p_1}\left({\rm LR}(B_{m,\kappa\log m})\right)& = &1- \P_{1-p_1}\left({\rm UD}(B^*_{m,\kappa\log m})\right)\\
 &\geq & 1-(m+1)\P_{1-p_1}\left(\mathbf{0}\longleftrightarrow \partial R_{\lfloor\kappa\log m\rfloor}\right)\\
 &\geq & 1-(m+1) e^{-\psi(1-p_1)\lfloor\kappa\log m\rfloor }.
\end{eqnarray*}
On $\Z^2$, for any $\omega\in\{0,1\}^{E(\Z^2)},$ define
\[
E_r(\omega):=\left\{\omega'\in \{0,1\}^{E(\Z^2)}:\ \mbox{there are at most}\ r\ \mbox{edges}\ e\ \mbox{such that}\ \omega'(e)\neq \omega(e)\right\},
\]
and for any measurable set $A$, let
\[
I_r(A):= \left\{\omega\in \{0,1\}^{E(\Z^2)}:\ E_r(\omega)\subseteq A\right\}.
\]
Then by the max-flow min-cut theorem and the Menger's theorem, $I_{r}\left({\rm LR}(B_{m,\kappa\log m})\right)$ is the event that there are at least $r+1$ edge-disjoint open horizontal crossings from left boundary to right boundary of box $B_{m,\kappa\log m}$.
By \cite[Theorem~2.45]{GG1999}, for any $r\in\N,$
\[
1- \P_p\left(I_r\left({\rm LR}(B_{m,\kappa\log m})\right)\right)\leq \left(\frac{p}{p-p_1}\right)^r \left(1- \P_{p_1}\left({\rm LR}(B_{m,\kappa\log m})\right)\right).
\]
Let $p_1= p/2+ 1/4$. Notice $\kappa >\frac{2}{\psi\left(3/4- p/2\right)}$. We can choose $\delta=\delta(\kappa,p)>0$ small enough such that
\begin{eqnarray*}
\delta\log \frac{p}{p-p_1}+1-\kappa\psi(1- p_1)< -1.
\end{eqnarray*}
Therefore, as $m\rightarrow\infty,$
\begin{eqnarray*}
&&\P_{p}(\mbox{Number of edge-disjoint open horizontal crossings in}\ B_{m,\kappa \log m}\ \mbox{is at most}\ \delta\log m)\\
&&\leq\left(\frac{p}{p-p_1}\right)^{\lfloor\delta\log m\rfloor} (m+1) e^{-\psi(1-p_1)\lfloor\kappa\log m\rfloor }\\
&&\sim m^{\delta\log\frac{p}{p-p_1}+1-\kappa\psi(1-p_1)}.
\end{eqnarray*}
By the Borel-Cantelli lemma, we obtain the lemma immediately.
\qed
\vskip 2mm
%
\vskip 2mm
\noindent\textbf{Proof of Theorem \ref{generald}.} {\bf (a)} Verify $p_c^{*}= p_c$ for $0<\lambda_1<1<\lambda_2$ and $d\geq 2.$ We only prove $p_c^{*}= p_c$ for $d=2$ due to similarity for the case $d\geq 3.$

{\bf (a.1)} When $p>p_c=1/2$, almost surely, there is an infinite open cluster $\mathcal{C}$ on which all edges $e$ have conductances $\mathbf{C}_{\lambda_1}(e)$. Fix such a percolation configuration and choose a simple infinite path $\gamma=x_0x_1x_2\cdots$ in $\mathcal{C}.$
Define a unit flow $\theta$ on $\gamma$ such that
$$\theta(x_ix_{i+1})=1,\ \theta(x_{i+1}x_i)=-1,\ i\in\mathbb{Z}_{+}.$$
Then the energy $\mathscr{E}(\theta)$ of $\theta$ on network $(\gamma,\,\mathbf{C}_{\lambda_1})$ satisfies that
\begin{eqnarray*}
\mathscr{E}(\theta)=\sum\limits_{i\in\Z_{+}}\vert\theta(x_ix_{i+1})\vert^2\frac{1}{\mathbf{C}_{\lambda_1}(\{x_i,x_{i+1}\})}=\sum\limits_{i\in\Z_+}\lambda_1^{\vert x_i\vert\wedge\vert x_{i+1}\vert}
<\sum\limits_{e\in E(\Z^2)}\lambda_1^{\vert e\vert}<\infty.
\end{eqnarray*}
By Lemma \ref{lem-energy}, $(\gamma,\,\mathbf{C}_{\lambda_1})$ is transient. Then by the Rayleigh's monotonicity principle (Lemma \ref{lem-ray}), both $(\mathcal{C},\mathbf{C}_{\lambda_1})$ and $(\Z^2,\,\mathbf{C}_{\lambda_1},\,\mathbf{C}_{\lambda_2},\,p)$ are transient.

{\bf (a.2)} When $p< 1/2$, let $A_{n}$ be the event that $\partial R_{n^2}$ and $\partial R_{(n+1)^2}$ are connected by an open path for any $n\in \N$. Then with probability $1$, $\{A_n\}_n$ occurs for only finite times. {\it In fact}, by Lemma \ref{clustersize}, there exists $\psi(p)>0$ such that for a positive constant $C,$
\begin{eqnarray*}
\P_p\left(A_n\right)&\leq & \sum\limits_{x\in \partial R_{n^2}} \P_p\left(x\longleftrightarrow \partial R_{(n+1)^2}\right)\\
  &\leq & \sum\limits_{x\in \partial R_{n^2}} \P_p\left(\mathbf{0} \longleftrightarrow \partial R_{(n+1)^2- n^2}\right)\leq C n^2 e^{-\left(2n+1\right)\psi(p)},
\end{eqnarray*}
which verifies the claim by the Borel-Cantelli lemma. Thus, almost surely, there exists a random natural number $N$ such that for any $n\geq N$, there is a closed cutset $\Pi_n$ in $R_{(n+1)^2}\setminus R_{n^2-1}$ separating $\bf{0}$ and $\infty$; and further for some positive constant $C',$
\begin{eqnarray*}
\mathscr{R}(\bf{0}\leftrightarrow \infty)&\geq &\sum\limits_{n=N}^{\infty}\left(\sum\limits_{e\in \Pi_{n}} \mathbf{C}_{\lambda_2}(e)\right)^{-1}\geq C'\sum\limits_{n= N}^{\infty} \left( (n+1)^4 \lambda_2^{-n^2}\right)^{-1}\\
&=&C'\sum\limits_{n=N}^{\infty} \frac{\lambda_2^{n^2}}{(n+1)^4}= \infty.
\end{eqnarray*}
By the Nash-Williams recurrence criterion (Lemma \ref{lem-Nash}), $(\Z^2,\,\mathbf{C}_{\lambda_1},\,\mathbf{C}_{\lambda_2},\,p)$ is recurrent almost surely.
\vskip 2mm

{\bf (b)} Prove $p_c^*=p_c$ for $d\geq 3$ and $\lambda_1=1<\lambda_2.$

The recurrence for $p< p_c$ can be proved similarly to {\bf (a.2)}. For the supercritical case $p>p_c$, note that simple random walk on the infinite open cluster is transient almost surely (\cite[Theorem 1]{GHY1993}). Thus the transience for $p> p_c$ holds directly from this conclusion by the Rayleigh's monotonicity principle.
\vskip 2mm

{\bf (c)} Show that critical $\left(\Z^d,\,\mathbf{C}_{\lambda_1},\,\mathbf{C}_{\lambda_2},\,p_c\right)$ with $0<\lambda_1<1<\lambda_2$ and $d=2$ or $0<\lambda_1\leq 1<\lambda_2$ and $d\geq 11$ is recurrent almost surely.

{\bf (c.1)} Assume $d=2$. Then $p_c=1/2$, by Lemma \ref{clustersize}, for a positive constant $\alpha,$
\[
\frac{1}{2n}\leq \P_{1/2}(\mathbf{0}\longleftrightarrow \partial R_n)\leq \frac{1}{n^{\alpha}}.
\]
For any $n\in \N$, we take $A_n$ to be the event that there is no open path connecting $\partial R_{K^{l^n}}$ with $\partial R_{K^{l^{n+1}}}$ with $1<K\in\N$ and $1/\alpha< l\in\N$. Then
for some positive constant $C_1$,
\begin{eqnarray*}
\P_{1/2}\left(A_n\right)&\leq & \sum\limits_{x\in \partial R_{K^{l^n}}} \P_{1/2}\left(x\longleftrightarrow \partial R_{K^{l^{n+1}}}\right)\leq \sum\limits_{x\in \partial R_{K^{l^n}}} \P_{1/2}\left(\mathbf{0} \longleftrightarrow \partial R_{K^{l^{n+1}}- K^{l^{n}}}\right)\\
&\leq & C_1 K^{l^n} \left(K^{l^{n+1}}- K^{l^{n}}\right)^{-\alpha}\leq C_1 \left(\frac{K}{K-1}\right)^{\alpha} K^{l^n(1-\alpha l)}.
\end{eqnarray*}
Then similarly to {\bf (a.2)}, almost surely, there is a sequence $\{\Pi_n\}_{n=N}^{\infty}$ of disjoint closed cutsets with each $\Pi_n$ being in $R\left(K^{l^{n+1}}\right)\setminus R\left(K^{l^{n}}-1\right)$ and separating
$\bf{0}$ from $\infty;$ and for a positive constant $C_2,$
\begin{eqnarray*}
\mathscr{R}(\bf{0}\leftrightarrow\infty)&\geq &\sum\limits_{n=N}^{\infty}\left(\sum\limits_{e\in \Pi_{n}}\mathbf{C}_{\lambda_2}(e)\right)^{-1}\geq C_2 \sum\limits_{n=N}^{\infty} \left(K^{2 l^{n+1}} \lambda_2^{-K^{l^{n}}}\right)^{-1}\\
&=&C_2 \sum\limits_{n=N}^{\infty} \frac{\lambda_2^{K^{l^{n}}}}{K^{2 l^{n+1}}} = \infty,
\end{eqnarray*}
which shows $(\Z^2,\, \mathbf{C}_{\lambda_1},\, \mathbf{C}_{\lambda_2},\, 1/2)$ is recurrent almost surely.

{\bf (c.2)} Similarly to {\bf (c.1)}, one can prove that critical $\left(\Z^d,\,\mathbf{C}_{\lambda_1},\,\mathbf{C}_{\lambda_2},\,p_c\right)$ with $0<\lambda_1\leq 1<\lambda_2$ and $d\geq 11$ is recurrent almost surely.
\vskip 2mm

{\bf (d)} Prove $p_c^*=p_c=1/2$ for $d=2$ and $0< \lambda_1< 1=\lambda_2.$

The transience when $p>1/2$ can be proved similarly as in {\bf (a.1)}.

Assume $p\in [0,\, 1/2)$. we claim $(\Z^2,\, \mathbf{C}_{\lambda_1},\, \mathbf{C}_1,\, p)$ is recurrent almost surely. Simple random walk (SRW) on supercritical infinite open cluster in $\Z^2$ is recurrent by the Rayleigh's monotonicity principle.
Compared with SRW on supercritical infinite open cluster, to prove the recurrence of $(\Z^2,\, \mathbf{C}_{\lambda_1},\, \mathbf{C}_1,\, p)$ here is more interesting.
We prove this in two steps.
\vskip 2mm
\begin{figure}[!htp]
\centering
\includegraphics[width=7.5cm, height= 7cm]{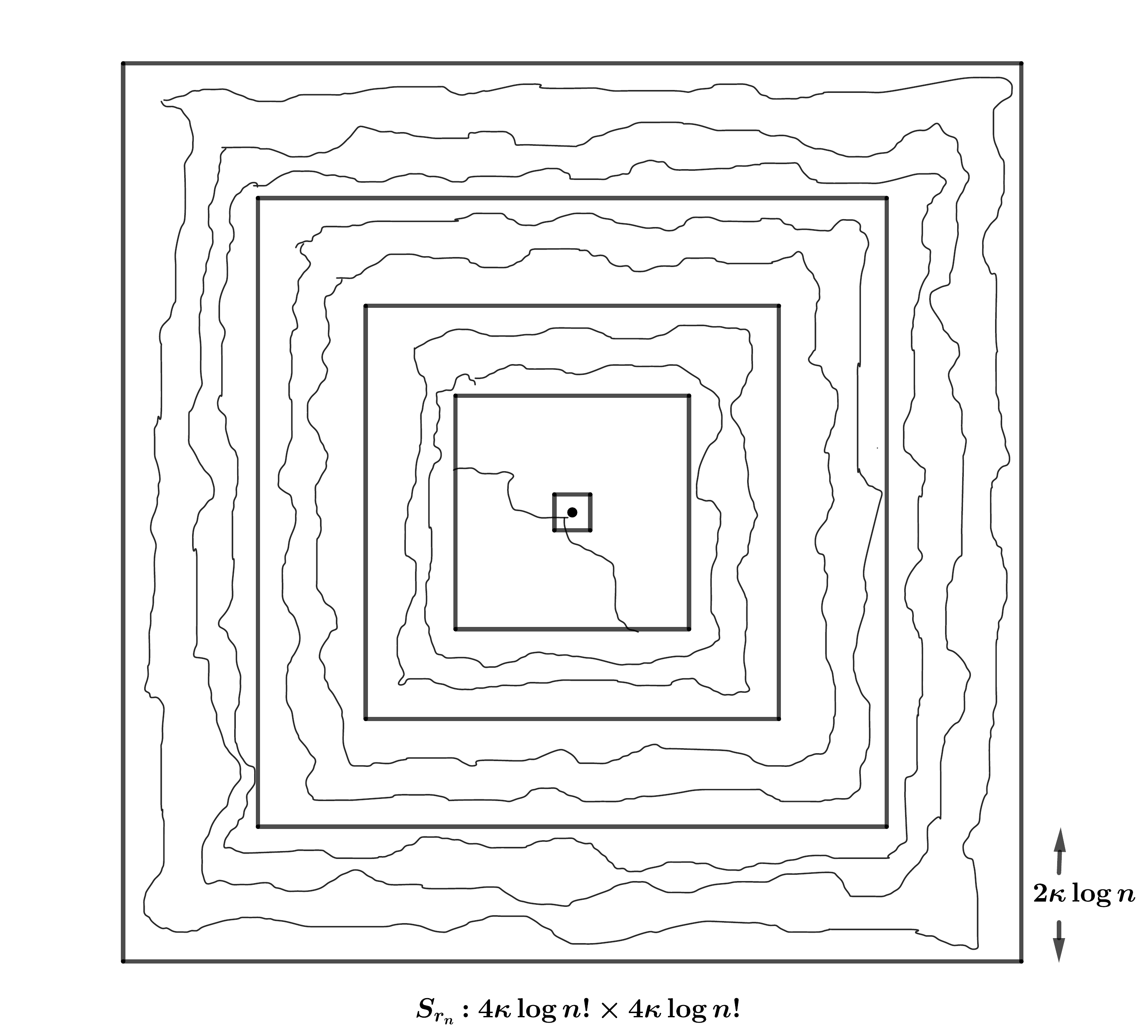}
\caption{Open circuits in annulus of box $S_{r_n}$}
\label{graphopencircuits}
\end{figure}
\textbf{Step (d.1). Preparation.} Let $S_{n}= [-n,\, n]\times [-n,\, n]$ be a square in $\Z^2$ for any $n\in (0,\,\infty).$ Fix $p_0\in (1/2,\,1)$ and $\kappa>2/\psi(3/4- p_0/2)$.
We aim to prove that there exists a constant $\delta=\delta(\kappa,p_0)>0$ such that when $n$ is large enough, number of edge-disjoint open circuits in the annuli $S_{r_{n+1}}\backslash S_{r_n}$ (c.f. Figure \ref{graphopencircuits}) with $r_n=  2\kappa\log{n!}$ is no less than $\delta \log(4 \kappa\log{(n+1)!})$ almost surely.

{\it In fact}, $S_{r_{n+1}}\backslash S_{r_{n}}$ is the union of
\begin{eqnarray*}
&&A_n^u:=[-2\kappa \log{(n+1)!},\, 2\kappa \log{(n+1)!}]\times [2\kappa \log{n!},\, 2\kappa \log{(n+1)!}],\\
&&A_n^d:=[-2\kappa \log{(n+1)!},\, 2\kappa \log{(n+1)!}]\times [-2\kappa \log{(n+1)!},\, -2\kappa \log{n!}],\\
&&A_n^l:=[-2\kappa \log{(n+1)!},\, -2\kappa \log{n!}]\times [-2\kappa \log{(n+1)!},\, 2\kappa \log{(n+1)!}],\\
&& A_n^r:=[2\kappa \log{n!},\, 2\kappa \log{(n+1)!}]\times [-2\kappa \log{(n+1)!},\, 2\kappa \log{(n+1)!}].
\end{eqnarray*}
Define
\begin{eqnarray*}
&&C_n^{u}=\{\mbox{There are at least}\ \delta \log(4 \kappa\log{(n+1)!})\ \mbox{edge-disjoint open horizontal crossings of}\ A_n^{u}\},\\
&&C_n^{d}=\{\mbox{There are at least}\ \delta \log(4 \kappa\log{(n+1)!})\ \mbox{edge-disjoint open horizontal crossings of}\ A_n^{d}\};\\
&&C_n^{l}=\{\mbox{There are at least}\ \delta \log(4 \kappa\log{(n+1)!})\  \mbox{edge-disjoint open vertical crossings of}\ A_n^{l}\},\\
&&C_n^{r}=\{\mbox{There are at least}\ \delta \log(4 \kappa\log{(n+1)!})\ \mbox{edge-disjoint open vertical crossings of}\ A_n^{r}\}.
\end{eqnarray*}
Note that as $n\rightarrow \infty$,
\[
2\kappa\log (n+1)!- 2\kappa\log n!> \kappa \log(4\kappa\log (n+1)!).
\]
Then similarly to Lemma \ref{supdense}, we can prove that for some constant $\delta=\delta(\kappa,p_0)>0,$ almost surely, $C_n^u$, $C_n^d$, $C_n^l$
and $C_n^r$ occur when $n$ is large enough.

Thus, by taking intersections of vertical crossings in $A_n^l$ and $A_n^r$ with horizontal crossings in $A_n^u$ and $A_n^d$, we have that almost surely,
number of edge-disjoint open circuits in the annuli $S_{r_{n+1}}\backslash S_{r_n}$ is no less than $\delta \log(4 \kappa\log{(n+1)!})$ for large enough $n$.
\vskip 2mm

\textbf{Step (d.2). Completing proof.} On $(\Z^2)^*$, each edge is open independently with probability $1-p$ by the coupling introduced in proof of Lemma \ref{supdense}.
Recall an edge $e_*$ in $(\Z^2)^*$ is open iff its dual edge $e$ in $\Z^2$ is closed. Note that the set of dual edges of an open circuit in $(\Z^2)^*$ is a closed cutset separating $\bf{0}$ from $\infty$ in $\Z^2$, and
edges $e$ in this cutset are of conductances $\mathbf{C}_{2}(e)=1$.

Choose $n_0\in\N$ such that
$$\delta\log(4\kappa\log(n_0+1)!)\geq 2.$$
For any $n_0\leq n\in \N,$ on $(\Z^2)^*$, let $c_n^1,\, c_n^2,\, \cdots,\, c_n^{m(n)}$ be edge-disjoint open concentric circuits in $S_{r_{n+1}}\backslash S_{r_n}$ with $c_n^i$ being in the interior of $c_n^{i+1},\, 1\leq i\leq m(n)-1$. And $\pi_n^1,\, \pi_n^2,\, \cdots,\, \pi_n^{m(n)}$ denote dual edge sets of $c_n^1,\ c_n^2,\, \cdots,\, c_n^{m(n)}$ which is a sequence of disjoint cutsets  on $\Z^2$. By the Nash-Williams inequality (Lemma \ref{lem-Nash}) on $\Z^2$ and Step {\bf (d.1)}, almost surely,
\begin{eqnarray*}
\mathscr{R}(\bf{0}\leftrightarrow\infty)&\geq &\sum\limits_{n=n_0}^{\infty} \sum\limits_{i=1}^{m(n)-1} |\pi_n^{i}|^{-1}=\sum\limits_{n=n_0}^{\infty} \sum\limits_{i=1}^{m(n)-1} |c_n^{i}|^{-1}\\
&\geq &\sum\limits_{n=n_0}^{\infty} \sum\limits_{i=1}^{\lfloor\delta \log(4\kappa\log{(n+1)!})\rfloor- 1}  |c_n^{i}|^{-1}.
\end{eqnarray*}
Here $|c_n^i|$ denotes the number of edges in $c_n^i$. Note that
\[
\frac{(2\kappa\log (n+1)!)^2- (2\kappa\log n!)^2}{\{\delta \log(4\kappa\log{(n+1)!})\}^2}\asymp n,
\]
and for some positive constant $C_3,$
\[|c_n^1|+ |c_n^2|+ |c_n^3|+ \cdots+ |c_n^{m(n)}| \leq  C_3\left\{(2\kappa\log (n+1)!)^2- (2\kappa\log n!)^2\right\},
\]
and by the Jensen's inequality, with $a_n=\lfloor\delta \log(4\kappa\log{(n+1)!})\rfloor- 1\leq m(n)-1,$
\begin{eqnarray*}
\sum\limits_{i=1}^{a_n}  |c_n^{i}|^{-1}\geq a_n\left\{\frac{1}{a_n}\sum\limits_{i=1}^{a_n}  |c_n^{i}|\right\}^{-1}\geq a_n^2\left\{\sum\limits_{i=1}^{m(n)}  |c_n^{i}|\right\}^{-1}.
\end{eqnarray*}
Thus, for some positive constant $C_4,$  almost surely,
\begin{eqnarray*}
\mathscr{R}(\bf{0}\leftrightarrow \infty)&\geq & \sum\limits_{n=n_0}^{\infty}\frac{\{\lfloor\delta \log(4\kappa\log{(n+1)!})\rfloor-1\}^2}{C_3\left\{(2\kappa\log (n+1)!)^2- (2\kappa\log n!)^2\right\}}\\
&\geq & C_4\sum\limits_{n=n_0}^{\infty} \frac{1}{n}=\infty,
\end{eqnarray*}
which implies that $(\Z^2,\, \mathbf{C}_{\lambda_1},\, \mathbf{C}_1,\, p)$ is recurrent almost surely.\qed\\

\noindent\textbf{Proof of Theorem \ref{thm-current-unique}.} {\bf Step 1.} Some preliminaries. Any recurrent infinite network $(G,\,\mathbf{c})$ has unique currents. So by the Rayleigh's monotonicity principle, Theorem \ref{thm-current-unique} is trivial for $1\leq \lambda_1<\lambda_2.$ Recall the following two criteria on current uniqueness from \cite[Chapter 9]{LP2017}.

Consider infinite network $(G,\, \mathbf{c})$. Note $\mathbf{r}=1/\mathbf{c}.$ Let $\ell^2\left(\overrightarrow{E},\,\mathbf{r}\right)$ be the space of anti-symmetric functions $\theta$ satisfying $\sum\limits_{e\in \overrightarrow{E}} \theta^2(e)\mathbf{r}(e)<\infty$. For any function $f$ on $V$, define its gradient as follows:
\[
\nabla f(e)=\mathbf{c}(\{x,y\}) {\rm d}f (e)=\mathbf{c}(\{x,y\}) \left(f(x)- f(y)\right),\ e=xy\in \overrightarrow{E}.
\]
Call ${\bf D}=\left\{f:\ \nabla f\in \ell^2\left(\overrightarrow{E},\,\mathbf{r}\right)\right\}$ the Dirichlet space.
Write {\bf HD} for the harmonic Dirichlet space consisting of harmonic functions in {\bf D}. Then $(G,\,\mathbf{c})$ is current unique iff
${\bf HD}=\mathbb{R}.$

For any finite subgraph $A=(V(A),\,E(A))$ of $G$, define
\[
{\rm RD}(A)= \sup\left\{\mathscr{R}\left(x\leftrightarrow y;\ A\right): \ x,\, y\in V(A)\right\},
\]
where $\mathscr{R}(x\leftrightarrow y;\ A)$ is the effective resistance between $x$ and $y$ in finite network $(A,\,\mathbf{c})$.
Then $(G,\,\mathbf{c})$ has unique currents if there is a sequence $\{W_n\}_{n=1}^{\infty}$ of pairwise edge-disjoint finite subnetworks of $(G,\,\mathbf{c})$ such that each (equivalently, some) vertex is separated from $\infty$ by all but finitely many $W_n$ and
\begin{equation}\label{currentu}
\sum\limits_{n=1}^{\infty}\frac{1}{{\rm RD}(W_n)}= \infty.
\end{equation}
\vskip 2mm
\noindent{\bf Step 2.} Assume $0< \lambda_1<\lambda_2\leq 1$. For any $n\geq 1,$ let $W_n$ be the set of edges $e$ of $\Z^2$ with $\vert e\vert=n.$ View naturally each $W_n$ as a subnetwork $W_n(p)$ of any given $(\Z^2,\,\mathbf{C}_{\lambda_1},\,\mathbf{C}_{\lambda_2},\,p).$
Let $\mathbf{c}_p$ (resp.\,$\mathbf{r}_p$) be the conductance (resp. resistance) function of $(\Z^2,\,\mathbf{C}_{\lambda_1},\,\mathbf{C}_{\lambda_2},\,p).$ Denote by  $\mathcal{P}_n(x,y)$ any edge-disjoint path in $W_n$ starting from $x$ and ending at $y.$ Clearly the length of $\mathcal{P}_n(x,y)$ is at most $4(n+1).$ Then
\begin{eqnarray*}
{\rm RD}(W_n(p))&=&\sup\left\{\mathscr{R}(x\leftrightarrow y;\ W_n(p)):\ x,y\in V(W_n)\right\}\\
&\leq &\sup\left\{\sum\limits_{e\in \mathcal{P}_n(x,y)}\mathbf{r}_p(e):\ x,y\in V(W_n)\right\}
     \leq 4(n+1)(\lambda_1 \vee \lambda_2)^{n}\leq 4(n+1),
\end{eqnarray*}
and thus
\begin{equation*}
\sum\limits_{n=1}^{\infty}\frac{1}{{\rm RD}(W_n(p))}= \infty,\ p\in [0,\,1].
\end{equation*}
Therefore, almost surely, for all $p\in [0,\,1],$ $(\Z^2,\,\mathbf{C}_{\lambda_1},\,\mathbf{C}_{\lambda_2},\,p)$ is current unique.
\vskip 2mm

\noindent{\bf Step 3.} Suppose $0<\lambda_1<1<\lambda_2.$ Note that almost surely,
for any $p\in [0,\,1/2],$ $(\Z^2,\,\mathbf{C}_{\lambda_1},\,\mathbf{C}_{\lambda_2},\,p)$ is recurrent.
It suffices to prove the theorem for all $p>1/2$ by proving that
$$\mbox{almost surely},\ \mathbf{HD}(p)= \R,\ \forall\, p\in (1/2,1].$$
Here $\mathbf{HD}(p)$ is the harmonic Dirichlet space of $(\Z^2,\,\mathbf{C}_{\lambda_1},\,\mathbf{C}_{\lambda_2},\,p).$

{\it In fact}, for any fixed $p>1/2,$ almost surely, there exist infinitely many open circuits $\{C_n\}_{n=1}^{\infty}$ in $\omega_p$ and thus in all $\omega_q$ with $q\geq p$ such that each $[-n,\, n]^2$ is within the interior of $C_n,$
and $C_n$ is within the interior of $C_{n+1}.$ Fix such a percolation environment $\omega.$ Clearly on $C_n$,
$$\mathbf{c}_q(e)=\lambda_1^{-|e|},\ e\in C_n,\ q\geq p. $$
For any $f\in \mathbf{HD}(q)$ with $q\geq p$, suppose $x_1^n, x_2^n$ are respectively the maximal and minimal points of $f$ in the area \[B(C_n)=\{x\in \Z^2:\ x\ \mbox{is on or within}\ C_n\}\]
whose boundary is $C_n$. Namely,
\begin{eqnarray*}
f(x_1^n)= \max\{f(y):\ y\in B(C_n)\},\ f(x_2^n)= \min\{f(y):\ y\in B(C_n)\}.
\end{eqnarray*}
Note $x_1^n,x_2^n\in C_n$ since $f$ is harmonic on $B(C_n)$. Let $\mathcal{P}_n(x_1^n,x_2^n)$ be the path in $C_n$ connecting $x_1^n$ and $x_2^n.$ Then as $n\rightarrow\infty,$
\begin{eqnarray*}
\left|f(x_1^n)- f(x_2^n)\right|&\leq &\sum\limits_{\{u,v\}\in \mathcal{P}_n(x_1^n,x_2^n)} |f(u)- f(v)|\\
&\leq & \sum\limits_{\{u,v\}\in \mathcal{P}_n(x_1^n,x_2^n)}\mathbf{r}_q(\{u,v\})
\sum\limits_{\{u,v\}\in \mathcal{P}_n(x_1^n,x_2^n)} \mathbf{c}_q(\{u,v\}) |f(u)- f(v)|^2\\
&\leq & \sum\limits_{\{u,v\}\in \mathcal{P}_n(x_1^n,x_2^n)}\lambda_1^{\vert u\vert\wedge\vert v\vert}
  \left(\sum\limits_{\vert z-w\vert=1}\mathbf{c}_q(\{z,w\})\vert f(z)-f(w)\vert^2\right)\\
&\leq & \sum\limits_{\stackrel{\vert u-v\vert=1}{\vert u\vert,\vert v\vert\geq n}}\lambda_1^{\vert u\vert\wedge\vert v\vert}
  \left(\sum\limits_{\vert z-w\vert=1}\mathbf{c}_q(\{z,w\})\vert f(z)-f(w)\vert^2\right)\\
&\rightarrow & 0.
\end{eqnarray*}
Namely as $n\rightarrow\infty,$
$$\max\limits_{x,y\in B(C_n)}\vert f(x)-f(y)\vert\rightarrow 0.$$
Therefore, $f$ is constant, and further ${\bf HD}(q)=\R.$ So far we have proved that for any fixed $p>1/2,$ almost surely,
for all $q\geq p,\ {\bf HD}(q)=\R;$
which completes proving the theorem.
\qed
\vskip 2mm

\subsection{Proof of Theorem \ref{regulartree1}}\label{rrt}
\noindent

\begin{lem}[Heathcote, Seneta and Vere-Jones (1967) \cite{HSV1967}]\label{heightgw}
For a Galton-Watson tree $T$, let $T_n$ be the number of vertices in $n$-th generation from the origin $o$ and $L$ the random number of offsprings of $o$. Assume $m:= \E[L]\in (0,\,\infty)$. Then $\{\P(T_n> 0)/ m^{n}\}$ is a decreasing sequence. When $m< 1$, the following are equivalent:
\begin{eqnarray*}
{\bf (i)}\ \lim\limits_{n\rightarrow \infty} \P[T_n> 0]/ m^{n} >0;\ {\bf (ii)}\ \sup \E[\left. T_n\ \right\vert\ T_n> 0]< \infty;\ {\bf (iii)}\ \E[L \log^{+} L]< \infty.
\end{eqnarray*}
\end{lem}

Before proving Theorem \ref{regulartree1}, we calculate the growth rate of infinite open cluster of supercritical percolation on $\mathbb{T}^d$ with $d\geq 3$ and root $o.$
Let $\mathcal{T}_{\omega_p}$ be the open cluster of root $o$ in Bernoulli-$p$ bond percolation $\omega_p$ on $\mathbb{T}^d$. Write $B(n,p)$ for the binomial distribution with trial number $n\in \N$ and success probability $p\in [0,\,1].$
Then $\mathcal{T}_{\omega_p}$ is a Galton-Watson tree with a slight modification, namely offspring number of root $o$ obeys $B(d,p)$ and the other vertices in $\mathcal{T}_{\omega_p}$ have an i.i.d.\,$B(d-1,p)$ family of offsprings. So many limit properties for Galton-Watson trees are still applicable for $\mathcal{T}_{\omega_p}$.

Let $|T_n|$ be the number of edges at distance $n$ from $o$ in $\mathcal{T}_{\omega_p}$. Similarly to \cite[Proposition~5.5]{LP2017}, one can prove that
$(|T_n|/[d(d-1)^{n-1}p^n])_{n\geq 1}$ is a non-negative martingale, and has a finite limit $\mathscr{W}$ a.s. by the martingale convergence theorem. Similarly to the Kesten-Stigum Theorem (\cite[Section 12.2]{LP2017}), one
can show that when $(d-1)p>1$ namely $p>\frac{1}{d-1}$,
\[
\P(\mathscr{W}=0)=\P[\vert \mathcal{T}_{\omega_p}\vert <\infty]=:q.
\]
Thus for $p>\frac{1}{d-1}$, on the event $\left\{\vert \mathcal{T}_{\omega_p}\vert =\infty\right\}$, almost surely,
\begin{eqnarray*}
\lim\limits_{n\rightarrow \infty} \frac{|T_n|}{d(d-1)^{n-1}p^n} =\mathscr{W}> 0\ \mbox{and}\
{\rm gr}(\mathcal{T}_{\omega_p})=\lim\limits_{n\rightarrow \infty} {\sqrt[n]{|T_n|}} = (d-1)p.
\end{eqnarray*}

Note the percolation process $(\omega_p)_{p\in [0,\,1]}$ on $\mathbb{T}^d$ is constructed by the grand coupling. By the indistinguishability of infinite open clusters of Bernoulli percolation on $\mathbb{T}^d$ (see Subsection \ref{sec-discussion-growth}), almost surely, as $\frac{1}{d-1}<p\uparrow 1$, each infinite open cluster of $\omega_p$ has growth rate $(d-1)p\uparrow d-1={\rm gr}(\mathbb{T}^d)$.\\

\noindent{\bf Proof of Theorem \ref{regulartree1}.} {\bf (i)} When $0<\lambda_1\leq 1<d-1\leq\lambda_2$ and $d\geq 3,$ $p_c^*=p_c=\frac{1}{d-1}.$

{\bf (i.1)} For any $p>\frac{1}{d-1}$, by \eqref{eq-p_c-br-tree} and Lemma \ref{gpercolation}, condition on $\left\{\vert \mathcal{T}_{\omega_p}\vert =\infty\right\},$ almost surely,
\[({\rm br}(\mathcal{T}_{\omega_p}))^{-1}= p_{c}(\mathcal{T}_{\omega_p})=\frac{1}{(d-1)p},\ \mbox{and thus}\ {\rm br}(T_{\omega_p})= (d-1)p={\rm gr}(\mathcal{T}_{\omega_p}).\]
By \cite[Theorem~3.5]{LP2017}, when $p>\frac{\lambda_1\vee 1}{d-1}=\frac{1}{d-1}$, condition on $\left\{\vert \mathcal{T}_{\omega_p}\vert =\infty\right\},$ almost surely, $\RW_{\lambda_1}$ constrained on $\mathcal{T}_{\omega_p}$ is transient. So by the Rayleigh's monotonicity principle and $0$-$1$ law for recurrence/transience (see {\bf (i)} in proving Theorem \ref{generalgraph01}), for any $p>\frac{1}{d-1},$
$(\mathbb{T}^d,\,\mathbf{C}_{\lambda_1},\, \mathbf{C}_{\lambda_2},\,p)$ is transient almost surely.

{\bf (i.2)} Assume $p\leq\frac{1}{d-1}$. For any $n\in\N,$ let $L_n$ be the edge set in which all edges are at distance $n$ from the root $o$ on $\mathbb{T}^d$. Note $p_c=\frac{1}{d-1}$ and critical percolation $\omega_{p_c}$ on $\mathbb{T}^d$ does not percolate almost surely (\cite[Theorem 8.21]{LP2017}).
Define a random sequence $\{X_n\}_{n=0}^{\infty}$ as follows: Let $X_0= 0$ and $L_{X_0}= \{o\}$. For any $i\in\N,$ define $X_i$ to be the smallest natural number $n> X_{i-1}$ such that $L_{X_{i-1}}$ cannot connect to $L_{n}$ in $\omega_p.$ Here $L_{i}$ and $L_{j}$ $(j> i)$ are connected if there exists an open path in $\omega_p$ with its two endpoints belonging to $L_{i}$ and $L_{j}$ respectively.
Clearly, almost surely, $\{X_n\}_{n=0}^{\infty}\subseteq\Z_{+}$ is a strictly increasing sequence.

For any subtree $T\subset \mathbb{T}^d$, its outer boundary is defined by
\[
\partial T:= \left\{e\in E(\mathbb{T}^d):\ e=\{x,y\},\ x\in T,\ y\notin T,\ \vert y\vert>\vert x\vert\right\}.
\]
For any $x\in\mathbb{T}^d,$ let $\mathcal{T}(x)$ be the open subtree of $\mathbb{T}^d$ consisting of $x$ and its offsprings connected to $x$ by an open path in $\omega_p$.
Clearly, $\mathcal{T}(x)$ is a Galton-Watson tree with root $x$ and offspring distribution $B(d-1,p).$ Let
$$W(x):=\sum_{e\in \partial \mathcal{T}(x)} \mathbf{C}_{\lambda_2}(e).$$
And code the individuals of $n$-th generation in $\mathbb{T}^d$ by
\[x^{n}_1,\, x^{n}_2,\, \cdots,\, x^{n}_{d(d-1)^{n-1}}.\]
Then we claim that when $p<\frac{1}{d-1}$,
\begin{eqnarray}\label{eq-claim}
\P\left[\sum\limits_{n=1}^{\infty} \left(W\left(x^{X_{n}}_1\right)+ W\left(x^{X_{n}}_2\right)+\cdots +W\left(x^{X_{n}}_{d(d-1)^{X_{n}-1}}\right)\right)^{-1}=\infty\right]=1.
\end{eqnarray}


\begin{figure}[!htp]
\centering
\includegraphics[width=13cm, height= 9cm]{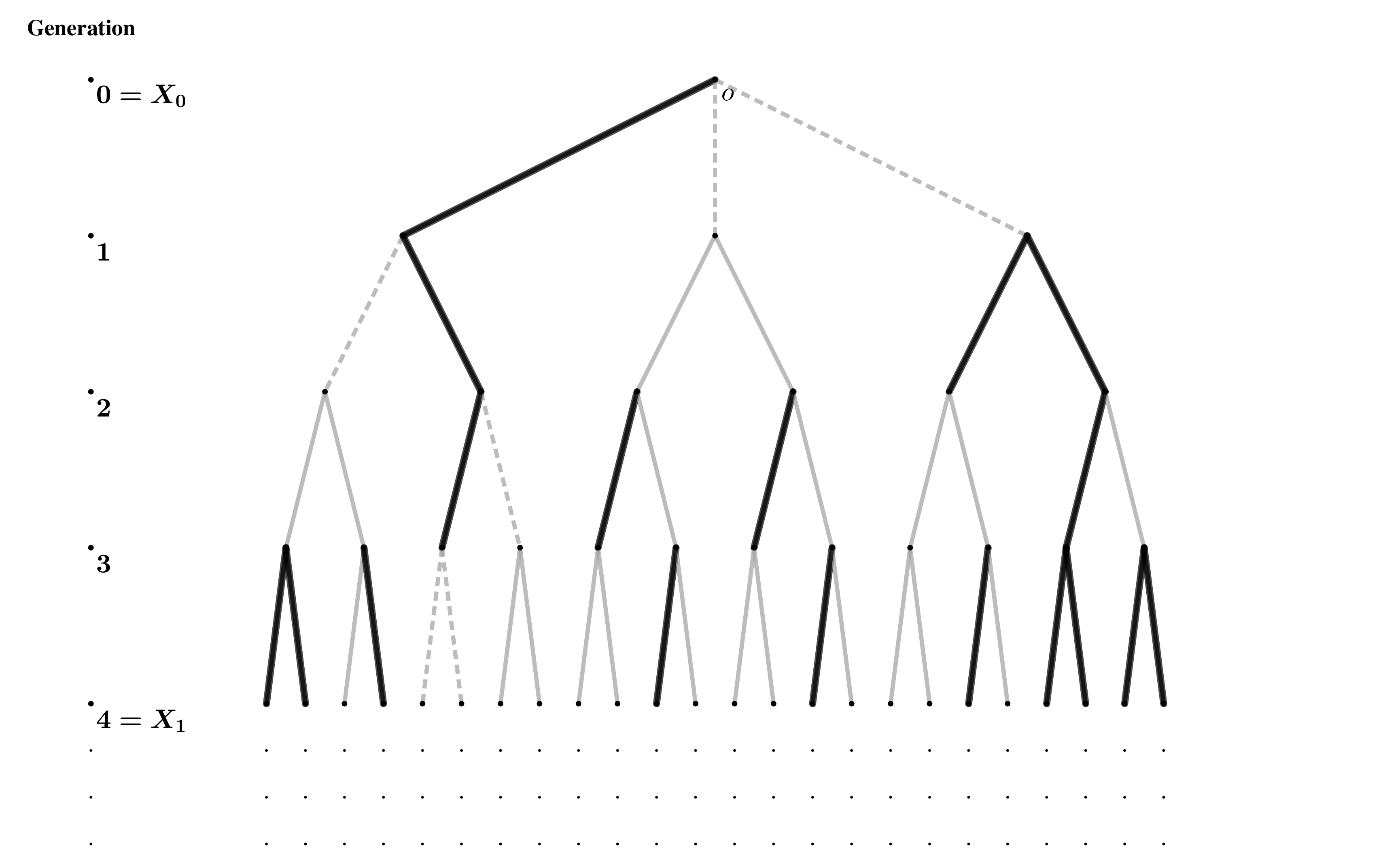}
\caption{Bernoulli percolation ($p= 0.4$) on binary tree: open edges are in black line; closed edges are in grey line; $\partial \mathcal{T}(o)$ are in dash line.}
\label{graphregulartree}
\end{figure}

To see this, for $j\in\N,$ let $R_j= d(d-1)^{X_j- 1},$ and $Y_{i}^{j}= \lambda_2^{X_j} W\left(x_{i}^{X_j}\right),\ 1\leq i\leq R_j,$ and
$$S_j:=\frac{\sum\limits_{i= 1}^{R_j} \left[Y_{i}^{j}- \E\left(Y_{1}^{1}\right)\right]}{R_j}.$$
By Lemma \ref{heightgw}, we have that when $p<\frac{1}{d-1}$,
\begin{eqnarray}\label{eq-2nd-moment}
\E\left[(Y_1^1)^2\right]= \lambda_2^2 \E\left[W^2\left(x_1^{1}\right)\right]\leq \lambda_2^2\sum\limits_{k=0}^{\infty} p_k \left(\sum\limits_{i=0}^{k} \lambda_2^{-i} d (d-1)^{i}\right)^2=:C < \infty,
\end{eqnarray}
where $p_k= \P\left[\mbox{The height of}\ \mathcal{T}\left(x_1^{1}\right)\ \mbox{is}\ k\right]$. Then when $p<\frac{1}{d-1}$, for any $j\in\N,$
\begin{eqnarray*}
\E\left[S_j^2\right]= \E\left[\E\left[\left. S_j^2\ \right\vert\ R_j\right]\right]= \E\left[\frac{\E[(Y_{i}^{j}- \E(Y_{1}^{1}))^2]}{R_j}\right]\leq \frac{C}{d(d-1)^{j-1}},
\end{eqnarray*}
and by the Chebyshev inequality, for any $\epsilon>0,$
\begin{eqnarray*}
\P[|S_j|> \epsilon]\leq \frac{\E\left[S_j^2\right]}{\epsilon^2}\leq \frac{C}{\epsilon^2 d(d-1)^{j-1}};
\end{eqnarray*}
which leads to
\[
\lim\limits_{j\rightarrow \infty} S_j= 0 \ \ a.s.
\]
by the Borel-Cantelli lemma. Namely, when $p<\frac{1}{d-1}$,
\[
\frac{W\left(x^{X_n}_1\right)+ W\left(x^{X_n}_2\right)+\cdots + W\left(x^{X_n}_{d(d-1)^{X_n-1}}\right)}{d(d-1)^{X_n-1} \lambda_2^{- X_n}}\xrightarrow[n\rightarrow \infty]{a.s.} \lambda_2 \E\left[W\left(x^{1}_1\right)\right]
\in (0,\,\infty).
\]
Since $\P\left[\sum\limits_{n=1}^{\infty} d(d-1)^{-(X_{n}-1)} \lambda_2^{X_{n}}=\infty\right]=1$, we obtain \eqref{eq-claim}.

Therefore, when $p< \frac{1}{d-1}$, by the Nash-Williams inequality and recurrence criterion (Lemma \ref{lem-Nash}), almost surely,
\begin{align*}
\mathscr{R}\left(o\leftrightarrow \infty\right)\geq \sum\limits_{j=1}^{\infty} \left(W\left(x^{X_{n_j}}_1\right)+ W\left(x^{X_{n_j}}_2\right)+\cdots +W\left(x^{X_{n_j}}_{d(d-1)^{X_{n_j}-1}}\right)\right)^{-1}
  =\infty,
\end{align*}
and $(\mathbb{T}^d,\, \mathbf{C}_{\lambda_1},\, \mathbf{C}_{\lambda_2},\, p)$ is recurrent almost surely.
\vskip 2mm

{\bf (ii)} Critical $\left(\mathbb{T}^d,\,\mathbf{C}_{\lambda_1},\,\mathbf{C}_{\lambda_2},\,\frac{1}{d-1}\right)$ with $0<\lambda_1\leq 1<d-1<\lambda_2$ and $d\geq 3$ is recurrent almost surely.

{\it In fact}, note that \eqref{eq-2nd-moment} still holds in this case due to $\lambda_2>d-1.$ Similarly to {\bf (i.2)}, we can prove that almost surely $\left(\mathbb{T}^d,\,\mathbf{C}_{\lambda_1},\,\mathbf{C}_{\lambda_2},\,\frac{1}{d-1}\right)$ is recurrent.

\vskip 2mm

{\bf (iii)} Suppose $1<\lambda_1<d-1\leq\lambda_2$ and $d\geq 3.$ Then $p_c^*=\lambda_1p_c=\frac{\lambda_1}{d-1}.$

{\bf (iii.1)} For $p>\frac{\lambda_1}{d-1},$ almost sure transience of $\left(\mathbb{T}^d,\,\mathbf{C}_{\lambda_1},\,\mathbf{C}_{\lambda_2},\, p\right)$ can be proved similarly to {\bf (i.1)}.

{\bf (iii.2)} When $p< \frac{1}{d-1}$, similarly to {\bf (i.2)}, one can verify that $\left(\mathbb{T}^d,\,\mathbf{C}_{\lambda_1},\,\mathbf{C}_{\lambda_2},\, p\right)$ is recurrent almost surely.
To prove that when $\frac{1}{d-1}< p< \frac{\lambda_1}{d-1}$, $\left(\mathbb{T}^d,\,\mathbf{C}_{\lambda_1},\,\mathbf{C}_{\lambda_2},\, p\right)$ is recurrent almost surely. Note by the Rayleigh's monotonicity principle,
this implies that $\left(\mathbb{T}^d,\,\mathbf{C}_{\lambda_1},\,\mathbf{C}_{\lambda_2},\, q\right)$ is a.s.\,recurrent for any $q<\frac{\lambda_1}{d-1}.$

{\bf Assume firstly $\lambda_2> d-1$.} Define a sequence $\{a_n\}_{n=0}^{\infty}\subset\Z_+$ by letting
$$a_0= 0,\ a_1= d-1,\ \mbox{and}\ a_{i+1}=(d-1)^{a_i},\ i\in\mathbb{N}.$$
Then define a sequence of cutsets $\{\Pi_i\}_{i=1}^{\infty}$ in $\mathbb{T}^d$ as follows. On tree $\mathbb{T}^d,$ call an edge $e'$ in the geodesic path connecting root $o$ and edge $e$ an ancestor edge of $e$.
Recall $L_n=\{e\in E(\mathbb{T}^d):\ \vert e\vert =n\},\ n\in\N.$ For any $i\in \N,$ let
\begin{equation*}
\begin{aligned}
&\Pi^1_i= \left\{e^*:\ e\in E(\mathbb{T}^d),\ |e|= a_{i+1}-1,\ e\ \mbox{is not connected to}\ L_{a_i}\ \mbox{in}\ \omega_p,\  e^*\  \mbox{is a closed ancestor}\right.\\
&\hskip 1.75cm \left.\mbox{of}\ e\ \mbox{with}\ \vert e^*\vert=\min\left\{\vert e'\vert:\ \vert e'\vert\geq a_i,\ e'\ \mbox{is a closed ancestor of}\ e\right\}\right\},\\
&\Pi^2_i= \left\{e\in E(\mathbb{T}^d):\ e\ \mbox{is open and connected to}\ L_{a_i}\ \mbox{in}\ \omega_p,\ |e|=a_{i+1}-1\right\},\\
&\Pi_i= \Pi_i^1\cup \Pi_i^2.
\end{aligned}
\end{equation*}
Then all $\Pi_i$s are edge-disjoint cutsets separating $o$ from $\infty$ in $\mathbb{T}^d$. Clearly
\begin{eqnarray}\label{eq-conductance-sum}
C_1:=\sum\limits_{i=1}^{\infty} \sum\limits_{e\in\Pi_i^1} \mathbf{C}_{\lambda_2}(e)\leq \sum\limits_{i=1}^{\infty} d(d-1)^{i} \lambda_2^{-i}<\infty.
\end{eqnarray}

To continue, recall the following useful facts. For the Galton-Watson tree $T$ with offspring distribution $B(d-1,p)$, write
$$m=(d-1)p\ \mbox{and}\ \sigma^2=(d-1)p(1-p).$$
Let $Z_n(T)$ be the number of individuals in the $n$th generation of $T.$ It is easy to check that
\[
\E\left[Z^2_{n+1}(T)\right]= m^2 \E\left[Z_n^2(T)\right]+ \sigma^2 m^n.
\]
Note $m>1.$ Then
$$\lim\limits_{n\rightarrow\infty}\E\left[\left(\frac{Z_n(T)}{m^n}\right)^2\right]=1+\sum\limits_{k=0}^{\infty}\frac{\sigma^2}{m^{k+2}}\in (0,\infty).$$
It is well-known that almost surely, $\left\{\frac{Z_n(T)}{m^n}\right\}$ converges to a random variable $W\in [0,\,\infty)$ with mean $1.$
Thus by the Doob's maximum inequality and the dominated convergence theorem,
 \begin{eqnarray}\label{eq-uniform-2nd-moment}
 \E\left[W^2\right]=\lim\limits_{n\rightarrow\infty}\E\left[\left(\frac{Z_n(T)}{m^n}\right)^2\right]=\sup\limits_{n} \E\left[\left(\frac{Z_n(T)}{m^n}\right)^2\right]<\infty.
 \end{eqnarray}

Notice that from {\bf (i.2)}, $x^i_{k}$ is the $k$-th $(1\leq k\leq d(d-1)^{i-1})$ individual in the $i$-th generation in $\mathbb{T}^d$; and $\mathcal{T}(x^i_{k})$ denotes random open descendant subtree rooted at $x^i_{k}$,
a Galton-Watson tree with offspring distribution $B(d-1,p).$ Write $Z_m\left(x^i_{k}\right):=Z_m\left(\mathcal{T}(x^i_{k})\right)$ for the number of individuals of the $m$-th generation of random tree $\mathcal{T}(x^i_{k})$. Then for any $i\in\N,$
\begin{equation*}
\begin{aligned}
\left|\Pi_i^2\right|=\sum\limits_{k=1}^{d(d-1)^{a_i-1}} Z_{a_{i+1}-a_{i}}(x^{a_i}_{k}).
\end{aligned}
\end{equation*}
For any $j\in\N,$ let $\mathcal{R}_j= d(d-1)^{a_j- 1}$ and
$$\mathcal{Y}_{i}^{j}=\frac{Z_{a_{j+1}-a_{j}}(x^{a_j}_{i})}{m^{a_{j+1}-a_{j}}},\, 1\leq i\leq\mathcal{R}_j.$$
By the standard theory of Galton-Watson branching processes, almost surely, for any $j\in\N$ and $1\leq i\leq\mathcal{R}_j,$
$$W_i^j=\lim\limits_{k\rightarrow\infty}\frac{Z_k(x_i^{a_j})}{m^k}\in [0,\,\infty)\ \mbox{exists};$$
and all $W_i^j$s have a common distribution as that of $W$. Let
$$\mathcal{S}_j:=\frac{\sum\limits_{i= 1}^{\mathcal{R}_j} \left(\mathcal{Y}_{i}^{j}- 1\right)}{\mathcal{R}_j}.$$
By \eqref{eq-uniform-2nd-moment},
\begin{eqnarray*}
\sup\limits_{j\in\N,\,1\leq i\leq \mathcal{R}_j}\E\left[\left(\mathcal{Y}_i^j\right)^2\right] < +\infty.
\end{eqnarray*}
Since for any $j\in\N,$ $\{\mathcal{Y}_i^j-1\}_{1\leq i\leq\mathcal{R}_j}$ is an i.i.d.\,family with $\E\left[\mathcal{Y}_i^j-1\right]=0,$ then for some constant $C_2\in (0,\,\infty),$
\begin{eqnarray*}
\E\left[\mathcal{S}_j^2\right]=\E\left[\frac{\E\left[\left(\mathcal{Y}_{1}^{j}- 1\right)^2\right]}{\mathcal{R}_j}\right]\leq \frac{C_2}{d(d-1)^{a_j-1}}.
\end{eqnarray*}
By the Chebyshev inequality, for any $\epsilon\in (0,\,\infty),$
\begin{eqnarray*}
\P[|\mathcal{S}_j|>\epsilon]\leq \frac{\E\left[\mathcal{S}_j^2\right]}{\epsilon^2}\leq \frac{C_2}{\epsilon^2 d(d-1)^{a_j-1}}
\end{eqnarray*}
which leads to
\[
\lim\limits_{j\rightarrow \infty} \mathcal{S}_j= 0\ a.s.
\]
by the Borel-Cantelli lemma. Namely, almost surely,
\[
\lim\limits_{j\rightarrow \infty} \frac{\sum\limits_{i=1}^{\mathcal{R}_j} Z_{a_{j+1}-a_{j}}(x^{a_j}_{i})}{m^{a_{j+1}-a_{j}}d(d-1)^{a_j-1}}=1.
\]
Since
\[
\sum\limits_{j=1}^{\infty} \lambda_1^{-a_{j+1}+1} m^{a_{j+1}-a_{j}}d(d-1)^{a_j-1}=\sum\limits_{j=1}^{\infty}\frac{d\lambda_1}{d-1}\left(\frac{(d-1)p}{\lambda_1}\right)^{a_{j+1}}
\left(\frac{1}{p}\right)^{a_j}<\infty,
\]
we have that almost surely
\begin{equation}\label{eq-conductance-sum-2}
\begin{aligned}
C_3:=\sum\limits_{j=1}^{\infty} \sum\limits_{e\in \Pi_j^2} \mathbf{C}_{\lambda_1}(e)&= \sum\limits_{j=1}^{\infty} \lambda_1^{-a_{j+1}+1} |\Pi_{j}^2|\leq \sum\limits_{j=1}^{\infty} \lambda_1^{-a_{j+1}+1} \sum\limits_{i=1}^{\mathcal{R}_j} Z_{a_{j+1}-a_{j}}(x^{a_j}_{i})< \infty.
\end{aligned}
\end{equation}
By \eqref{eq-conductance-sum} and \eqref{eq-conductance-sum-2}, and the Nash-Williams inequality and recurrence criterion (Lemma \ref{lem-Nash}), when $\frac{1}{d-1}< p< \frac{\lambda_1}{d-1}$, almost surely,
\begin{align*}
\mathscr{R}\left(o\leftrightarrow \infty\right)&\geq \sum\limits_{j=1}^{\infty} \left(\sum\limits_{e\in \Pi_{j}^1}\mathbf{C}_{\lambda_2}(e)+\sum\limits_{e\in\Pi_j^2}\mathbf{C}_{\lambda_1}(e)\right)^{-1}
 \geq \sum\limits_{j=1}^{\infty}(C_1+C_3)^{-1}=\infty,
\end{align*}
and $(\mathbb{T}^d,\, \mathbf{C}_{\lambda_1},\, \mathbf{C}_{\lambda_2},\, p)$ is recurrent.
\vskip 2mm

{\bf Assume secondly $\lambda_2=d-1$.} Recall the following facts. Let $T$ be a finite tree with root $o$ and leaf set $A$ such that
\begin{eqnarray}\label{eq-degree-condition}
\deg_T(o)=d-1,\ \deg_T(x)=d,\ x\in T\setminus (A\cup\{o\}).
\end{eqnarray}
Here $\deg_T(\cdot)$ is the degree function of vertices. For any $x\in T\setminus\{o\},$ let $x_*$ be the parent vertex of $x$:
$$x_*\sim x,\ \vert x_*\vert=\vert x\vert-1.$$
Define a unit flow $\theta_T$ on $T$ from $o$ to $A$: For any $x\in T\setminus\{o\},$
$$\theta_T(x_*x)=(d-1)^{-(\vert x_*\vert+1)},\ \theta_T(xx_*)=-(d-1)^{-(\vert x_*\vert+1)}.$$
By the flow conservation property (\cite[Lemma 2.8]{LP2017}),
$$\sum\limits_{x\in A}\theta_T(x_*x)=\sum\limits_{x\sim o}\theta_T(ox)=1.$$
Let $T'$ be a subtree of $T$ with the root $o$ and leaf set $A'$ that is a subset of $A.$ Then
\begin{eqnarray*}
\sum\limits_{x\in A'}\theta_{T}(x_*x)\leq \sum\limits_{x\in A}\theta_T(x_*x)=1,
\ \mbox{namely}\ \sum\limits_{x\in A'}(d-1)^{-\vert x\vert}\leq 1.
\end{eqnarray*}
Note any finite tree $T'$ with root $o$ and leaf set $A'$ satisfying
\begin{eqnarray}\label{eq-degree-condition-1}
\deg_{T'}(o)\leq d-1,\ \deg_{T'}(x)\leq d,\ x\in T'\setminus (A'\cup\{o\}),
\end{eqnarray}
can be embedded into a finite tree $T$ with root $o$ and leaf set $A$ such that both $A'\subseteq A$ and \eqref{eq-degree-condition} are true.
Therefore, for any finite tree $T'$ with root $o$ and leaf set $A'$ satisfying \eqref{eq-degree-condition-1},
\begin{eqnarray}\label{eq-leaf-mass}
\sum\limits_{x\in A'}(d-1)^{-\vert x\vert}\leq 1.
\end{eqnarray}

For any $i\in \N$ and $x\in L_{a_i},$ let $\mathcal{T}_i(x)$ be the open subtree of $\mathbb{T}^d$ consisting of $x$ and its offsprings $y$ connected to $x$ by an open path in $\omega_p$
such that $\vert y\vert \in [a_i+1,a_{i+1}-1].$ Write $T_i(x)=\mathcal{T}_i(x)\cup\partial\mathcal{T}_i(x).$
Let $A_i(x)$ be the leaf set of finite tree $T_i(x).$ Then
\begin{eqnarray*}
\sum\limits_{e\in\Pi_i^1}\mathbf{C}_{d-1}(e)&\leq &\sum\limits_{x\in L_{a_i}}\sum\limits_{e\in\partial\mathcal{T}_i(x)}\mathbf{C}_{d-1}(e)=\sum\limits_{x\in L_{a_i}}\sum\limits_{e\in\partial\mathcal{T}_i(x)}(d-1)^{-\vert
   e\vert}\\
&=& \sum\limits_{x\in L_{a_i}}\sum\limits_{y\in A_i(x)}(d-1)^{-(\vert y\vert-1)}\\
&=&(d-1)^{-a_i+1}\sum\limits_{x\in L_{a_i}}\sum\limits_{y\in A_i(x)}(d-1)^{-(\vert y\vert-a_i)},
\end{eqnarray*}
and by \eqref{eq-leaf-mass},
\begin{eqnarray}\label{eq-sum-closed}
\sum\limits_{e\in\Pi_i^1}\mathbf{C}_{d-1}(e)\leq (d-1)^{-a_i+1}\sum\limits_{x\in L_{a_i}}1=(d-1)^{-a_i+1}d(d-1)^{a_i-1}=d.
\end{eqnarray}

By \eqref{eq-conductance-sum-2} and \eqref{eq-sum-closed}, and the Nash-Williams inequality and recurrence criterion, when $\frac{1}{d-1}< p< \frac{\lambda_1}{d-1}$, almost surely,
\begin{align*}
\mathscr{R}\left(o\leftrightarrow \infty\right)&\geq \sum\limits_{j=1}^{\infty} \left(\sum\limits_{e\in \Pi_{j}^1}\mathbf{C}_{d-1}(e)+\sum\limits_{e\in\Pi_j^2}\mathbf{C}_{\lambda_1}(e)\right)^{-1}
 \geq \sum\limits_{j=1}^{\infty}(d+C_3)^{-1}=\infty,
\end{align*}
and $(\mathbb{T}^d,\, \mathbf{C}_{\lambda_1},\, \mathbf{C}_{d-1},\, p)$ is recurrent.
\qed
\vskip 2mm

\section{Concluding remarks and problems}\label{sec-concluding}
\setcounter{equation}{0}
\subsection{Biased disordered random networks $(G,\, \mathbf{C}_{\lambda_1},\, \mathbf{C}_{\lambda_2},\, p)$}
\noindent
Recall Theorem \ref{generald} and $p_c\in (0,\,1/2)$ on $\Z^d$ with $d\geq 3.$ Note that on $\Z^d$ with $d\geq 2,$ the non-existence of infinite cluster at critical percolation is a well-known conjecture, which holds for $d=2$ and $d\geq 11.$ While $p_c=1/2$ for $\Z^2,$ the competing behavior of biased $(\Z^2,\, \mathbf{C}_{\lambda_1},\, \mathbf{C}_{\lambda_2},\, 1/2)$ with $0< \lambda_1<\lambda_2=1$ is more subtle; and in this case transience seems to dominate recurrence.
These lead to the following
 \begin{conj}\label{conj-critical-recurrent/transient}
 When $3\leq d\leq 10$ and $0<\lambda_1\leq 1<\lambda_2,$  almost surely $(\Z^d,\, \mathbf{C}_{\lambda_1},\, \mathbf{C}_{\lambda_2},\, p_c)$ is recurrent. And biased $(\Z^2,\, \mathbf{C}_{\lambda_1},\, \mathbf{C}_{\lambda_2},\, 1/2)$ with $0< \lambda_1<\lambda_2=1$ is transient almost surely.
 \end{conj}
\vskip 2mm

The following conjecture on having unique currents arises naturally.
\begin{conj}\label{conj-current-unique}
For $d\geq 3$ and $0<\lambda_1\leq 1<\lambda_2,$ almost surely, all transient biased $(\Z^d,\, \mathbf{C}_{\lambda_1},\, \mathbf{C}_{\lambda_2},\, p)$ are current unique.
 \end{conj}
\vskip 2mm

For $d$-regular tree $\mathbb{T}^d$ with $d\geq 3,$ $p_c=\frac{1}{d-1}$ and $\lambda_c=d-1.$
  \begin{conj}\label{conj-critical-recurrent-tree}
For any $d\geq 3,$ biased $(\mathbb{T}^d,\, \mathbf{C}_{\lambda_1},\, \mathbf{C}_{\lambda_2},\, \frac{\lambda_1\vee 1}{d-1})$ with $1<\lambda_1<d-1\leq \lambda_2$ or
$0<\lambda_1\leq 1<d-1= \lambda_2$ is recurrent almost surely.
 \end{conj}
\vskip 2mm
Recall for any $x=(x_1,\,\ldots,\, x_d)\in\Z^d,$ $\vert x\vert=\sum\limits_{i=1}^d\vert x_i\vert.$
Note the invariance principle and the large deviation for $\RW_\lambda$ with $\lambda\in (0,1)$ on $\Z^d$ was proved in \cite{LS2020}. And from \cite{SSSWX2018}, on $\Z^d$,
$\RW_\lambda$ $(X_n)_{n=0}^{\infty}$ with $\lambda\in (0,1)$ almost surely has positive speed, i.e., $\lim\limits_{n\rightarrow\infty}\frac{\vert X_n\vert}{n}=\frac{1-\lambda}{1+\lambda};$
and the heat kernel of $(X_n)_{n=0}^{\infty}$ decays exponentially.

\begin{prob}\label{plimit}
Let $(X_n)_{n=0}^{\infty}$ be the random walk associated to any biased $(\Z^d,\, \mathbf{C}_{\lambda_1},\, \mathbf{C}_{\lambda_2},\, p)$ with $0<\lambda_1<\lambda_2<\infty.$

{\bf (i)} For transient biased $(\Z^d,\, \mathbf{C}_{\lambda_1},\, \mathbf{C}_{\lambda_2},\, p)$, does almost surely the quenched speed $\lim\limits_{n\rightarrow\infty}\frac{\vert X_n\vert}{n}$ exist? If yes, is it positive and constant almost surely?

{\bf (ii)} Prove the quenched invariance principle for random walk $(X_n)_{n=0}^{\infty}$.

{\bf (iii)} Assume $0<\lambda_1<1<\lambda_2.$ Almost surely, is there an exponential vs polynomial decay in time
phase transition for the quenched heat kernel of $(X_n)_{n=0}^{\infty}$ when $p$ varies from $0$ to $1$?
\end{prob}
\vskip 2mm

Let $K_{o}(\omega_p)$ denote the open cluster of $o$ in Bernoulli-$p$ bond percolation $\omega_p$ on graph $G.$ Define
$$p_u=\inf\{p\in [0,\,1]:\ \omega_p\ \mbox{has a unique infinite open cluster with positive probability}\}.$$
To characterize $p_c^*$ for biased disordered random networks $(G,\, \mathbf{C}_{\lambda_1},\, \mathbf{C}_{\lambda_2},\, p)$ with $G$ being a Cayley graph, we recall the following Lalley's conjecture.

\begin{conj}[Lalley \mbox{\cite[Conjecture 17]{SL2006}}]\label{conj-Lalley}
Given a transitive non-amenable graph $G$ with $p_c<p_u$ and fixed vertex $o$. If there exists a unique cluster a.s. at $p= p_u$,
then condition on $K_{o}(\omega_{p_u})$ is infinite,
\[\lim\limits_{p\uparrow p_u}\underline{\rm gr}\left(K_{o}(\omega_p)\right)= {\rm gr}(G)\ \mbox{almost surely}.\]
\end{conj}
\vskip 2mm

\begin{conj}\label{conj-p_c^*}
Let $G$ be an infinite Cayley graph with fixed vertex $o.$ Then the following hold.
\begin{enumerate}[{\bf (i)}]
\item If $G$ is amenable with $\lambda_c(G)=1,$ and $(G,\,\mathbf{C}_1)$ is transient (resp. recurrent), then $p_c^*=p_c$ for any $((G,\, \mathbf{C}_{\lambda_1},\, \mathbf{C}_{\lambda_2},\, p))_{p\in [0,\,1]}$ with $0<\lambda_1\leq 1<\lambda_2$ (resp. $0<\lambda_1< 1\leq \lambda_2 $).

\item When $G$ is amenable with $\lambda_c(G)>1,$ and $\left(G,\,\mathbf{C}_{\lambda_c(G)}\right)$ is transient (resp. recurrent),
the threshold $p_c^*$ of any $((G,\, \mathbf{C}_{\lambda_1},\, \mathbf{C}_{\lambda_2},\, p))_{p\in [0,\,1]}$ with $0<\lambda_1\leq\lambda_c(G)<\lambda_2$ (resp. $0<\lambda_1<\lambda_c(G)\leq\lambda_2$) satisfies that
\begin{eqnarray*}
p_c^*=p_c^*(\lambda_1)=\inf\left\{p\in (p_c,1]:\ \mathbb{P}\left[\left. \underline{\rm gr}(K_o(\omega_p))>(\lambda_1\vee 1)\,\right\vert \vert K_o(\omega_p)\vert=\infty\right]>0\right\}
\end{eqnarray*}
with convention $\inf\emptyset=p_c,$ and $p_c^*(\lambda_1)\in [p_c,\,1)$ is continuous in $\lambda_1\in(0,\,\lambda_c(G))$ and strictly increasing in $\lambda_1\in [1,\,\lambda_c(G)),$ and
$\lim\limits_{\lambda_1\uparrow\lambda_c(G)}p_c^*(\lambda_1)=1.$

\item For any non-amenable $G$ with $p_c<p_u$ such that $\omega_{p_u}$ has a.s.\,a unique infinite cluster and $\left(G,\,\mathbf{C}_{\lambda_c(G)}\right)$ is transient (resp. recurrent),
any $((G,\,\mathbf{C}_{\lambda_1},\, \mathbf{C}_{\lambda_2},\, p))_{p\in [0,\,1]}$
with $0<\lambda_1\leq\lambda_c(G)<\lambda_2$ (resp. $0<\lambda_1<\lambda_c(G)\leq\lambda_2$) has the threshold
  \begin{eqnarray*}
p_c^*=p_c^*(\lambda_1)=\inf\left\{p\in (p_c,1]:\ \mathbb{P}\left[\left. \underline{\rm gr}(K_o(\omega_p))>(\lambda_1\vee 1)\,\right\vert \vert K_o(\omega_p)\vert=\infty\right]>0\right\}\in [p_c,\,p_u),
\end{eqnarray*}
and $p_c^*(\lambda_1)$ is continuous in $\lambda_1\in(0,\,\lambda_c(G))$ and strictly increasing in $\lambda_1\in [1,\,\lambda_c(G)),$ and
$\lim\limits_{\lambda_1\uparrow\lambda_c(G)}p_c^*(\lambda_1)=p_u.$

\item Assume $G$ is non-amenable such that $p_c<p_u$ and $\omega_{p_u}$ has a.s.\, infinitely many infinite clusters, and $\left(G,\,\mathbf{C}_{\lambda_c(G)}\right)$ is transient (resp. recurrent).
Then for any $((G,\, \mathbf{C}_{\lambda_1},\, \mathbf{C}_{\lambda_2},\, p))_{p\in [0,\,1]}$ with $0<\lambda_1\leq \lambda_c(G)<\lambda_2$ (resp. $0<\lambda_1<\lambda_c(G)\leq\lambda_2$),
  \begin{eqnarray*}
p_c^*=p_c^*(\lambda_1)=\inf\left\{p\in (p_c,1]:\ \mathbb{P}\left[\left. \underline{\rm gr}(K_o(\omega_p))>(\lambda_1\vee 1)\,\right\vert \vert K_o(\omega_p)\vert=\infty\right]>0\right\}\in [p_c,\,p_u],
\end{eqnarray*}
and $p_c^*(\lambda_1)$ is continuous in $\lambda_1\in (0,\,\lambda_c(G))$, and there is a $\lambda_1^*\in (1,\lambda_c(G))$ such that $p_c^*(\lambda_1)=p_u$ for $\lambda_1\in [\lambda_1^*,\,\lambda_c(G))$
and $p_c^*(\lambda_1)$ is strictly increasing in $\lambda_1\in [1,\,\lambda_1^*].$
\end{enumerate}
\end{conj}
\vskip 2mm

For $(G,\, \mathbf{C}_{\lambda_1},\, \mathbf{C}_{\lambda_2},\, p)$, let $\theta(p)$ be the probability of event $\{X_0=o,\, X_n\not=o,\,\forall n\geq 1\}$ given $\omega_p.$ Note $\theta(p)$ is a measurable function of $\omega_p$ and thus of whole process $(\omega_q)_{q\in [0,1]}.$ So if $G$ is quasi-transitive, then by the ergodic property of $(\omega_q)_{q\in [0,1]},$ almost surely, $\theta(p)$ is a constant for any $p\in [0,1].$

\begin{prob}\label{prob-exit-probability}
Assume $G$ is quasi-transitive with $p_c\in (0,1).$ Is $\theta(p_c^*)$ zero? Is $\theta(p)$ continuous in $p\in [p_c^*,1]?$ Is $\theta(p)$ strictly increasing in $p\in [p_c^*,1]?$ And if $\theta(p)$ is right continuous at $p_c^*\in (0,1),$ is there a critical exponent $\alpha>0$ such that
$\theta(p)-\theta(p_c^*)\asymp \vert p-p_c^*\vert ^{\alpha+o(1)}$ as $p\downarrow p_c^*?$
\end{prob}


\subsection{Other disordered random networks}
\noindent
\begin{prob}\label{pphasetransition}
Given a quasi-transitive infinite graph $G$, and any two conductance functions $\mathbf{c}_1$ and $\mathbf{c}_2$ such that $\mathbf{c}_1(e)\geq \mathbf{c}_2(e)$ for any $e$ of $G$, $(G,\, \mathbf{c}_1)$ is transient and $(G,\, \mathbf{c}_2)$ is recurrent.
Almost surely, when is there a nontrivial recurrence/transience phase transition for $(G,\, \mathbf{c}_1,\, \mathbf{c}_2,\, p)$ as $p$ varies from $0$ to $1$? Does there exist a class of graphs $G$ such that almost surely there is always a nontrivial recurrence/transience phase transition of $(G,\, \mathbf{c}_1,\, \mathbf{c}_2,\, p)$ for any fixed $\mathbf{c}_1$ and $\mathbf{c}_2$?
\end{prob}
\vskip 2mm

\begin{prob}\label{prob-PhaseTransition}
Let $G$ be a quasi-transitive infinite graph. {\bf (i)} Assume $(G,\,\mathbf{c}_1)$ is not current unique while $(G,\,\mathbf{c}_2)$ is. Study current uniqueness/nonuniqueness phase transition for $(G,\,\mathbf{c}_1,\,\mathbf{c}_2,\,p)$ when $p$ varies from $0$ to $1.$

{\bf (ii)} Suppose the free (wired) Ising model on $G=(V,\,E)$ with coupling constants $(\mathbf{c}_1(e))_{e\in E}$ (resp. $(\mathbf{c}_2(e))_{e\in E}$) and inverse temperature $1$ is in a ferromagnet (resp. paramagnet) regime. Investigate paramagnet/ferromagnet
phase transition for the free (wired) Ising model on $G$ with coupling constants $\left(\mathbf{c}_1(e)I_{\{e\in\omega_p\}}+\mathbf{c}_2(e)I_{\{e\notin\omega_p\}}\right)_{e\in E}$ and inverse temperature $1$ when $p$ varies from $0$ to $1.$
\end{prob}
\vskip 2mm

\begin{prob}\label{prob-synchronize}
Given a quasi-transitive infinite graph $G$. {\bf (i)} When $(G,\, \mathbf{c}_1)$ and $(G,\, \mathbf{c}_2)$ are both recurrent (resp. transient) networks, is $(G,\, \mathbf{c}_1,\, \mathbf{c}_2,\, p)$ recurrent (resp. transient) almost surely for $p\in [0,\, 1]$? {\bf (ii)} If both $(G,\, \mathbf{c}_1)$ and $(G,\, \mathbf{c}_2)$ have current uniqueness, does $(G,\, \mathbf{c}_1,\, \mathbf{c}_2,\, p)$ have current uniqueness almost surely for $p\in [0,\, 1]$?
\end{prob}
\vskip 2mm

\begin{remark}[Discussion on Problem \ref{prob-synchronize} (i)] Assume $(G,\, \mathbf{c}_1)$ and $(G,\, \mathbf{c}_2)$ are both recurrent, and one of the following conditions is true: {\bf (a)} $(G,\, \mathbf{c}_1\wedge \mathbf{c}_2)$ is recurrent or $\mathbf{c}_1\asymp \mathbf{c}_2$; {\bf (b)} $(G,\, \mathbf{c}_1)$ and $(G,\, \mathbf{c}_2)$ are both positive recurrent. Then $(G,\, \mathbf{c}_1,\, \mathbf{c}_2,\, p)$ is recurrent almost surely for $p\in [0,\, 1]$. It is unknown whether this is true generally.

There are examples that $(G,\, \mathbf{c}_1)$ and $(G,\, \mathbf{c}_2)$ are both transient, but $(G,\, \mathbf{c}_1,\, \mathbf{c}_2,\, p)$ is recurrent almost surely for $p\in (0,\, 1)$.
On $\Z$, define for $0< \lambda< 1$,
$$\mathbf{c}_1(\{i,i-1\}) =\lambda^{i}I_{\{i\leq 0\}}+I_{\{i>0\}},\ \mathbf{c}_2(\{i,i+1\}) =\lambda^{-i}I_{\{i\geq 0\}}+I_{\{i<0\}},\ i\in\Z.$$
Then $(\Z,\, \mathbf{c}_{1})$ and $(\Z,\, \mathbf{c}_{2})$ are transient; and by the Nash-Williams criterion, $(\Z,\, \mathbf{c}_1,\, \mathbf{c}_2,\, p)$ is a.s. recurrent for $p\in(0,\, 1).$
To obtain an example on $\Z^2$, define for $0< \lambda< 1$, $\mathbf{c}_1(e)=\lambda^{-\vert e\vert}$ when $e$ is in the negative $x$-axis and $1$ otherwise; and similarly $\mathbf{c}_2(e)=\lambda^{-\vert e\vert}$ when $e$ is in the positive $x$-axis and $1$ otherwise. Then $(\Z^2,\, \mathbf{c}_{1})$ and $(\Z^2,\, \mathbf{c}_{2})$ are transient. However, $(\Z^2,\,\mathbf{c}_1,\, \mathbf{c}_2,\, p)$ is a.s. recurrent for $p\in (0,\, 1)$,
which can be proved by the Nash-Williams criterion,
\[
\mathscr{R}(0\leftrightarrow \infty)\geq \sum\limits_{k=1}^{\infty}\left(\sum\limits_{e\in \Pi_{k}} c(e)\right)^{-1}= \infty,
\]
where each $\Pi_{k}$ is taken as follows: Let $\{e_{k,+}\}_{k=1}^{\infty}$ (resp. $\{e_{k,-}\}_{k=1}^{\infty}$) be all open (resp. closed) edges in the positive (resp. negative) axis such that $\vert e_{k,+}\vert$ (resp.
$\vert e_{k,-}\vert$) is strictly increasing in $k$. Let $\Pi_k$ be the set of all edges in $[-\vert e_{k,-}\vert-1,\vert e_{k,+}\vert+1]^2$ with exactly one endpoint in $[-\vert e_{k,-}\vert,\vert e_{k,+}\vert]^2.$
To apply the Nash-Williams criterion, one needs to note that almost surely,
$$\limsup\limits_{k\rightarrow\infty}\frac{\vert e_{k,+}\vert}{k}<\infty\ \mbox{and}\ \limsup\limits_{k\rightarrow\infty}\frac{\vert e_{k,-}\vert}{k}<\infty.$$
\end{remark}
\vskip 2mm

Note for $(G,\, \mathbf{c}_1,\, \mathbf{c}_2,\, p)$, the random environment is given by Bernoulli-$p$ bond percolation. Instead of Bernoulli bond percolation on graph $G=(V,\,E),$ one can define disordered random networks with random environments described by point processes on $E$ such as Poisson point process (PPP), determinant point process (DPP) and random cluster model through an obvious ways. To study phase transition, one needs to let the edge density in these environments vary from $0$ to $1$ via some natural means.

Additionally, one can also define disordered random network in random environment provided respectively by PPP, DPP and Bernoulli-$p$ site percolation on $V$ as follows: Let random subset $V_1$ of $V$ follow one of the distributions of PPP, DPP and Bernoulli-$p$ site percolation on $V.$ Write $E(V_1)$ (resp. $E_o(V_1)$) for the set of edges of $G$ whose one endpoint intersects (resp. two endpoints intersect) $V_1$. By letting edges in $E(V_1)$ (resp. $E_o(V_1)$) take conductance $\mathbf{c}_1$ and other edges take conductance $\mathbf{c}_2,$ one gets the desired disordered random network.

Finally, let $\mathbf{c}_{V,1}$ and $\mathbf{c}_{V,2}$ be two nonnegative weight functions on $V$ such that $(G,\,\mathbf{c}_1)$ is transient and $(G,\,\mathbf{c}_2)$ is recurrent, where
$$\mathbf{c}_1(\{x,y\})=\mathbf{c}_{V,1}(x)\mathbf{c}_{V,1}(y)\ \mbox{and}\ \mathbf{c}_2(\{x,y\})=\mathbf{c}_{V,2}(x)\mathbf{c}_{V,2}(y)\
  \mbox{for any edge}\ \{x,y\}\ \mbox{of}\ G.$$
For the just mentioned random subset $V_1$, let $V_2=V\setminus V_1$ and
$$\mathbf{c}_{1,2}(\{x,y\})=\mathbf{c}_{V,i}(x)\mathbf{c}_{V,j}(y)\ \mbox{if}\ x\in V_i,\, y\in V_j\ \mbox{and}\ \{x,y\}\in E.$$
Thus one obtain a disordered random network $(G,\,\mathbf{c}_{1,2}).$

It is very interesting to study various typical properties for the above disordered random networks. Recall from \cite{MA1960, FM2019} that there is an interesting random resistor network built in a different way in a homogeneous Poisson point process environment, which is called Miller-Abrahams random resistor network.

\subsection{Disordered random walks}\label{sec-disorderedRW}
\noindent Let $\mathbf{p}_1=\mathbf{p}_1(\cdot,\,\cdot)$ and $\mathbf{p}_2=\mathbf{p}_2(\cdot,\,\cdot)$ be respectively two 1-step transition probabilities of two Markov chains on graph $G=(V,\,E).$ Let $V_1=V_1(p)$ be the set of open vertices in Bernoulli-$p$ site percolation on $G$ and $V_2=V\setminus V_1.$ Use $(G,\,\mathbf{p}_1,\,\mathbf{p}_2,\,p)$ to denote the Markov chains on $G$ with $1$-step transition probability given by
$$\mathbf{p}(x,\,y)=\mathbf{p}_1(x,\,y)I_{\{x\in V_1\}}+\mathbf{p}_2(x,\,y)I_{\{y\in V_2\}},\ x,y\in V.$$
Like the disordered random network $(G,\,\mathbf{c}_1,\,\mathbf{c}_2,\, p)$, we also define Bernoulli site percolation process $(V_1(p))_{p\in [0,\,1]}$ on $G$ by the grand coupling, and also study recurrence/transience phase transition (and other ones) for
disordered random walk $(G,\,\mathbf{p}_1,\,\mathbf{p}_2,\,p)$ as $p$ varies from $0$ to $1.$

Notice that recurrence/transience phase transition for biased $(G,\,\mathbf{p}_1,\,\mathbf{p}_2,\,p)$ is more subtle than that of biased $(G,\,\mathbf{c}_1,\,\mathbf{c}_2,\, p)$ due to that unlike the latter, there is no the Rayleigh's monotonicity principle for the former. For example, similarly to Remark \ref{remark-PT}\,(ii), one can also define $p_c^*$ and $\widehat{p}_c^*$ for biased $(G,\,\mathbf{p}_1,\,\mathbf{p}_2,\,p)$, but $p_c^*=\widehat{p}_c^*$ may not hold and transient regime may not be $(p_c^*,1]$ or $(\widehat{p}_c^*,1].$ The disordered random walks $(G,\, \mathbf{p}_1,\, \mathbf{p}_2,\, p)$ are a new kind of RWREs. Replacing Bernoulli site percolation by a PPP or DPP on $V$, one gets another disordered random walk.

\begin{prob}
Suppose $G$ is a quasi-transitive infinite graph, $\mathbf{p}_1$-random walk is transient and $\mathbf{p}_2$-random walk is recurrent. Characterize the recurrent and transient regimes for the disordered random walk family $\left((G,\, \mathbf{p}_1,\, \mathbf{p}_2,\, p)\right)_{p\in [0,1]}$. When are these regimes almost surely a single subinterval of $[0,\,1]$ such that almost surely, there is a threshold $p_c^*$ (may be random) satisfying $(G,\, \mathbf{p}_1,\, \mathbf{p}_2,\, p)$ is recurrent for any $p<p_c^*$ and transient for any $p>p_c^*?$ Study various typical properties for the disordered random walks (particularly the biased ones).
\end{prob}
\vskip 2mm


\begin{thm}\label{thm-sz1}
Let $0<\lambda_1<\lambda_2<\infty$ and each $\mathbf{p}_i$-random walk be ${\RW}_{\lambda_i}$ on $\Z,$ and
$$p_c^*=\left(\frac{\log \lambda_2}{\log \lambda_2- \log \lambda_1}\vee 0\right) \wedge 1.$$
Then almost surely, $(\Z,\, \mathbf{p}_1,\, \mathbf{p}_2,\, p)$ is recurrent for any $p\in [0,p_c^*]$ and transient for any
$p\in (p_c^*,1].$
\end{thm}
\vskip 2mm

For the recurrence/transience phase transition of biased $\left((\Z,\, \mathbf{p}_1,\, \mathbf{p}_2,\, p)\right)_{p\in [0,1]}$ with each $\mathbf{p}_i$-random walk being ${\RW}_{\lambda_i}$ and $0<\lambda_1<1\leq\lambda_2<\infty,$ $p_c^*=\frac{\log \lambda_2}{\log \lambda_2- \log \lambda_1}$, which is not related to the phase transition of the Bernoulli site percolation on $\Z.$ And hence it has a novel nature
different from that of biased random networks $\left((\Z,\, \mathbf{c}_1,\, \mathbf{c}_2,\, p)\right)_{p\in [0,1]}.$

Note that on regular trees $\mathbb{T}^d$ with $d\geq 3,$ the biased disordered random walk has i.i.d.\,random transition probabilities on all vertices but the root. Recall from \cite{LP1992} and \cite{RY1995}, Lyons, Pemantle and Peres gave complete recurrent/transience criteria for RWREs with some special i.i.d.\,random environments on trees in 1990s. By these criteria and some monotonicity similarly to \eqref{stoc-monotone}, we have the following clear picture of phase transition for the biased disordered random walks on $\mathbb{T}^d$:
\begin{enumerate}
\item[] On $\mathbb{T}^d$ with $d\geq 3$, let each $\mathbf{p}_i$-random walk be ${\RW}_{\lambda_i}$ with $0<\lambda_1<\lambda_2<\infty.$
Then when $\lambda_1\geq d-1$ (resp. $\lambda_2<d-1$), almost surely, $(\mathbb{T}^d,\, \mathbf{p}_1,\, \mathbf{p}_2,\, p)$ is recurrent (resp. transient) for all $p\in [0,\, 1];$ when $\lambda_2=d-1,$ almost surely, $(\mathbb{T}^d,\, \mathbf{p}_1,\, \mathbf{p}_2,\, p)$ is transient for all $p\in (0,\, 1]$ (note $(\mathbb{T}^d,\, \mathbf{p}_1,\, \mathbf{p}_2,\, 0)$ is recurrent). And when $\lambda_1<d-1<\lambda_2$, there is a non-trivial recurrence/transience phase transition for the biased disordered random walks $(\mathbb{T}^d,\, \mathbf{p}_1,\, \mathbf{p}_2,\, p)$ such that for a constant $p^{*}_c\in (0,\, 1)$, almost surely, $(\mathbb{T}^d,\, \mathbf{p}_1,\, \mathbf{p}_2,\, p)$ is recurrent for any $p\leq p_c^{*}$ and transient for $p> p_c^{*};$ where $p_c^*$ is the unique solution to
$$f(p)=\min\limits_{0\leq x\leq 1}\left\{\left(\frac{1}{\lambda_1}\right)^xp+\left(\frac{1}{\lambda_2}\right)^x(1-p)\right\}=\frac{1}{d-1},\ p\in [0,1];$$
and explicitly $p_c^*=\frac{\frac{1}{d-1}-\frac{1}{\lambda_2}}{\frac{1}{\lambda_1}-\frac{1}{\lambda_2}}$ when $1\leq\lambda_1<d-1<\lambda_2.$
\end{enumerate}

\noindent{\bf Proof of Theorem \ref{thm-sz1}.}
{\bf Step 1.} Recall some preliminaries on random walks.
%

\begin{lem}[\cite{FS1975}]\label{Cprop}
Let $\{X_n\}_{n=0}^{\infty}$ be a sequence of i.i.d.\,non-degenerate finite-value random variables and $S_n= \sum\limits_{i=1}^{n} X_i,\, n\in \N$. Then {\bf (i)} $\sum\limits_{n=1}^{\infty} n^{-1}\P\left(S_n >0\right)< \infty \Longleftrightarrow\lim\limits_{n\rightarrow \infty} S_n= -\infty \ \mbox{a.s.},$ and under this condition $\sum\limits_{n=1}^{\infty} e^{S_n}< \infty$ a.s.; and {\bf (ii)}
$\sum\limits_{n=1}^{\infty} n^{-1}\P\left(S_n >0\right)= \infty=\sum\limits_{n=1}^{\infty} n^{-1}\P\left(S_n <0\right)$
is equivalent to
\[
 -\infty= \liminf\limits_{n\rightarrow \infty} S_n< \limsup\limits_{n\rightarrow \infty} S_n= \infty\ a.s.,
\]
and in this case $\sum\limits_{n=1}^{\infty} e^{S_n}= \infty= \sum_{n=1}^{\infty} e^{-S_n}$ a.s..
\end{lem}

For any countably infinite set $\Sigma$, say a sequence $(x_n)_{n\geq 1}\subseteq \Sigma$ converges to $\infty$ if for any finite subset $A$ of $\Sigma,$ $x_n\notin A$ for large enough $n.$

\begin{lem}[Lyapunov recurrence/transience criterion \mbox{(\cite[Chapter~2.5]{MSA2015}, \cite{FMS1998})}]\label{Lcriterion} Suppose $\{X_n\}_n$ is an irreducible Markov chain on a countably infinite state space $\Sigma$. Fix $o\in\Sigma.$
\begin{enumerate}[{\bf (i)}]
\item $\{X_n\}_n$ is recurrent if and only if there are a function $f:\, \Sigma \rightarrow \R_{+}$ and a finite non-empty set $A$ of $\Sigma$ such that $f(x)\rightarrow \infty$ as $x\rightarrow \infty$, and
\[
\E[f(X_{n+1})- f(X_n)  \, \vert\,  X_n= x]\leq 0\ \mbox{for all}\ x\in \Sigma \backslash A\ \mbox{and}\ n\geq 0.
\]
$\{X_n\}_n$ is transient if and only if there exist a function $f:\, \Sigma \rightarrow \R_{+}$ and a non-empty set $A\subset\Sigma$
satisfying $f(y)<\inf\limits_{x\in A} f(x)$ for at least one $y\in \Sigma\backslash A,$ and
\[
\E[f(X_{n+1})- f(X_n)  \ \vert\  X_n= x]\leq 0\ \mbox{for all}\ x\in \Sigma \backslash A \ \mbox{and}\ n\geq 0.
\]

\item $\{X_n\}_n$ is recurrent if there is a function $f:\, \Sigma\rightarrow \R$ with that $f(x)\rightarrow\infty$ as $x\rightarrow\infty$ and
\[\E\left[f(X_{n+1})- f(X_n)  \, \vert\,  X_n=x\right]\leq 0,\ x\not=o,\ n\geq 0;\]
and is transient if there is a bounded and non-constant function $f:\, \Sigma \rightarrow \mathbb{R}$ such that
\[\E\left[f(X_{n+1})- f(X_n)  \, \vert\,  X_n=x\right]= 0,\ x\not=o,\ n\geq 0.\]
\end{enumerate}
\end{lem}

\noindent{{\bf Step 2. Turn to prove Theorem \ref{thm-sz1}.}

The proof is routine. Let $\omega_p=\Z(p)$ be the Bernoulli-$p$ site percolation on $\Z$, and $\omega=(\omega_p)_{0\leq p\leq 1}$ the Bernoulli site percolation process constructed by the grand coupling. Fix a random environment $\omega$, for $x>0$, let
\begin{eqnarray*}
&&p_x(p)=\mathbf{p}_1(x,x-1)I_{\{\omega_p(x)=1\}}+\mathbf{p}_2(x,x-1)I_{\{\omega_p(x)=0\}}
     =\frac{\lambda_1}{1+\lambda_1}I_{\{\omega_p(x)=1\}}+\frac{\lambda_2}{1+\lambda_2}I_{\{\omega_p(x)=0\}},\\
&&q_x(p)=\mathbf{p}_1(x,x+1)I_{\{\omega_p(x)=1\}}+\mathbf{p}_2(x,x+1)I_{\{\omega_p(x)=0\}}
     =\frac{1}{1+\lambda_1}I_{\{\omega_p(x)=1\}}+\frac{1}{1+\lambda_2}I_{\{\omega_p(x)=0\}};
\end{eqnarray*}
and for $x<0,$ let
\begin{eqnarray*}
&&p_x(p)=\mathbf{p}_1(x,x+1)I_{\{\omega_p(x)=1\}}+\mathbf{p}_2(x,x+1)I_{\{\omega_p(x)=0\}}
     =\frac{\lambda_1}{1+\lambda_1}I_{\{\omega_p(x)=1\}}+\frac{\lambda_2}{1+\lambda_2}I_{\{\omega_p(x)=0\}},\\
&&q_x(p)=\mathbf{p}_1(x,x-1)I_{\{\omega_p(x)=1\}}+\mathbf{p}_2(x,x-1)I_{\{\omega_p(x)=0\}}
     =\frac{1}{1+\lambda_1}I_{\{\omega_p(x)=1\}}+\frac{1}{1+\lambda_2}I_{\{\omega_p(x)=0\}}.
\end{eqnarray*}
To construct a nonconstant bounded function $f=f_{\omega,p}:\ \Z\rightarrow\R$ which is harmonic on $\Z\backslash\{0\}$ for $(\Z,\,\mathbf{p}_1,\,\mathbf{p}_2,\,p).$ Define $f(0)= 0$ and $f(1)= f(-1)= 1$. Since for any $x> 0$,
$$\left(f(x-1)- f(x)\right)p_x(p)+  \left(f(x+1)- f(x)\right)q_x(p)= 0,$$
we have that
\[
\frac{f(x+1)- f(x)}{f(x)- f(x-1)}= \frac{p_x(p)}{q_x(p)}\ \mbox{and}\
f(x+1)= 1+ \sum\limits_{i=1}^{x} e^{\sum_{k=1}^i Y_k.}
\]
Here $\{Y_k\}_{k\in \Z\setminus\{0\}}$ is a random sequence with
$$Y_k:=Y_k(p)=\log(p_k(p)/q_k(p))=\log\left(\lambda_1I_{\{\omega_p(k)=1\}}+\lambda_2I_{\{\omega_p(k)=0\}}\right).$$
Similarly for any $x< 0$,
$$f(x-1)= 1+ \sum_{i=1}^{x} e^{\sum_{k=1}^i Y_{-k}}.$$
Clearly $\{p_x(p)/q_x(p)\}_{x\in\Z\setminus\{0\}}$ is a sequence of i.i.d.\,random variables.

By Lemma \ref{Cprop}\,(i),
$f$ is almost surely a non-constant bounded function on $\Z$ if
 \[\E\left(Y_x\right)=p\cdot \log \lambda_1+ (1-p) \cdot \log \lambda_2 < 0,\ x\in\Z\setminus\{0\},\ \mbox{namely}\ p>p_c^*.\]
Note that
\begin{eqnarray}\label{stoc-monotone}
\mbox{each}\ Y_k(p)\ \mbox{is a decreasing function in}\ p\in [0,\,1]\ \mbox{for any environment}\ \omega,\ \mbox{so is every}\ f(k).
\end{eqnarray}
Then by \eqref{stoc-monotone} and Lemma \ref{Lcriterion}\,(ii), almost surely, $(\Z,\, \mathbf{p}_1,\, \mathbf{p}_2,\, p)$ is transient for all $p\in \left(p_c^*,\, 1\right].$ When each $\E\left(Y_x \right)>0$, by the law of large numbers, $\lim\limits_{x\rightarrow\infty}f(x)=\lim\limits_{x\rightarrow -\infty}f(x)=\infty$ a.s.. Thus by \eqref{stoc-monotone} and Lemma \ref{Lcriterion}\,(ii) again, almost surely, $(\Z,\, \mathbf{p}_1,\, \mathbf{p}_2,\, p)$ is recurrent for all $p\in \left[0,\,p_c^*\right).$

When $p_c^*= \frac{\log \lambda_2}{\log \lambda_2- \log \lambda_1}\in (0,\, 1)$, for any $n\in\mathbb{N},$ define
$$Z_n=\frac{Y_n(p_c^*)- \log \lambda_1}{\log \lambda_2- \log \lambda_1},\ S_n=\sum_{i=1}^{n}Y_i(p_c^*),\ S'_n=\sum_{i=1}^{n} Z_i.$$
Then $S'_n$ has binomial distribution $B(n,\, 1-p_c^*)$, and
\begin{equation*}
\begin{split}
\sum\limits_{n=1}^{\infty} n^{-1} \P\left(S_n > 0\right)&=\sum\limits_{n=1}^{\infty} n^{-1} \P\left(S'_n > \frac{-n \log \lambda_1}{\log \lambda_2- \log \lambda_1}\right)= \sum\limits_{n=1}^{\infty} n^{-1} \P\left(S'_n > n(1-p_c^*)\right)=\infty,\\
\sum\limits_{n=1}^{\infty} n^{-1} \P\left(S_n < 0\right)&= \sum\limits_{n=1}^{\infty} n^{-1} \P\left(S'_n< n(1-p_c^*)\right)=\infty.
\end{split}
\end{equation*}
Here we have used that by the central limit theorem,
$$\lim\limits_{n\rightarrow\infty}\P\left(S'_n > n(1-p_c^*)\right)=\lim\limits_{n\rightarrow\infty}\P\left(S'_n < n(1-p_c^*)\right)=\frac{1}{2}.$$
By Lemma \ref{Cprop}\,(ii), almost surely, $f_{\omega,p_c^*}(x+1)= 1+ \sum_{i=1}^{x} e^{S_i}\rightarrow \infty,\ x\rightarrow\infty.$
Similarly,
$$\mbox{almost surely},\ f_{\omega,p_c^*}(x-1)\rightarrow \infty,\ x\rightarrow-\infty.$$
Thus $\left(\Z,\, \mathbf{p}_1,\, \mathbf{p}_2,\, p_c^*\right)$ is recurrent almost surely by Lemma \ref{Lcriterion}\,(ii). The proof is done.\qed

\subsection{Discussion on volume growth rate of percolation clusters}\label{sec-discussion-growth}
\noindent
On a connected transitive graph $G$, the number of infinite open clusters is constant almost surely which can only be $0$, $1$, or $\infty$ (\cite{NS1981}). And the constant cannot be $\infty$ when $G$ is amenable (\cite{BK1989, GKN1992}). On non-amenable graphs, there might exists two phase transitions, namely there exist $0<p_c<p_u\leq 1$ such that there are infinitely many infinite open clusters when $p\in (p_c, p_u)$. Benjamini and Schramm  \cite{BS1996} conjectured that $p_u< 1$ when $G$ is quasi-transitive with one end and $p_c< p_u$ when $G$ is quasi-transitive and non-amenable. For recent progresses on `$p_c< p_u$'-conjecture, see Hutchcroft \cite{HT2019, HT2020}.
\vskip 2mm

Here we can verify Conjecture \ref{conj-Lalley} on transitive non-amenable graphs $G$ with $0< p_c< p_u=1$.

Let $\mathcal{T}$ be the geodesic spanning tree of $G$ in Lemma \ref{brgr}. Note the basic property of trees indicates that the lower growth rate is no less than its branching number. Then by \eqref{ppercolation} and \cite[Theorem~5.15]{LP2017}, almost surely, when $p_c< p\rightarrow p_u=1$,
\[
\sup\limits_{\sigma \in \mathcal{T}} \underline{\rm gr}\left(K_{\sigma}(\omega_{p,\mathcal{T}})\right)\geq \sup\limits_{\sigma \in \mathcal{T}} {\rm br}\left(K_{\sigma}(\omega_{p,\mathcal{T}})\right)= p\cdot {\rm br}(\mathcal{T})\rightarrow
  {\rm br}(\mathcal{T})={\rm gr}(G),
\]
\[\sup\limits_{\sigma \in G} \underline{\rm gr}\left(K_{\sigma}(\omega_p)\right) \geq \sup\limits_{\sigma \in G} \underline{\rm gr}\left(K_{\sigma}(\omega_{p,\mathcal{T}})\right)\rightarrow K \geq {\rm br}(\mathcal{T});\]
where $\omega_{p,\mathcal{T}}$ is the restriction of Bernoulli percolation $\omega_p$ to $\mathcal{T}.$ Together with $\sup\limits_{\sigma \in G} \underline{\rm gr}\left(K_{\sigma}(\omega_p)\right)\leq {\rm gr}(G)$, we obtain that almost surely, as $p_c< p\rightarrow 1$,
\begin{eqnarray}\label{rate-limit}
\sup\limits_{\sigma \in G} \underline{\rm gr}\left(K_{\sigma}(\omega_{p})\right)\rightarrow {\rm br}(\mathcal{T})={\rm gr}(G).
\end{eqnarray}

Recall that a set of subgraphs of $G$ is automorphism-invariant if and only if it is invariant under action of any element of ${\rm Aut}(G)$; and a percolation measure $\P$ has indistinguishable infinite clusters if and only if for any automorphism-invariant property $\mathcal{A}$ of $subgraphs$, $\P$ a.s. either all infinite clusters satisfy $\mathcal{A}$, or all infinite clusters don't satisfy $\mathcal{A}$. Note that Lyons and Schramm \cite{RO1999}
proved that on any quasi-transitive unimodular graph, an insertion-tolerant invariant bond percolation has indistinguishable infinite clusters.

Since the (lower) volume growth rate is an automorphism-invariant property and Cayley graph is a class of transitive and unimodular graphs, we have that $\P$ almost surely, all open infinite clusters have the same lower volume growth rate. This together with \eqref{rate-limit} imply that Conjecture \ref{conj-Lalley} holds in the mentioned case.

\end{document}